\documentclass[a4paper,11pt]{amsart}
\usepackage{epsf,psfig,amscd,amsmath,amssymb,amsxtra,secthm}
1
\def\str{\operatorname{str}}
\def\({\left (}
\def\){\right )}

\def\rb{{\includegraphics{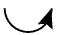}}}
\def\lb{{\includegraphics{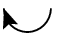}}}
\def\rd{{\includegraphics{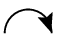}}}
\def\ld{{\includegraphics{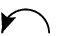}}}
\def\RT{\protect\operatorname{\cal R\cal T}\protect}
\def\A{\protect\operatorname{\cal A}\protect}
\def\CC{\protect\operatorname{\cal C}\protect}
\def\lk{\operatorname{lk}}
\def\fig#1{\vcenter{\hbox{\includegraphics{fig/#1}}}}
\def\UD{\underline\Delta}
\newcommand{\GD}{\Delta}
\newcommand{\Gd}{\delta}
\newcommand{\Ga}{\alpha}
\newcommand{\Gb}{\beta}
\newcommand{\Gt}{\tau}
\newcommand{\Ge}{\epsilon}
\newcommand{\GG}{\Gamma}
\newcommand{\GF}{\Phi}
\newcommand{\GS}{\Sigma}
\newcommand{\cal}{\mathcal}
\newcommand{\p}{\partial}
\newcommand{\empt}{\varnothing}
\newcommand{\sminus}{\smallsetminus}
\newcommand{\Z}{\mathbb Z}

\newcommand{\C}{\mathbb C}
\newcommand{\R}{\mathbb R}
\newcommand{\Int}{\operatorname{Int}}
\newcommand{\Ker}{\operatorname{Ker}}
\newcommand{\Hom}{\operatorname{Hom}}
\newcommand{\id}{\operatorname{id}}
\newcommand{\tr}{\operatorname{tr}}
\def\GT{\mathrm{T}}
\def\GL{\Lambda}
\def\im{\boldsymbol{i}}
\def\placeY#1#2#3{${}^{#1}\raisebox{-3ex}{\includegraphics{fig/#2}}\ 
{}^{#3}$} \def\placeA#1#2#3{${}_{#1}{\includegraphics{fig/#2}}\ {}_{#3}$}
\def\[{\left [}
\def\]{\right ]}
\swapnumbers
\theoremstyle{plain}
\newtheorem{thm}{\itshape\ Theorem}[subsection]

\theoremstyle{free}
\newtheorem{prop}[thm]{}
\theoremstyle{remark}
\newtheorem{rem}{\itshape\ Remark}[thm]
\theoremstyle{free}

\newtheorem{mynumthm}{}[subsection]
\def\numthm#1{\def\themynumthm{{\itshape #1}}}
\usepackage[xdvi]{graphicx}

\title{Quantum Relatives of Alexander Polynomial}

\author{Oleg Viro}
\dedicatory{Uppsala University, Uppsala, Sweden\break
POMI, St.\ Petersburg, Russia}
\address{Department of Mathematics, Uppsala University,
Box 480, S-751 06 Uppsala, Sweden
}
\email{oleg@@math.uu.se}
\begin{document} 

\begin{abstract}The multivariable Conway function is generalized to 
oriented framed trivalent graphs equipped with additional 
structure (coloring). This is done via refinements of 
Reshetikhin-Turaev functors based on irreducible representations of 
quantized $gl(1|1)$ and $sl(2)$. The corresponding face state sum 
models for the generalized Conway function are presented.
\end{abstract}

\maketitle

\section*{Introduction}

\subsection*{Alexander Polynomial and Conway Function}The Alexander
polynomial is one of the most classical topological invariants. It
was defined  \cite{A} as early as in 1928. 
A recent paper by Fintushel and Stern \cite{FS} has again drawn attention
to the Alexander polynomial by relating it to the Seiberg-Witten
invariants of 4-dimensional manifolds.

The Alexander polynomial can be thought of in many different ways. There 
are a homological definition via free abelian covering spaces, a
definition via Fox's free differential calculus, a  definition via
Reidemeister torsion, and a diagrammatical Conway's  definition. See 
Turaev's survey \cite{Tu3}. As a matter of fact, Conway \cite{C} has enhanced
the original notion by eliminating an indeterminacy. He introduced also a
link invariant which encodes the same information as the enhanced
Alexander polynomial, but is more convenient from several points  of
view. Conway called it a {\em potential function} of a link. Following 
Turaev \cite{Tu3}, we call it the {\it Conway function.}\footnote{Some 
authors call it {\it the Alexander-Conway polynomial,\/} see e.g. 
\cite{KS}, \cite{RS1} and \cite{RS2}. However there are two reasons not
to  use this term.  First, this is not a polynomial even in the sense
in which  the Alexander polynomial is: this is not a Laurent
polynomial, but rather  a rational function. Say, for the unknot it is
$\frac1{t-t^{-1}}$. Second,  there is another invariant deserving the
name of Alexander-Conway  polynomial. It is a polynomial knot
invariant, which was introduced by  Conway in the same paper \cite{C}:
the Alexander polynomial normalized and  represented as a polynomial in
$z=t^{\frac12}-t^{-\frac12}$, where $t$ is  the original variable.
Therefore following Turaev \cite{Tu3}, we use the  term of {\it Conway
function.}} 

The definitions of the (enhanced) Alexander polynomial and Conway function 
are given below in Section \ref{sec7.7}. Here let me recall only their 
formal appearance. For an 
oriented link $L\subset S^3$ with connected components $L_1$, \dots, $L_k$ 
the Alexander polynomial $\GD_L(t_1,\dots,t_k)$ is a Laurent polynomial in 
variables $t_1^{\frac12}$, \dots, $t_k^{\frac12}$, and the Conway function 
$\nabla(L)(t_1,\dots,t_k)$ is a rational function in variables $t_1$, 
\dots, $t_k$. They are related by formulas 
$\nabla(L)(t_1,\dots,t_k)=\GD_L(t_1^2,\dots,t_k^2),\text{ if }k>1 $ and 
$\nabla(L)(t)=\frac{\GD_L(t^2)}{t-t^{-1}}, \text{ if $L$ is a knot, i.e. 
$k=1$.} $ 

\subsection*{Alexander Polynomial and Quantum Topology} The intervention 
of the quantum field theory into the low-dimensional topology, which was 
initiated in the mid eighties by discovery of the Jones polynomial, added a 
new point of view on the Alexander polynomial. It was included 
into series
of other quantum polynomial link invariants such as the Jones, HOMFLY and 
Kauffman polynomials. 

Invariants of a topological object which are studied by the quantum 
topology are presented by explicit formulas inspired by the quantum field 
theory. A versatility of these formulas allows one to generalize the 
invariants to wider classes of objects and find counterparts of the 
invariants for completely different objects. For example, the Jones 
polynomial was first generalized to colored links, then to trivalent 
framed knotted graphs, then its counterparts (Reshetikhin-Turaev-Witten 
invariants) of closed oriented 3-manifolds  were discovered, and they were 
upgraded to TQFT'es, i.e. functors from categories of closed surfaces and 
their cobordisms with an appropriate additional structures to the category 
of finite-dimensional vector spaces.  The invariant of 3-manifolds 
(Witten-Reshetikhin-Turaev invariant) was defined technically even as an 
invariant of 4-manifolds with boundary, see \cite{TuBk}. Although its 
value depends only of the boundary, the construction deals with the 
4-manifold. To some extent the same scheme works for any quantum link 
polynomial. 

However, the quantum theory of the Alexander polynomial has not evolved to 
the extents which have been achieved by the theory based on the Jones 
polynomial. Only the very first steps of this program have been made: it 
has been shown that the Alexander polynomial can be defined via 
representations 
of quantum groups. The multivariable Conway function was obtained  along 
the lines of the construction of link invariants related to a quantum 
group by Jun Murakami \cite{Mu1}, \cite{Mu2}, Rozansky and Saleur 
\cite{RS1} and Reshetikhin \cite{Re}. Rozansky, Saleur and 
Reshetikhin applied this construction to  quantum 
supergroup $gl(1|1)$ (i.e., a quantum deformation $U_qgl(1|1)$ of the 
universal enveloping algebra of the superalgebra $gl(1|1)$, see 
Appendix A to \cite{RS2}). Murakami did this with the quantum group $sl(2)$ at 
$q=\sqrt{-1}$. Deguchi and Akutsu \cite{DA} included the Conway function 
in a family of invariants related to representations of the quantum 
$sl(2)$ at 
roots of unity, generalized these invariants to colored framed trivalent 
graphs, and found face state sum models for the generalizations. However 
the case of the Conway function has not been studied separately from the 
generalizations. 

Thus, the multivariable Conway function was studied by the methods of 
quantum topology along two lines: based on quantum supergroup $gl(1|1)$ by 
Reshetikhin, Rozansky and Saleur and quantum group $sl(2)$ by Jun 
Murakami, Deguchi and Akutsu, respectively. The $sl(2)$ approach was 
developed up to explicit generalization to colored framed trivalent graphs 
and face state sum models, although the results were given in the form of 
quite complicated formulas, which have not been analyzed from a geometric 
point of view. In the frameworks of $gl(1|1)$ approach, a generalization to 
colored framed trivalent graph has not been considered. 

\subsection*{Results}The goal of this paper is a further development 
of quantum invariants related to the Alexander polynomial. 
We make  the next two steps in the $gl(1|1)$ direction: generalize the 
(multivariable) Conway function to trivalent graphs equipped with some 
additional structures (colorings) and find a face state sum 
presentation for this generalization (and, in particular, the Conway 
function of classical links).  In the $sl(2)$ direction the 
corresponding steps seem to have been done in \cite{DA}, but, for the 
sakes of a subsequent  development, we do them all over again from
scratch, in effort to get  simpler formulas and geometric
presentations. In particular, this  allows us to compare the results in
$gl(1|1)$ and $sl(2)$ directions. 

None of them could be reduced to another one.  Colorings used in 
$gl(1|1)$ approach are based on a larger palette.  Colorings based on 
a part of $gl(1|1)$-colors can be turned into $sl(2)$-colorings, but 
not all the $sl(2)$-colorings can be obtained in this way.  Larger 
palette makes $gl(1|1)$ approach somehow more flexible, but leads to 
more complicated face state sum models.  Both directions seem to 
deserve investigation, since they may be useful for different 
purposes.  

The relatives of the Alexander polynomial studied in this paper
are: \begin{itemize}
\item the Reshetikhin-Turaev functors $\RT^1$ and $\RT^2$ based on the
irreducible representations of quantum $gl(1|1)$ and $sl(2)$, 
respectively,
\item a modification $\A^c$ of $\RT^c$ for $c=1,2$,
\item invariants $\UD^c$ of closed colored framed generic graphs 
similar
to the Conway function $\nabla$ of links and, in a sense, generalizing it.
\end{itemize}

Further steps may be a construction of state sum invariants for shadow 
polyhedra and, eventually, 4-manifolds. However, there are several 
obstacles, which make this, at least, not completely straightforward. 
It is not clear yet, which family of generalizations of the Conway 
function is more suitable for that. 

Although the possibility of the future four-dimensional applications 
was the main motivation for this work, I believe the results are 
interesting on their own. The formulas involved in the face model 
appear to be much simpler than their counterparts even for the Jones 
polynomial. 

\subsection*{Reshetikhin-Turaev Functors}
As it was shown by Reshetikhin and Turaev \cite{RT1}, the construction of
link invariants related to a quantum group is upgraded to a construction
of functor from a category, whose morphisms are colored framed generic
graphs, to a category of finite-dimensional representations of a quantum
group. A restriction of this functor to a subcategory of tangles (graphs
consisting of disjoint circles and intervals) is described in \cite{RT1}
explicitly, while for graphs with trivalent vertices Reshetikhin and
Turaev \cite{RT1} leave a choice.

As far as I know, in the case of quantum $gl(1|1)$ the choice has never
been made, although the Reshetikhin-Turaev functor from the category of
tangles to the category of representations of quantum $gl(1|1)$ and its
relation to the Conway function was alluded in several papers, see e.g.
Rozansky and Saleur \cite{RS2}.

In this paper a Reshetikhin-Turaev functor is constructed explicitly, by
presenting Boltzmann weights (see Tables \ref{tab:BW1} and \ref{tab:BW2})
on the whole category of colored framed generic graphs. The choice left
by Reshetikhin and Turaev \cite{RT1} can be made in many ways, but I
believe the one made here deserves a special attention, since the resulted
functor $\RT^1$ has special nice properties. For example, the coefficient
in skein relation \ref{Aphi-gr-skein} is very simple
and the Boltzmann weights are polynomial in powers of $q$.

Special features of quantized $gl(1|1)$ make the construction
interesting and not completely straightforward. Since $gl(1|1)$
is a superalgebra, all the formulas change in a peculiar way, which looks
strange for anybody but experienced super-mathe\-ma\-ti\-ci\-an (i.e., 
a
mathematician working with superobjects). In turn, the unusual algebra
causes unusual topological features. For example, colorings of graphs
involve orientations of 1-strata in more subtle way, and at some vertices 
a cyclic order of the adjacent 1-strata must also be included into a 
coloring.

\subsection*{Modifications of Reshetikhin-Turaev Functor} The obvious 
polynomial nature of the Boltzmann weights in Tables \ref{tab:BW1} and 
\ref{tab:BW2} and the polynomial nature of the Alexander polynomial 
suggest a modification of $\RT^1$ which gives rise to a functor $\A^1$ 
which acts from almost 
the same category of colored framed generic graphs to the category of 
finite-dimensional free modules over some commutative ring $B$. For 
example $B$ may be $\Z[M]$ where $M$ is a free abelian multiplicative 
group. In this case $B$	is 
the ring of Laurent polynomials. Thus $\A^1$ is closer to the 
Alexander polynomial. 

Roughly speaking, $\A^1$ is obtained from $\RT^1$ by eliminating the 
Planck
constant $q$ (parameter of the quantum deformation in the quantization
$U_qgl(1|1)$) by replacing the powers of $q$ with independent variables
in the Boltzmann weights. The Boltzmann weights for $\A^1$ are 
presented in Tables \ref{tab:BW3} and \ref{tab:BW4}.

The transition from $\RT^1$ to $\A^1$ does not eliminate 
the quantum
nature together with $q$. Although $U_qgl(1|1)$ does not act in the 
$B$-modules which are the values of $\A$ 
on objects, there is a Hopf subalgebra of $U_qgl(1|1)$ 
such that
$\A^1$  can be upgraded to functors from the same category to 
the categories of modules over this subalgebra.

\subsection*{Alexander Invariant} The scheme which was used by Rozansky 
and Saleur \cite{RS1} and \cite{RS2} for relating the Conway function with 
the quantum $gl(1|1)$ version of the Reshetikhin-Turaev functor for 
tangles, transforms $\CC^1$ to the {\em Alexander invariant} $\UD$ 
of closed colored framed generic graphs in $\R^3$. 

$\UD^1$ assigns to a closed colored framed generic graph $\GG$ in 
$\R^3$ an 
element $\UD^1(\GG)$ of $B$. The coloring of $\GG$ consists of 
orientation 
of the 1-strata of $\GG$, assigning to each of the 1-strata a {\em 
multiplicity} taken from $\{t\in M\mid t^4\ne1\}$ and an integer 
weight,
and fixing a cyclic order of 1-strata adjacent to some vertices. 
Dependence of the Alexander invariant on the weights, multiplicities and 
cyclic orders is completely understood and described below. There exist 
universal multiplicities on $\GG$ with $M=H^1(\GG;\Z)$ such that, for 
given orientations, the Alexander invariant of $\GG$ colored with these 
orientations, arbitrary cyclic orders, weights and multiplicities can be 
recovered from the Alexander invariant of $\GG$ colored with the same 
orientations, the universal multiplicities and any cyclic orders and 
weights. Orientations constitute the most subtle part of colorings. This 
is demonstrated in case of the 1-skeleton of a tetrahedron considered in 
Section \ref{sec7.2}. However the Alexander invariant of this graph 
with any 
coloring is much simpler than the Jones polynomial of the same graph, 
which is basically the $6j$-symbol and plays an important role in the TQFT 
based on quantum $su(2)$.

For a link $L$ colored in an appropriate way, 
$$\nabla(L)(t_1,\dots,t_k)=\UD^1(t_1^2,\dots,t_k^2).$$ Thus $\UD^1$ 
can be 
considered as a generalization of the Conway function to graphs. Exactly 
as in the transition from the Alexander polynomial to the Conway function, 
here we have to replace the previous variables by their square roots. 
Hence, the variables in the Alexander invariant are quartic roots of the 
variables in the Alexander polynomial. Compare with the relation between 
the Kauffman brackets and Jones polynomial. 

\subsection*{Face Models} All formulas of quantum topology can be divided 
roughly into two classes: vertex and face state sums. Face state sums seem 
to be more versatile. At least, they are more uniform for topological 
objects of different nature. The formulas in the definition of the 
Alexander invariant are of the vertex type. 

In this paper face state sums representing the Alexander invariant are 
also obtained. This is done via a version of ``transition to the shadow 
world'' invented by Kirillov and Reshetikhin \cite{KR}  to obtain a face 
state sum model for the quantum $sl(2)$ invariant of colored framed 
graphs, generalizing the Jones polynomial. The construction used here is 
another special version of a more general construction, which  I found 
analyzing Kauffman's ``quantum spin-network'' construction \cite{Ka7} of 
Turaev-Viro invariants \cite{TV}, while presenting this in my UC San Diego 
lectures in the Spring quarter of 1991. I presented the general 
construction in several talks, but never published, since in the full 
generality it looks too cumbersome, while in the special cases, which I 
knew before this work and in which it looks nice, the result had been 
already known.

\subsection*{Thanks} I am grateful to Lev Rozansky for his valuable 
consultations concerning super-mathematics, $gl(1|1)$ and papers 
\cite{RS1} and \cite{RS2} and Alexander Shumakovitch for pointing out a 
mistake in a preliminary version of this paper. 

\tableofcontents

\section{Geometric Preliminaries on Knotted Graphs}\label{sec1}

\subsection{Generic Graphs}\label{sec1.1} By a {\it generic graph} we mean a 
1-dimensional {\it CW}-complex such that each of 
its points has a neighborhood homeomorphic either to $\R$, or half-line 
$\R_+$, or the union of three copies of $\R_+$ meeting in their common end 
point. Each generic graph is naturally stratified. Strata of dimension 1 
are the connected components of the set of points which have neighborhoods 
homeomorphic to $\R$. The 0-strata of a generic graph $\GG $  are the 
3-valent and 1-valent vertices of $\GG$. The former are called {\it 
internal vertices\/} and the latter {\it end points\/} or {\it boundary 
points\/} of the graph. The set of boundary points of $\GG$ is called the 
{\it boundary\/} of $\GG$ and denoted by $\p \GG$, the complement 
$\GG\sminus\p\GG$ is called the {\it interior\/} of $\GG$ and denoted by 
$\Int \GG$. A generic graph $\GG$ with empty boundary is said to be {\it 
closed.\/} The 1-strata homeomorphic to $\R$ are called {\it edges.\/} 
(Thus a component of a generic graphs may contain no vertex.) 

A generic graph $\GG$ is said to be {\it properly embedded\/} into a
3-manifold $M$, if $\GG\subset M$ and $\p\GG=\GG\cap\p M$.

\subsection{Additional Structures on Generic Graphs}\label{sec1.2}Generic 
graphs embedded in a 3-manifold are usually equipped with various
additional structures. The most common of the structures are orientations
and framings. An {\it orientation\/} of a generic graph $\GG$ is an
orientation of all the 1-strata of $\GG$.

Let $\GG$ be a generic graph properly  embedded into a 3-manifold $M$. An 
extension of $\GG$ to an embedded compact surface $F$ in which $\GG$ is 
sitting as a deformation retract is called a {\it framing of $\GG$\/} (in 
$M$), if $\p F\cap\p\GG= F\cap\p\GG$ and each component of $F\cap \p M$ is 
an arc containing exactly one boundary point of $\GG$. These arcs comprise 
a framing of $\p\GG$ in $\p M$. By a framing of a discrete subset of a 
surface we understand a collection of disjoint arcs embedded in the 
surface in such a way that each arc contains exactly one point of the 
subset and each point is contained in one of the arcs. 

By an isotopy of a framed generic graph we mean an isotopy of the
graph in the ambient 3-manifold extended to an isotopy of the framing. In
the case of properly embedded framed generic graph, the isotopy is assumed
to be fixed on the intersection with the boundary of ambient manifold.

The notion of framed generic graph generalizes the notion of framed
link. Indeed, a collection of circles can be considered as a generic
graph, and, if the framing surface is orientable, the framing is defined
up to isotopy by a non-zero vector field normal to the circles and
tangent to the surface. Recall that usually by a framing of a link one
means a non-zero normal vector field on it. In the case when the framing
surface is not orientable, it is related in a similar way to a framing
with a field of normal lines, rather than normal vectors.

The third additional structure that is involved below is {\it orientation 
at vertices.\/} A generic graph is said to be {\it oriented at a 
(trivalent) vertex} if  the germs of edges adjacent to the vertex are 
cyclically ordered. We speak here of germs rather than edges, because an 
edge may be adjacent to a vertex twice and then there are only two 
adjacent edges (recall that a set of two elements may not be {\it 
cyclically} ordered). If the graph is framed then its orientation at 
a vertex defines a local orientation of the framing surface at the vertex. 
                                     
\subsection{Diagrams}\label{sec1.3}To describe a generic graph in $\R^3$
up to isotopy, one uses a natural generalization of link diagrams. Let
$\GG$ be a generic graph embedded into $\R^3$. A projection of $\GG$ to
$\R^2$ is said to be {\it generic,\/} if
\begin{enumerate}
\item its restriction to each 1-stratum of $\GG$ is an immersion,
\item it has no point of multiplicity $\ge3$,
\item no double point is the image of a vertex of $\GG$,
\item at each double point the  images of 1-strata intersect each other
transversally,
\item at the image of a vertex no two branches of 1-strata are tangent
to each other.
\end{enumerate}
It is clear that the term {\it generic\/} is appropriate here in the sense
that embedded generic graphs with non-generic projection to $\R^2$ make a
nowhere dense set in the space of all embedded generic graphs: by
arbitrarily small isotopy one can make any embedded generic graph having
generic projection and any generic projection cannot be made non-generic
by a sufficiently small isotopy.

To get a description up to isotopy of a generic graph embedded in
$\R^3$, one should enhance its generic projection with information which 
branch is over and which is under the other one at each double point.
If it is done by breaking the lower branch,  the picture obtained is
called a {\it diagram\/} of the graph, like in the case of links.

\subsection{Strata of Diagram}\label{sec1.4}The image of a generic graph 
$\GG$ under a generic projection is a graph
embedded in $\R^2$. Its vertices are the double points of the projection
and the images of the vertices of $\GG$. The former are called also {\it
crossing points.\/} The natural stratification of this graph is extended
to a stratification of $\R^2$ with 2-strata being connected components
of the complement of the projection. Strata of this stratification are
called  {\it strata of the diagram,\/} the 2-strata are also called
{\it faces\/} of the diagram.

\subsection{Moves of Diagrams}\label{sec1.5}
\begin{figure}[ht]
\centerline{\includegraphics[width=5in]{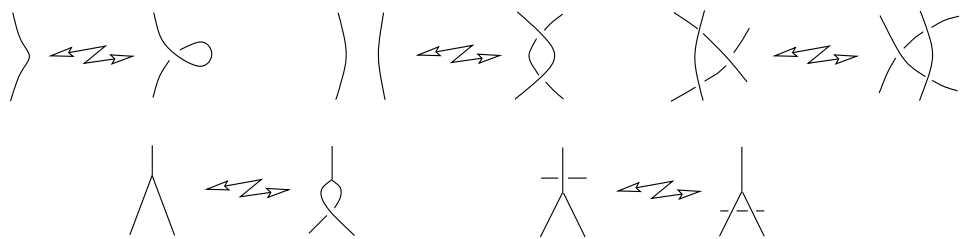}}
\caption{Elementary isotopies of generic graphs.}
\label{f1}
\end{figure}

Diagrams of isotopic generic 
graphs embedded in $\R^3$ can be obtained from each other by a sequence of 
transformations of 5 types shown in Figure \ref{f1}. Here two 
transformations are considered of the same type if the modifications 
happen in disks, which can be mapped to each other by a homeomorphism 
moving one picture to another one followed maybe by a simultaneous change 
of all the crossings (this change corresponds to reflection in the plane 
parallel to the plane of the diagram).  The types correspond to the 
conditions of the definition of generic projection: each of the 
transformation of Figure \ref{f1} can be realized by an isotopy under 
which exactly one of the conditions is violated exactly once and in the 
simplest manner. The first 3 transformations are the well-known 
Reidemeister moves. They do not involve vertices and can be used to obtain 
diagrams of isotopic links from each other. 

\subsection{Showing Additional Structures on Diagram}\label{sec1.6}
Any generic graph with generic projection to $\R^2$ can be equipped with a
framing such that the projection to $\R^2$ restricted to the framing
is an immersion (recall that the framing is a surface). This framing is
unique up to isotopy. It is called  {\it black-board\/} framing.

Similarly, a generic graph generically projected to $\R^2$ can be equipped 
with the {\it black-board\/} orientation at a vertex: counter-clockwise 
cyclic ordering of the adjacent germs of edges. An arbitrary orientation 
at vertices can be shown in diagram by circular arrows around the vertices 
following the cyclic orders. For the sake of simplicity, we will skip all 
the counter-clockwise arrows. 

Let $\GG$ be a framed generic graph with a generic projection to $\R^2$.  
On each 1-stratum of $\GG$ one can compare two framings: the black-board  
and original ones. The difference is an integer or half-integer number.   
It is the number of  (full) right twists which should be added to the 
black-board framing on the  1-stratum to get the framing involved into the 
original framing of $\GG$.  To describe a framed generic graph up to 
isotopy it is enough to draw its diagram  equipped with these numbers 
assigned to all 1-strata of the graph. Another way to show the difference 
between the framing and the black-board one is to put on the arc of the 
diagram the fragment ${\fig{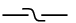}}$ at a positive half-twist 
and 
the fragment ${\fig{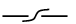}}$ at a negative half-twist. Of 
course, 
the isotopy type of a framing is described by the total number of the 
half-twists on each component. However sometimes, especially when 
describing changes of a diagram, it is convenient to localize the places 
where the framing is orthogonal to the plane of the diagram and specify 
the local behavior of the framing at these places, as the pictures 
${\fig{plhalftw.eps}}$ and ${\fig{nghalftw.eps}}$ do. By {\it a 
diagram of 
framed graph\/} we call a diagram enhanced by the differences between the 
framing and the black-board framing, which are presented either by the 
numbers or by pictures of half-twists. 

\begin{figure}[htb]
\centerline{\includegraphics[width=5in]{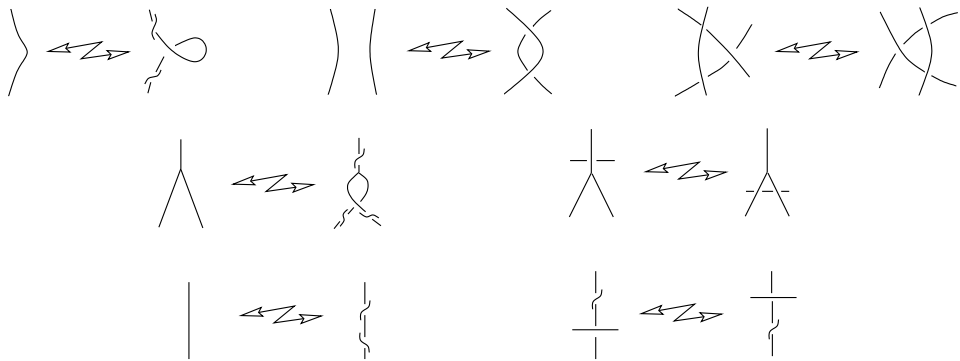}
} \caption{Elementary isotopies of framed generic graphs.}
\label{f2}
\end{figure}

Diagrams of isotopic framed generic graphs can be  obtained from each 
other by a sequence of transformations of 7 types\footnote{Up to the same 
rough equivalence as above, where isotopies of unframed graphs were 
considered.} shown in Figure \ref{f2}. The first 5 of them are the 
transformations of Figure \ref{f1} equipped with description of behavior 
of the framing. The last 2 describe emerging and moving of half-twists. 

The notion of diagram is extended obviously to generic graphs and
framed generic graphs lying in a product of a surface $P$ by $\R$
or $[0,1]$. The diagrams are drawn on $P$. To avoid situations
when diagrams are not quite understandable, let us agree to assume that
the boundary points of graphs are contained over the boundary of $P$, there
is no double point on the boundary of $P$, and framings at the end
points of graphs are black-board.

\subsection{Category of Graphs}\label{sec1.7}Framed graphs gives rise to
a category $\cal G$ which is defined as follows. Its objects are in
one-to-one correspondence with natural numbers. The object
corresponding to a natural number $k$ is the pair consisting of $\R^2$
and the set $\Pi_k$ of $k$ points $(1,0)$, $(2,0)$, \dots, $(k,0)$
equipped with a framing which consists of small intervals of the
$x$-axis. 
A morphism from the $k$-th to $l$-th object is the isotopy
class of a framed generic graph $\GG$ embedded in $\R^2\times[0,1]$
with $\p\GG=\Pi_k\times0\cup\Pi_l\times1=\p(\R^2\times[0,1])\cap\GG$. A 
framed
generic graph of this sort is called a {\it (framed) $(k,l)$-graph.\/}
The composition of two morphisms is defined by attaching them one over
another. The identity morphism of the object $(\R^2,\Pi_k)$ is
$(\R^2\times[0,1], \Pi_k\times[0,1])$.

There is another basic operation with framed $(k,l)$-graphs, which is called
{\it tensor product.\/} The tensor product of a framed $(k,l)$-graph $\GG$
by a $(m,n)$-graph $\GG'$ is an $(k+m,l+n)$-graph obtained
by placing $\GG'$ to the left of $\GG$.

\subsection{Generators and Relations of the Category}\label{sec1.8}
Framed $(k,l)$-graphs are represented by diagrams on the
strip $\R\times[0,1]$.  Applying general position arguments to the diagrams
and their deformations under isotopies, one can easily find the
natural generators and relations for the category of framed generic 
graphs described in the following theorems.

\begin{figure}[hbt]
\centerline{\includegraphics{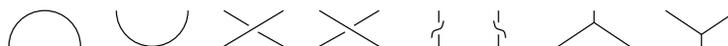}}
\caption{Generators of $\cal G$.}\label{f3}                   
\end{figure}

\begin{prop}\label{1.8.A} Any morphism of $\cal G$ can be presented as 
a composition of morphisms each of which
is a tensor product of an identity morphisms (i.~e., the isotopy class
of graph $\Pi_i\times[0,1]$ black-board framed) by one of the morphisms 
shown in Figure \ref{f3}.
\end{prop}

Indeed, by a small isotopy one can make the restriction of the second
coordinate to each 1-stratum of the plane projection of the graph to
have only isolated nondegenerated critical points. Further, the values
of the second coordinate at these points, half-twist signs and vertices
of the projection can be made pairwise distinct. After this draw
horizontal lines separating from each other the critical points,
half-twist signs and vertices. The pieces between the lines next to
each other are graphs of the desired types.\qed

Similar arguments, but applied to a one-parameter family of graphs and
combined with representing some of the elementary relations as
consequences of other elementary relations, prove the following theorem.

\begin{figure}[bht]
\centerline{\includegraphics[height=3.8in]{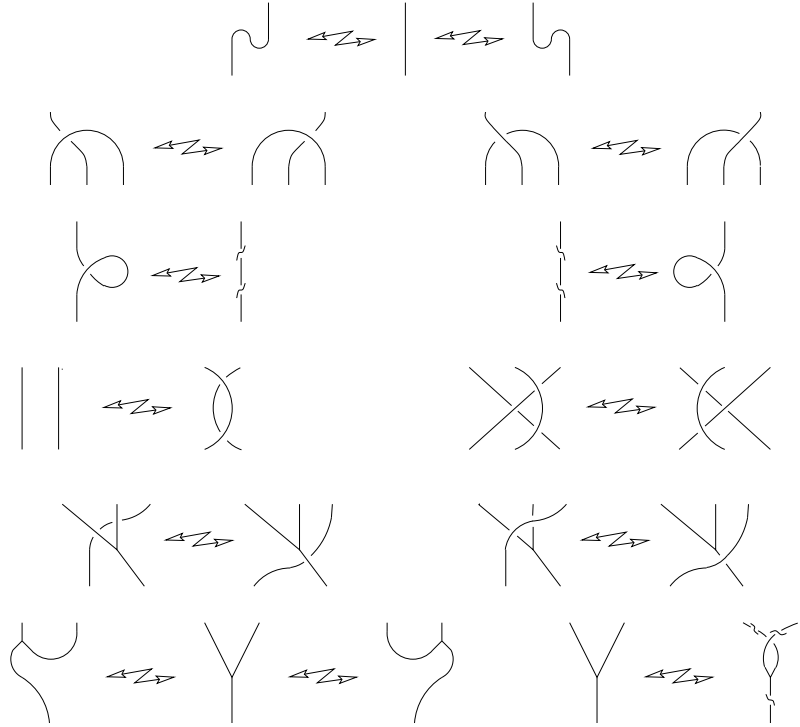}} 
\caption{Relations of $\cal G$.} \label{f4}
\end{figure}

\begin{prop}\label{1.8.B} Two products of the graphs generating $\cal G$
by Theorem \ref{1.8.A} are isotopic, iff they can be transformed to 
each other by a sequence of elementary moves shown in Figure \ref{f4}.
\end{prop}

\section{Reshetikhin-Turaev Functor Based on $gl(1|1)$
}\label{sec2}

\subsection{Category of Colored Framed Graphs. 
Objects}\label{sec2.1}
Recall that according to Reshetikhin-Turaev's scheme \cite{RT1} the
invariants are constructed via a functor from a category of colored 
framed graphs, whose objects and morphisms are the objects and morphisms
of $\cal G$ equipped with additional structures (colored), to a  
category of finite-dimensional representations of some quantum group.
Here Reshetikhin-Turaev's setup has to be slightly modified.

Let $h$ be a non-zero complex number. Below we define a category $\cal 
G_h^1$. It is denoted also by $\cal G^1$, when dependence of $h$ is 
not emphasized. Its objects are
finite sequences of triples consisting of two complex numbers, say $j$ 
and
$J$, and a sign (i.~e., $\Gd=+$ or $-$). The first number must be 
such that $2jh\sqrt{-1}/\pi$ is not an integer.
The first number is called the {\it multiplicity,\/} the second {\it 
weight.\/} 

The triples $(j_1,J_1,\Gd_1)$, \dots, $(j_k,J_k,\Gd_k)$ comprising
an object are convenient to associate with points $(1,0)$, \dots,
$(k,0)$ of the plane $\R^2$.  The object
$\{(j_1,J_1,\Gd_1)$, \dots, $(j_k,J_k,\Gd_k)\}$
is to be thought of as the set $\Pi_k$ whose points are colored with
the triples $(j_1,J_1,\Gd_1)$, \dots, $(j_k,J_k,\Gd_k)$. The signs
$\Gd_1$, \dots, $\Gd_k$ are to be interpreted as orientations of the
corresponding points.

In fact, the triples stand for irreducible $(1|1)$-dimensional 
$U_qgl(1|1)$-modules. Recall (see Appendix 1, Section \ref{sL.4}) that 
irreducible $(1|1)$-dimensional $U_qgl(1|1)$-modules are numerated by a 
sign $\pm$ and a pair of complex parameters $(j,J)$ with 
$j\in\C\sminus\{\pi n\sqrt{-1}/2h : n\in\Z\}$ and $J\in\C$. Here $h$ is 
the parameter of quantization (the Planck constant).
 
\subsection{Category of Colored Framed Graphs. 
Morphisms}\label{sec2.2} A morphism of $\cal G_h^1$ 
$$\{(i_1,I_1,\Ge_1),\dots (i_k,I_k,\Ge_k)\}\to \{(j_1,J_1,\Gd_1),\dots 
(j_l,J_l,\Gd_l)\}$$ 
is a morphism of $\cal G$ (see Section \ref{sec1.7}) equipped with some 
additional structures described below. Recall that a morphism of $\cal G$ 
is a framed  generic graph (see Section \ref{sec1.1}) $\GG$ embedded in 
$\R^2\times[0,1]$ with 
$\p\GG=\Pi_k\times0\cup\Pi_l\times1=\p(\R^2\times[0,1])\cap\GG$ considered 
up to ambient isotopy fixed on $\p\GG$. Its framing at the end points is 
parallel to the $x$-axis. 

The additional structures on $\GG$ all together are called {\em 
$\cal G^1$-coloring}.
Each of the 1-strata is colored with an orientation and a pair of 
complex 
numbers the first of which does not belong to 
$\{\pi n\sqrt{-1}/2h : n\in\Z\}$. 
These numbers are also called the
{\it multiplicity\/} and {\it weight,\/} respectively. The orientations of 
the edges adjacent to the boundary points are determined by the signs of 
the points' colors: the edge adjacent to a point with $+$ is oriented at 
the point upwards (along the standard direction of the $z$-axis), 
otherwise it is oriented downwards. The multiplicity and weight of the 
edge adjacent to an end point coincide with the corresponding ingredients 
of the color of this end point. At each of the trivalent vertices the 
colors of the adjacent edges satisfy to the following 

\begin{prop}[Admissibility Conditions]\label{2.2.A}
Let  $(j_1,J_1)$, $(j_2,J_2)$,
$(j_3,J_3)$ be the multiplicity and weight components of the colors of
three edges whose germs are adjacent to the same vertex and
let $\Ge_i=-1$ if the $i^{th}$ of the edges is oriented towards the vertex
and $\Ge_i=1$ otherwise. Then
\begin{equation}\label{j-admiss}
\sum_{i=1}^3\Ge_ij_i=0
\end{equation}
\begin{equation}\label{J-admiss}
\sum_{i=1}^3 \Ge_iJ_i=-\prod_{i=1}^3\Ge_i
\end{equation}
\end{prop}

A vertex with all the adjacent edges oriented towards it (i.~e., 
$\sum_{i=1}^3\Ge_i=-3$) or in the opposite direction (i.~e., 
$\sum_{i=1}^3\Ge_i=3$) is said to be {\em strong}. At each strong vertex 
the graph $\GG$ is oriented. The orientation of $\GG$ at strong vertices 
is considered a part of the coloring of $\GG$, together with the 
orientations, multiplicities and weights of 1-strata.  

Equation \eqref{j-admiss} looks like the Kirchhoff equation on currents. 
More topological interpretation: the boundary of chain $\sum j_ie_i$, 
where the sum is taken over all the edges $e_i$ of the graph (with their 
orientations) and $j_i$ is the first component of the color of $e_i$, is 
the 
sum of all the boundary points of the graph $\GG$ taken with coefficients 
which are the multiplicities of these points with appropriate signs. In 
particular,  $\sum j_ie_i$ is a relative cycle with real coefficients of 
$\GG$ modulo $\p\GG$, so that it defines an element of 
$H_1(\GG,\p\GG;\C)$. 

\begin{rem}\label{2.2.B} In the Reshetikhin-Turaev setup the orientation 
of an edge could be reversed along with simultaneous change of the 
numerical component of the color. Here we cannot do this. The orientations 
are to be thought as inseparable ingredients of the colors. Reversing the 
orientation of a single 1-stratum cannot be compensated by changes of 
other components of the coloring. 

Orientations at strong vertices are not that significant. Reversing one 
of them causes multiplication by $-1$ of the image under the 
Reshetikhin-Turaev functor, see Section \ref{sec2.5}. Therefore many 
formulas below apparently do not involve orientations at vertices. 
Orientations at vertices are hidden by eliminating strong vertices with 
clockwise orientations. Recall that we have agreed not to show 
counter-clockwise orientations at vertices in a diagram.  \end{rem} 

A vertex with odd number of adjacent edges oriented towards it (i.~e., 
$\prod_{i=1}^3\Ge_i=-1$) is said to be {\em odd}. If one thought on 
weights also as on currents, then \eqref{J-admiss} would mean that each 
odd vertex is a source of capacity $1$ and an even vertex a sink of 
capacity~$1$. In other words, $\p\sum J_ie_i$ involves each odd 
vertex with coefficient $-1$ and each even vertex with coefficient 
$1$. 

The composition of two morphisms is defined by attaching them one over
another. There is another basic operation with morphisms, which is called
{\it tensor product.\/} The tensor product of a morphism $\GG$ by a
morphism $\GG'$ is a morphism obtained by placing $\GG'$ to the left of
$\GG$.

There is a natural forgetting functor $\cal G_h^1\to\cal G$. 
Morphisms of $\cal G_h^1$ mapped by this functor to the generators 
of $\cal G$ generate $\cal G_h^1$. Therefore Theorem \ref{1.8.A} 
implies the following.

\begin{prop}[Generators of $\cal G_h^1$]\label{2.2.C} Any morphism 
of $\cal G_h^1$ can be presented as a composition of morphisms each 
of which
is a tensor product of an identity morphisms (i.~e., the isotopy class
of graph $\Pi_i\times[0,1]$ black-board framed and appropriately
colored) by one of the morphisms whose underlying graph is shown in
Figure \ref{f3}. \end{prop}

Similarly, Theorem \ref{1.8.B} gives relations of $\cal G_h^1$.

\subsection{Refined Reshetikhin-Turaev Functor on Objects}\label{sec2.3} 
Denote by $\cal R_q$ the category of finite-dimensional modules over the
$q$-deformed universal enveloping algebra $U_qgl(1|1)$ of superalgebra
$gl(1|1)$. Our goal is to construct a functor $\RT^1:\cal G_h^1\to\cal 
R_q$.

For an object $\{(j_1,J_1,\Ge_1),\dots(j_k,J_k,\Ge_k)\}$ of $\cal 
G_h^1$, the image under $\RT^1$ is defined to be the tensor product
$(j_1,J_1)_{\Ge_1}\otimes\dots\otimes (j_k,J_k)_{\Ge_k}$ of the
irreducible $(1|1)$-dimensional $U_qgl(1|1)$-modules.

Each of the factors $(j_i,J_i)_{\Ge_i}$
is a 2-dimensional vector space with a canonical basis. The vectors $e_0$
and $e_1$ comprising the basis are distinguished by their nature: $e_0$ is
a boson, $e_1$ a fermion.\footnote{Whatever this means! See Appendix 1 
for brief explanations.} Thus
$(j_1,J_1)_{\Ge_1}\otimes\dots\otimes (j_k,J_k)_{\Ge_k}$ is a complex
vector space of dimension $2^k$.  Its basis consists of
vectors $e_{i_1}\otimes\dots\otimes e_{i_k}$, so the basis vectors may be
identified with sequences of zeros and ones placed at the points of
$\Pi_k$.

\subsection{Digression: Boltzmann Weights}\label{sec2.4}
To describe the images of morphisms under $\RT^1$, one uses geometric
interpretation borrowed from statistical mechanics. Recall how this works.

The image under $\RT^1$ of an arbitrary morphism of $\cal G_h^1$ is a 
linear map between the corresponding vector space.  The map can be 
presented by a matrix.  Elements of the matrix are complex numbers 
corresponding to the pairs consisting of basis vectors of the two 
vector spaces.  Recall that the basis vectors are interpreted as 
sequences of the zeros and ones.  These zeros and ones are placed at 
the end points of the colored framed graph representing the morphism.  
Hence, to define the image of this morphism under $\RT^1$, we need to 
associate a complex number with any distribution of zeros and ones at 
the end points of the graph.  

Of course, it suffices to define the images of all the generators of 
$\cal G_h^1$ (see Figure \ref{f3}) in such a way that the images 
would satisfy the relations shown in Figure \ref{f4}.  

For a graph underlying a generator, any distribution of zeros and ones 
at the end points naturally extends to the adjacent 1-strata of the 
diagram.  The matrix element of the morphism which is the image of the 
generator of $\cal G_h^1$ corresponds to such a coloring of 1-strata 
of the diagram.  

We deal now with two colorings: the first one is a $\cal G^1$-coloring, 
the second one associates zero or one to each 1-stratum of the 
diagram.  These colors stand for basis vectors of the module 
associated to the $\cal G^1$-color of the 1-stratum.  To distinguish 
these colorings from each other, we call the first of them (i.~e., a 
$\cal G^1$-coloring) the {\em primary coloring,} and the second one, 
the {\em secondary coloring}.  

Hence to define the image under $\RT^1$ of the generator of $\cal 
G^1$, we need to associate with any such a pair of colorings of the 
elementary graphs shown in Figure \ref{f3} an appropriate complex 
number.  In the statistical mechanics, where this graphical 
interpretation comes from, these numbers are called {\it Boltzmann 
weights.\/} 
 
Once Boltzmann weights are known for all the pair of colorings of 
elementary graphs, the matrix element for the image of an arbitrary 
morphism of $\cal G_h^1$ can be calculated as follows: choose a 
diagram
for the corresponding graph, consider all the extensions of the given 
distributions of zeros and ones at the end points of the graph to 
secondary colorings of the whole diagram.  For each of the secondary 
colorings make the product of all the Boltzmann weights of its 
elementary fragments.  The sum of the products over all these 
secondary colorings is the desired matrix element.  A secondary 
coloring in the terminology of statistical mechanics would be called a 
{\it state\/} (of the graph).  Hence the sum of the products is a sum 
over all the states, or a {\it state sum.\/} 

\subsection{Refined Reshetikhin-Turaev Functor on Morphisms}\label{sec2.5}
The Main Theorem of Reshetikhin-Turaev paper \cite{RT1} suggests the
images for the generators with the six underlying graphs shown on the
left hand side of Figure \ref{f3}. For the right two this should be
the Clebsch-Gordan morphisms scaled so that the relations corresponding
to the bottom row of Figure \ref{f4} were satisfied.

The choice of scaling is not prescribed by the Reshetikhin-Turaev Theorem
and, moreover, is not unique. However, properties of the functor
depend on the choice which is made here and justify it.

The images of the generators are given below in Tables \ref{tab:BW1} 
and \ref{tab:BW2} by the Boltzmann weights.  The secondary colors are 
shown there as follows.  The arcs colored with the boson basis vector 
$e_0$ are shown as dotted arcs, the arcs colored with the fermion 
basis vector $e_1$ are shown in solid.  (We follow Kauffman and Saleur 
\cite{KS}.) 

The leftmost graph in Figure \ref{f3} can be oriented in two ways. If 
it is oriented from right to left, its image under $\RT^1$ should map 
$(j,J)_-\otimes (j,J)_+$ to $\C$. Reshetikhin-Turaev Theorem suggests 
\footnote{This differs by the sign from Reshetikhin-Turaev due to the 
framework of supermathematics, see Appendix 1.} the pairing $\ld$ 
which acts by $(e_a\otimes e_b)\mapsto(-1)^{a}q^{-2ja}\Gd_{ab}$, see  
Section \ref{sL.5}. For the same graph with the opposite orientation, 
Reshetikhin-Turaev Theorem suggests ``quantum transposed'' pairing 
$\rd$ which acts by $(e_a\otimes e_b)\mapsto q^{2j(1-a)}\Gd_{ab}$, see 
Section \ref{sL.5}. 

\begin{table}[pbth]
\centerline{\renewcommand{\arraystretch}{1.6}
\begin{tabular}{|c||c|c|c|c|}
\hline
&${\fig{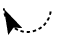}}_{(j,J)}$ & ${\fig{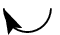}}_{(j,J)}$&
${\fig{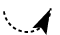}}_{(j,J)}$& ${\fig{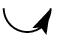}}_{(j,J)}$ \\
\hline &  $+1$ & $-q^{2j}$&$q^{-2j}$& $1$ \\ \hline\hline
&${\fig{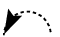}}_{(j,J)}$&
${\fig{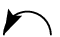}}_{(j,J)}$&${\fig{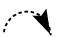}}_{(j,J)}$&
${\fig{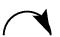}}_{(j,J)}$ \\                  \hline
&  $+1$ & $-q^{-2j}$&$q^{2j}$& $1$ \\ \hline
\hline
&$\fig{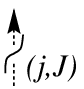}$&$\fig{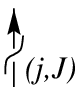}$&$\fig{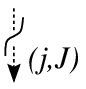}$&
$\fig{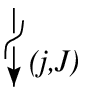}$\\ \hline
&$q^{-jJ}$&$q^{-jJ}$  &$q^{-jJ}$  &$q^{-jJ}$ \\ \hline\hline
&$\fig{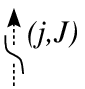}$&$\fig{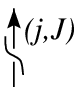}$&$\fig{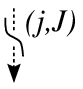}$&
$\fig{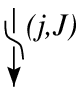}$\\ \hline
&$q^{jJ}$&$q^{jJ}$  &$q^{jJ}$  &$q^{jJ}$ \\ \hline
\hline
&${\fig{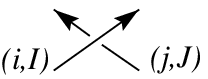}}$& ${\fig{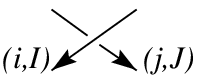}}$&
${\fig{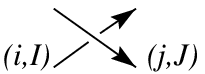}}$&
${\fig{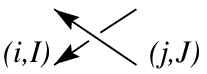}}$\\
\hline
${\fig{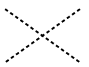}}$&  $q^{-iJ-jI}q^{i+j}$ &
 $q^{-iJ-jI}q^{i+j}$ &
 $q^{-iJ-jI}q^{i+j}$&$q^{-iJ-jI}q^{i+j}$\\
 \hline
${\fig{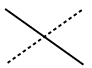}}$&  $q^{-iJ-jI}q^{j-i}$&
 $q^{-iJ-jI}q^{j-i}$&$q^{-iJ-jI}q^{j-i}$ &
 $q^{-iJ-jI}q^{j-i}$\\
 \hline
${\fig{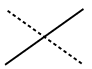}}$&  $q^{-iJ-jI}q^{i-j}$&
 $q^{-iJ-jI}q^{i-j}$&
  $q^{-iJ-jI}q^{i-j}$ & $q^{-iJ-jI}q^{i-j}$ \\
  \hline
${\fig{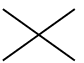}}$&  $-q^{-iJ-jI}q^{-i-j}$&
 $-q^{-iJ-jI}q^{-i-j}$&$-q^{-iJ-jI}q^{-i-j}$&
 $-q^{-iJ-jI}q^{-i-j}$ \\
\hline
${\fig{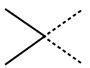}}$& 0 & $q^{-iJ-jI+i+j}(1-q^{-4i})$ &0
&0\\ \hline
${\fig{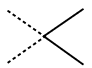}}$& $q^{-iJ-jI-i-j}(q^{4i}-1)$ &
0 &0&0\\
 \hline
${\fig{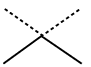}}$& 0 &0 &$q^{-iJ-jI-i+j}(1-q^{-4j})$ &0\\
 \hline
${\fig{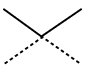}}$& 0 & 0 &0
&$q^{-iJ-jI+i-j}(1-q^{4j})$\\ \hline
\hline
&${\fig{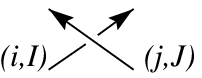}}$&
${\fig{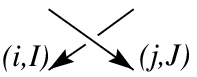}}$&
${\fig{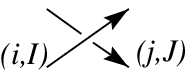}}$&
${\fig{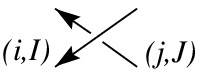}}$\\
\hline
${\fig{bx0000.eps}}$&  $q^{iJ+jI}q^{-i-j}$ &
 $q^{iJ+jI}q^{-i-j}$ &
 $q^{iJ+jI}q^{-i-j}$&$q^{iJ+jI}q^{-i-j}$\\
 \hline
${\fig{bx0110.eps}}$&  $q^{iJ+jI}q^{i-j}$&
 $q^{iJ+jI}q^{i-j}$&$q^{iJ+jI}q^{i-j}$ &
 $q^{iJ+jI}q^{i-j}$\\
 \hline
${\fig{bx1001.eps}}$&  $q^{iJ+jI}q^{j-i}$&
 $q^{iJ+jI}q^{j-i}$&
  $q^{iJ+jI}q^{j-i}$ & $q^{iJ+jI}q^{j-i}$ \\
  \hline
${\fig{bx1111.eps}}$&  $-q^{iJ+jI}q^{i+j}$&
 $-q^{iJ+jI}q^{i+j}$&$-q^{iJ+jI}q^{i+j}$&
 $-q^{iJ+jI}q^{i+j}$ \\
\hline
${\fig{bx1010.eps}}$&   $q^{iJ+jI-i-j}(1-q^{4j})$ &0&0
&0\\ \hline
${\fig{bx0101.eps}}$& 0&$q^{iJ+jI+i+j}(q^{-4j}-1)$
 &0&0\\  \hline
${\fig{bx1100.eps}}$& 0 &0&0 &$q^{iJ+jI+i-j}(1-q^{-4i})$ \\
\hline
${\fig{bx0011.eps}}$& 0 & 0 &$q^{iJ+jI-i+j}(1-q^{4i})$&0
\\ \hline\end{tabular}}\vskip.1in
\caption{}\label{tab:BW1}
\end{table}

\begin{table}[ptbh]
\centerline{\renewcommand{\arraystretch}{1.5}
\begin{tabular}{|c||c|c|c|c|}\hline
&${\fig{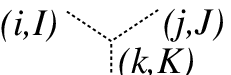}}$ &${\fig{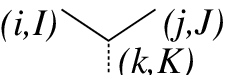}}$
&${\fig{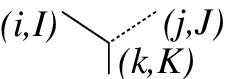}}$ &${\fig{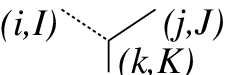}}$ \\
 \hline
${\fig{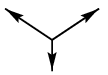}}$  &
$0$ & $1$ & $-q^{2i}$ &  $1$
\\ \hline
${\fig{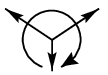}}$  &
$0$ & $-1$ & $q^{2i}$ &  $-1$
\\ \hline
${\fig{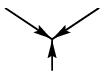}}$  &
$0$ & $q^{2k}-q^{-2k}$ & $q^{2j}-q^{-2j}$ &
$(q^{-2i}-q^{2i})q^{2j}$ \\ \hline
${\fig{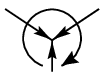}}$  &
$0$ & $q^{-2k}-q^{2k}$ & $q^{-2j}-q^{2j}$ &
$(q^{2i}-q^{-2i})q^{2j}$ \\ \hline
${\fig{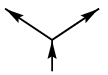}}$  &
$q^{2k}-q^{-2k}$ & $0$   & $(q^{2i}-q^{-2i})q^{-2j}$ &
$q^{2j}-q^{-2j}$ \\ \hline
${\fig{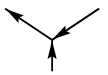}}$  &
$1$ & $-q^{2i}$ & $1$ & $0$
\\ \hline
${\fig{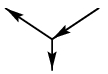}}$  &
$q^{2j}-q^{-2j}$ & $(q^{-2i}-q^{2i})q^{2j}$ & $0$ &
$q^{2k}-q^{-2k}$ \\ \hline
${\fig{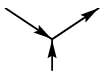}}$  &
$q^{-2i}$ & $1$ & $0$ & $1$
\\ \hline
${\fig{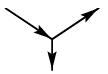}}$  &
$(q^{2i}-q^{-2i})q^{-2j}$ & $q^{2j}-q^{-2j}$ &
$q^{2k}-q^{-2k}$ & $0$ \\ \hline
${\fig{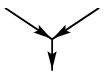}}$  &
$1$ & $0$ & $1$ & $q^{-2i}$\\
\hline\hline
&${\fig{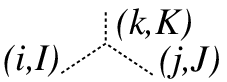}}$ &${\fig{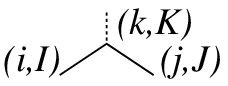}}$
&${\fig{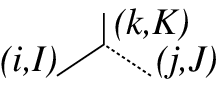}}$ &${\fig{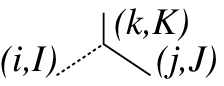}}$ \\
 \hline
${\fig{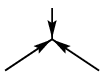}}$  &
$0$ & $q^{2k}-q^{-2k}$ & $q^{2j}-q^{-2j}$ &
$(q^{-2i}-q^{2i})q^{-2j}$ \\ \hline
${\fig{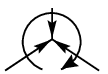}}$  &
$0$ & $q^{-2k}-q^{2k}$ & $q^{-2j}-q^{2j}$ &
$(q^{2i}-q^{-2i})q^{-2j}$ \\ \hline
${\fig{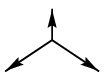}}$  &
$0$ & $1$ & $-q^{-2i}$ & $1$
\\ \hline
${\fig{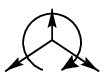}}$  &
$0$ & $-1$ & $q^{-2i}$ & $-1$
\\ \hline
${\fig{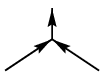}}$  &
$1$ & $0$ & $1$  & $q^{2i}$
\\ \hline
${\fig{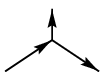}}$  &
$(q^{2i}-q^{-2i})q^{2j}$ & $q^{2j}-q^{-2j}$ &
$q^{2k}-q^{-2k}$ & $0$ \\ \hline
${\fig{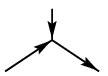}}$  &
$q^{2i}$ & $1$ & $0$ & $1$
\\ \hline
${\fig{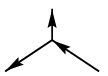}}$  &
$q^{2j}-q^{-2j}$ & $(q^{-2i}-q^{2i})q^{-2j}$ & $0$ &
$q^{2k}-q^{-2k}$ \\ \hline
${\fig{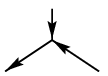}}$  &
$1$ & $-q^{-2i}$ & $1$ & $0$
\\ \hline
${\fig{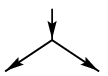}}$  &
$q^{2k}-q^{-2k}$ & $0$ & $(q^{2i}-q^{-2i})q^{2j}$ &
$q^{2j}-q^{-2j}$ \\ \hline
\end{tabular}}\vskip.1in
\caption{}\label{tab:BW2}
\end{table}

The next graph in Figure \ref{f3} oriented from right to left is 
mapped by $\RT^1$ to $\C\to(j,J)_-\otimes (j,J)_+$.  Reshetikhin-Turaev 
Theorem suggests the co-pairing $\lb:1\mapsto e_0\otimes 
e_0-e_1\otimes q^{2j}e_1$ (the minus sign is again due to the 
superenvironment).  For the same graph with the opposite orientation 
the formula is $$\rb:\C\to(j,J)_+\otimes (j,J)_-:1\mapsto 
q^{-2j}e_0\otimes e_0+e_1\otimes e_1.$$

The next two graphs in Figure \ref{f3} are diagrams of neighborhoods 
of crossing points.  At a crossing point where the strings are colored 
with modules $M$ and $N$, the Reshetikhin-Turaev Theorem suggests the 
composition of the action of the universal $R$-matrix in the tensor 
product $M\otimes N$ and the transposition $M\otimes N\to N\otimes M$.  
This composition is calculated in Appendix 1, Section 
\ref{sL.6}.  The (non-zero) Boltzmann weights are shown in 
Table \ref{tab:BW1}.  

Each of the next two graphs in Figure \ref{f3} consists of a single string
with a half-twist of the framing. We need to associate with them the
action of a square root of $v^{-1}$ and $v$, respectively,  in the module
corresponding to the color of the string. $v$ acts both in $(j,J)_+$ and
$(j,J)_-$ as multiplication by $q^{2jJ}$, see Section \ref{sL.4}.
Therefore the Boltzmann weights for ${\fig{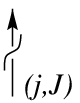}}$
and ${\fig{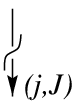}}$ are $q^{-jJ}$, and for
${\fig{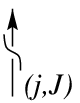}}$ and
${\fig{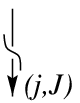}}$ they are $q^{jJ}$.

The last two graphs in Figure \ref{f3} are regular neighborhoods of a
trivalent vertex. They admit many orientations and the Boltzmann weights
heavily depend on them. As it was mentioned above, the weights
are Clebsch-Gordan coefficients scaled to make the state sum invariant with
respect to the moves shown in the last line of Figure \ref{f4}.
The Clebsch-Gordan morphisms are calculated in Section
\ref{sL.6}. 

The proof that the Boltzmann weights give rise to a well-defined 
functor can be either a reference to the Reshetikhin-Turaev Theorem 
\cite{RT1} and a routine check of invariance with respect to the moves 
of the bottom row of Figure \ref{f4}, or a more extensive routine 
check for the whole set of moves shown in Figure \ref{f4}.  Since 
either of these proofs is long and straightforward, I left them to the 
reader.

\subsection{Symmetry of Boltzmann Weights}\label{sec2.6} The Boltzmann 
weights of Tables \ref{tab:BW1} and \ref{tab:BW2} which involve 
neither weight components of colors nor clock-wise orientations at 
strong vertices satisfy a remarkable symmetry.  If one reverses 
orientations of two edges in an entry of Table \ref{tab:BW2} and 
switches simultaneously bosons and fermions on them, in the Boltzmann 
weight each expression of the form $q^{2j}$, where $j$ is the 
multiplicity of the arc involved in the change, changes to $-q^{-2j}$.  

For example, as the Boltzmann weight at $\fig{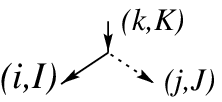}$ is 
$(q^{2i}-q^{-2i})q^{2j}$, the Boltzmann weight at 
$\fig{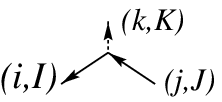}$	is 
$(q^{2i}-q^{-2i})(-q^{-2j})=(q^{-2i}-q^{2i})q^{-2j}$.

This symmetry allows one to recover the whole Table \ref{tab:BW2} from 
its first column.  The requirement to change two edges can be replaced 
by the requirement to preserve the parity of a trivalent vertex.  

The Boltzmann weights of Table \ref{tab:BW1} at points of maxima and 
minima satisfy a symmetry of the same sort: as the Boltzmann weight at 
${\fig{lb-f.eps}}_{(j,J)}$ is $-q^{2j}$, the Boltzmann weight at 
${\fig{rb-b.eps}}_{(j,J)}$ is $q^{-2j}$.

\subsection{How to Eliminate Planck Constant}\label{sec2.7}
All Boltzmann weights of Tables \ref{tab:BW1} and \ref{tab:BW2} are linear 
combinations of powers of $q$. The exponents are linear combinations of 
the multiplicities and products of a multiplicity by a weight. The number 
$q$ never appears without an exponent of this sort.  

The Boltzmann weights of Tables \ref{tab:BW1} and \ref{tab:BW2} 
satisfy identically the equations which express invariance of the 
state sums with respect to the moves of Figure \ref{f4}.  These 
equations have shape of integer polynomial equations in powers of $q$ 
with exponents $\pm j_i$ and $J_kj_l$, where $j_1$, \dots, $j_n$ are 
the multiplicities and $J_1$, \dots, $J_n$ weights of the colors 
involved.  The equations are satisfied identically for exponents 
satisfying the admissibility conditions.  The latter appear only in 
the equations corresponding to the moves of two bottom lines in Figure 
\ref{f4}.  Each of these moves involves a single triple vertex.  Hence 
the admissibility conditions allow one to eliminate multiplicity and 
weight of one of the edges.  After that the equations are satisfied 
identically by the powers of $q$ with algebraically independent 
exponents.  

Replace in them each $q^{j_i}$ with $t_i$
and restrict ourselves to integer $J_i$'s. Since the 
original equations are satisfied identically by $j_i$'s and $J_i$'s 
satisfying the admissibility conditions, the new equations are satisfied 
identically provided $t_i$'s and $J_i$'s satisfy the identities  
corresponding to the admissibility conditions.

This suggests to eliminate $q$ by replacing its powers with the new 
quantities.  The substitutes for $q^{j_i}$ can be taken from a 
multiplicative group, but must be different from the fourth 
roots of 1, and being as independent as the admissibility condition 
allows.  The group may be not fixed once forever, but chosen each time 
according to the needs.  To keep everything polynomial, the range for 
weights narrows to integers.  We shall follow a more general (but less 
polynomial) solution: take weights from an abelian additive subgroup
$W$ of $B$, 
which is equipped with pairing $M\times W\to M:(t,T)\mapsto t^T$ 
linear in each variable.  For admissibility conditions we need to 
assume that, firstly, $W$ contains the unity $1\in B$,
and, secondly, $t^1=t$ for $t\in M$.

The functors obtained in this way from $\RT^1$
are described from scratch below in the next three sections. 

\subsection{Categories of Colored Framed Graphs}\label{sec2.8} We 
shall call a {\em 1-palette} a quartiple $P$ consisting of a 
commutative associative unitary ring $B$, a subgroup $M$ of the 
multiplicative group of $B$, a subgroup $W$ of the additive group 
$B^+$ of $B$ and a (bilinear) pairing $M\times W\to  M:(m,w)\mapsto 
m^w$. We assume that $W$ contains the unity $1\in B$ and that the 
pairing sends $(m,1)$ to $m$, i.e., $m^1=m$ for any $m\in M$. 

In this section for each 1-palette $P$ we define a category $\cal 
G^1_{P}$.  Its objects are finite sequences of triples each of which 
consists of $t\in M$ with $t^4\ne1$, $T\in W$ and a sign.  The first 
and second elements of the triple are called as in the case of $\cal 
G_h^1$: {\em multiplicity} and {\em weight}.  The triples comprising 
an object are placed at the corresponding points on the line $\R^1$, 
cf.  Section \ref{sec2.1}.  

A morphism of $\cal G_{P}^1$ $$\{(t_1,T_1,\Ge_1),\dots 
(t_k,T_k,\Ge_k)\}\to 
\{(u_1,U_1,\Gd_1),\dots (u_l,U_l,\Gd_l)\}$$ is a framed  generic graph 
$\GG$ embedded in $\R^2\times[0,1]$ with                          
$\p\GG=\Pi_k\times0\cup\Pi_l\times1=\p(\R^2\times[0,1])\cap\GG$ and 
additional structures described below. It is considered up to ambient      
isotopy fixed on $\p\GG$. Its framing at the end points is parallel to the 
$x$-axis. 

The additional structures on $\GG$ all together are called {\em 
$\cal G^1_P$-coloring}. Each of the 1-strata is colored with an
orientation and a 
pair consisting of an element of $\{t\in M\mid t^4\ne1\}$  and an 
element of $W$.
The former is called {\em multiplicity}, the latter {\em weight}. The 
orientations of the edges adjacent to the boundary points are determined 
by the signs of the points' colors as in $\cal G_h^1$: the edge 
adjacent to a point with $+$ is oriented at the point upwards (along the 
standard direction of the $z$-axis), otherwise it is oriented downwards. 
The multiplicity and weight of the edge adjacent to an end point coincide 
with the corresponding ingredients of the color of this point. At each of 
the trivalent vertices the colors of the adjacent edges satisfy to the 
following 

\begin{prop}[Admissibility Conditions]\label{p2.8.A}
Let  $(t_1,T_1)$, $(t_2,T_2)$,
$(t_3,T_3)$ be the multiplicity and weight components of the colors of
three edges whose germs are adjacent to the same vertex and
let $\Ge_i=-1$ if the $i^{\text{th}}$ of the edges is oriented towards 
the vertex and $\Ge_i=1$ otherwise. Then
\begin{equation}\label{t-admiss}
\prod_{i=1}^3t_i^{\Ge_i}=1
\end{equation}
\begin{equation}\label{T-admiss}
\sum_{i=1}^3 \Ge_iT_i=-\prod_{i=1}^3\Ge_i
\end{equation}
\end{prop}

At each strong vertex the graph $\GG$ is oriented. The orientation of 
$\GG$ at strong vertices is considered a part of the $\cal 
G^1_P$-coloring of $\GG$, together with the orientations, 
multiplicities and weights of 1-strata.  

Equation \eqref{t-admiss} means that $\sum t_ie_i$, where the sum is taken 
over all the edges $e_i$ of the graph with their  orientations and $t_i$ 
is the first component of the color of $e_i$, is a relative 1-cycle with 
real coefficients of $\GG$ modulo $\p\GG$, so that it defines an element 
of $H_1(\GG,\p\GG;P)$. This element is called the {\em multiplicity 
homology class} of the coloring. 

\subsection{Functors}\label{sec2.9} Let $P=(B,M,W,M\times W\to  M)$
be a 1-palette (see the preceding section). Denote by $\cal M_{B}$ the 
category of finitely
generated free $B$-modules. Our goal is to construct a functor 
$\A^1:\cal G^1_{P}\to\cal M_{B}$, which will be called the {\em 
Alexander functor.} 

Denote by $\GL$ the $B$-module $B\oplus B$. Denote by $e_0$ 
and $e_1$ the natural basis elements of $\GL$. For an object 
$\{(t_1,T_1,\Ge_1),\dots(t_k,T_k,\Ge_k)\}$ of $\cal G^1_{P}$, the image 
under $\A^1$ is defined to be the tensor product 
$$\GL\otimes_{B}\dots\otimes_{B}\GL$$ with $k$ factors $\GL$. 
This is a free $B$-module of rank $2^k$. It has a canonical basis 
consisting of $e_{i_1}\otimes\dots\otimes e_{i_k}$. So, the basis elements 
may be identified with sequences of zeros and ones placed at the points of 
$\Pi_k$. 

The generating morphisms of $\cal G^1_{P}$ are just colored generators 
of $\cal G$, like in $\cal G_h^1$. Therefore functor $\A^1$ on 
morphisms can also be described by presenting Boltzmann weights, as it was 
done for $\RT^1$ in Tables \ref{tab:BW1} and \ref{tab:BW2}. The Boltzmann 
weights for $\A^1$ are given in Tables \ref{tab:BW3} and \ref{tab:BW4}. In 
these tables (as in \ref{tab:BW1} and \ref{tab:BW2}) a string with 
$e_0$ is shown as dotted arc, and a string with $e_1$ as solid arc. The 
entries of Tables \ref{tab:BW3} and \ref{tab:BW4} are obtained from the 
corresponding entries of Tables \ref{tab:BW1} and \ref{tab:BW2}: if an edge 
colored in a table of Section \ref{sec2.5} with $(j,J)$ has a 
counterpart
in the corresponding table of this section  colored with $(t,T)$ 
then
$q^j$ is replaced with $t$ and $J$ with $T$, as it was
promised 
in Section \ref{sec2.7}. Therefore isotopy invariance of 
state sums based on Boltzmann weights of Tables \ref{tab:BW1} and 
\ref{tab:BW2} implies isotopy invariance of state sums based 
on Tables \ref{tab:BW3} and \ref{tab:BW4}.

If $B=\Z[M]$ and $W=\Z$ then the morphism obtained by $\A^1$ from a 
colored framed graph is represented  by a matrix whose entries are 
Laurent polynomials in multiplicities  attached to the strings of the 
graph. Thus it looks closer to the  Alexander polynomial than 
morphisms provided by $\RT^1$. It is even too  nice, since we are not 
allowed to divide by polynomials, which is  quite desirable, cf. 
\eqref{p5.2.B4}. Therefore we are more interested in the case when $M$ 
is torsion free and $B$ the field $Q(M)$  of quotients of $\Z[M]$.

\begin{table}[pbth]
\centerline{\renewcommand{\arraystretch}{1.6}
\begin{tabular}{|c||c|c|c|c|}
\hline
&${\fig{lb-b.eps}}_{(t,T)}$ & ${\fig{lb-f.eps}}_{(t,T)}$&
${\fig{rb-b.eps}}_{(t,T)}$& ${\fig{rb-f.eps}}_{(t,T)}$ \\  
\hline &  $1$ & $-t^2$&$t^{-2}$& $1$ \\ \hline\hline               
&${\fig{ld-b.eps}}_{(t,T)}$& 
${\fig{ld-f.eps}}_{(t,T)}$&${\fig{rd-b.eps}}_{(t,T)}$& 
${\fig{rd-f.eps}}_{(t,T)}$ \\                  \hline 
&  $1$ & $-t^{-2}$&$t^{2}$& $1$ \\ \hline\hline
&$\fig{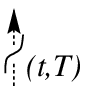}$&$\fig{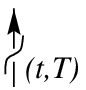}$&$\fig{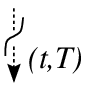}$& 
$\fig{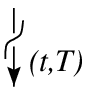}$\\ \hline
&$t^{-T}$& $t^{-T}$&$t^{-T}$&$t^{-T}$\\
\hline\hline                                                        
&$\fig{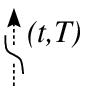}$&$\fig{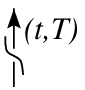}$&$\fig{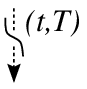}$& 
$\fig{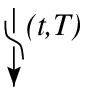}$\\ \hline
&$t^{T}$&$t^{T}$&$t^{T}$&$t^{T}$ \\ \hline\hline               
&${\fig{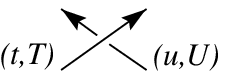}}$& ${\fig{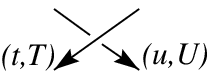}}$&
${\fig{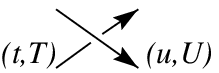}}$&
${\fig{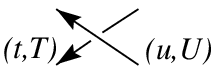}}$\\
\hline 
${\fig{bx0000.eps}}$& $t^{1-U}u^{1-T}$  &
 $t^{1-U}u^{1-T}$ & $t^{1-U}u^{1-T}$ & $t^{1-U}u^{1-T}$  \\ 
 \hline
${\fig{bx0110.eps}}$& $t^{-1-U}u^{1-T}$ & $t^{-1-U}u^{1-T}$
&$t^{1-U}u^{-1-T}$ & $t^{1-U}u^{-1-T}$ \\
 \hline
${\fig{bx1001.eps}}$&  $t^{1-U}u^{-1-T}$& $t^{1-U}u^{-1-T}$&  
$t^{-1-U}u^{1-T}$ & $t^{-1-U}u^{1-T}$ \\                             
  \hline
${\fig{bx1111.eps}}$& $-t^{-1-U}u^{-1-T}$ & $-t^{-1-U}u^{-1-T}$
&$-t^{-1-U}u^{-1-T}$& $-t^{-1-U}u^{-1-T}$ \\                            
\hline
${\fig{bx1010.eps}}$& 0 & 
$(1-t^{-4})t^{1-U}u^{1-T}$&0 &0\\ \hline
${\fig{bx0101.eps}}$& $(t^4-1)t^{-1-U}u^{-1-T}$ &
0 &0&0\\
 \hline
${\fig{bx1100.eps}}$& 0 &0 &$(1-u^{-4})t^{-1-U}u^{1-T}$ &0\\
 \hline
${\fig{bx0011.eps}}$& 0 & 0 &0 &
$(1-u^4)t^{1-U}u^{-1-T}$\\ \hline                           
\hline
&${\fig{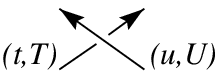}}$&
${\fig{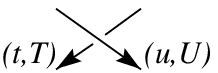}}$&
${\fig{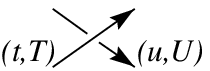}}$&
${\fig{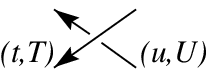}}$\\
\hline
${\fig{bx0000.eps}}$& $t^{-1+U}u^{-1+T}$ &
$t^{-1+U}u^{-1+T}$ & 
$t^{-1+U}u^{-1+T}$&$t^{-1+U}u^{-1+T}$\\                            
 \hline
${\fig{bx0110.eps}}$& $t^{1+U}u^{-1+T}$ &
$t^{1+U}u^{-1+T}$ &$t^{-1+U}u^{1+T}$  &                           
$t^{-1+U}u^{1+T}$ \\                           
 \hline
${\fig{bx1001.eps}}$&  $t^{-1+U}u^{1+T}$ &
$t^{-1+U}u^{1+T}$ &                           
$t^{1+U}u^{-1+T}$  & $t^{1+U}u^{-1+T}$  \\                            
  \hline
${\fig{bx1111.eps}}$& $-t^{1+U}u^{1+T}$  &
$-t^{1+U}u^{1+T}$ &$-t^{1+U}u^{1+T}$ &                            
$-t^{1+U}u^{1+T}$  \\                            
\hline
${\fig{bx1010.eps}}$&  $(1-u^4)t^{-1+U}u^{-1+T}$  &0&0
&0\\ \hline
${\fig{bx0101.eps}}$& 0& $(u^{-4}-1)t^{1+U}u^{1+T}$
 &0&0\\  \hline
${\fig{bx1100.eps}}$& 0 &0&0 &$(1-t^{-4})t^{1+U}u^{-1+T}$ \\
\hline
${\fig{bx0011.eps}}$& 0 & 0 &$(1-t^4)t^{-1+U}u^{1+T}$&0
\\ \hline
\end{tabular}}\vskip.1in
\caption{}\label{tab:BW3}
\end{table}

\begin{table}[tpbh]
\centerline{\renewcommand{\arraystretch}{1.5}
\begin{tabular}{|c||c|c|c|c|}\hline                  
&${\fig{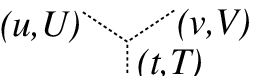}}$ &${\fig{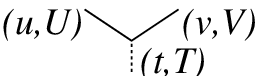}}$
&${\fig{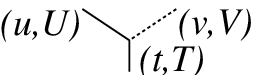}}$ &${\fig{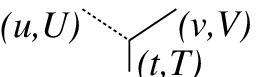}}$ \\
 \hline
${\fig{y001.eps}}$  &
$0$ & $1$ & $-u^{2}$ &  $1$
\\ \hline
${\fig{y001c.eps}}$  &
$0$ & $-1$ & $u^{2}$ &  $-1$
\\ \hline
${\fig{y110.eps}}$  &
$0$ & $t^2-t^{-2}$ & $v^2-v^{-2}$ & $(u^{-2}-u^2)v^{2}$ 
\\ \hline
${\fig{y110c.eps}}$  &
$0$ & $t^{-2}-t^{2}$ & $v^{-2}-v^{2}$ & 
$(u^{2}-u^{-2})v^{2}$ 
\\ \hline
${\fig{y000.eps}}$  &
$t^2-t^{-2}$ & $0$   & $(u^2-u^{-2})v^{-2}$ & $v^2-v^{-2}$ 
\\ \hline
${\fig{y010.eps}}$  &
$1$ & $-u^2$ & $1$ & $0$
\\ \hline
${\fig{y011.eps}}$  &
$v^2-v^{-2}$ & $(u^{-2}-u^2)v^2$ & $0$ & $t^2-t^{-2}$ 
\\ \hline
${\fig{y100.eps}}$  &
$u^{-2}$ & $1$ & $0$ & $1$
\\ \hline
${\fig{y101.eps}}$  &
$(u^2-u^{-2})v^{-2}$ & $v^2-v^{-2}$ & 
$t^2-t^{-2}$ & $0$ 
\\ \hline
${\fig{y111.eps}}$  &
$1$ & $0$ & $1$ & $u^{-2}$
\\ 
\hline\hline
&${\fig{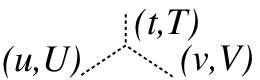}}$ &${\fig{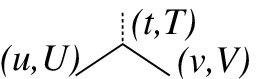}}$
&${\fig{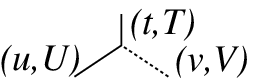}}$ &${\fig{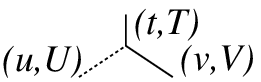}}$ \\
 \hline
${\fig{y-001.eps}}$  &
$0$ & $t^2-t^{-2}$ & $v^2-v^{-2}$ & 
$(u^{-2}-u^2)v^{-2}$ 
\\ \hline
${\fig{y-001c.eps}}$  &
$0$ & $t^{-2}-t^2$ & $v^{-2}-v^2$ & 
$(u^2-u^{-2})v^{-2}$ 
\\ \hline
${\fig{y-110.eps}}$  &
$0$ & $1$ & $-u^{-2}$ & $1$
\\ \hline
${\fig{y-110c.eps}}$  &
$0$ & $-1$ & $u^{-2}$ & $-1$
\\ \hline
${\fig{y-000.eps}}$  &
$1$ & $0$ & $1$  & $u^2$
\\ \hline
${\fig{y-010.eps}}$  &
$(u^2-u^{-2})v^2$ & $v^2-v^{-2}$ & 
$t^2-t^{-2}$ & $0$ 
\\ \hline
${\fig{y-011.eps}}$  &
$u^2$ & $1$ & $0$ & $1$
\\ \hline
${\fig{y-100.eps}}$  &
$v^2-v^{-2}$ & $(u^{-2}-u^2)v^{-2}$ & $0$ & 
$t^2-t^{-2}$ 
\\ \hline
${\fig{y-101.eps}}$  &
$1$ & $-u^{-2}$ & $1$ & $0$
\\ \hline
${\fig{y-111.eps}}$  &
$t^2-t^{-2}$ & $0$ & $(u^2-u^{-2})v^2$ & $v^2-v^{-2}$ \\ 
\hline \end{tabular}}\vskip.1in
\caption{}\label{tab:BW4}
\end{table}

\subsection{Still Representations of Hopf Algebra}\label{sec2.10} In 
definition of $\RT^1$ given in Sections \ref{sec2.3}, \ref{sec2.5}, the
action  of $U_qgl(1|1)$ is completely hidden behind the Boltzmann
weights, and at  first glance is not needed.  However if we forgot it 
completely, the arguments based on irreducibility of representations, like 
the ones used in the proofs of \ref{p5.1.A}, \ref{p5.1.B} and
\ref{p5.1.C} would be impossible.

The target of $\A^1$ is not yet equipped with a structure which
would be a counterpart of $U_qgl(1|1)$-representations. In fact, they 
cannot be equipped with an action of $U_qgl(1|1)$. However, $U_qgl(1|1)$ 
contains a Hopf-subalgebra $U^1$ (see Section \ref{sL.10}) such that
$\A^1$ can be  upgraded to a functor to the category of modules over
this subalgebra.

Now we can redefine the Alexander functor $\cal A^1:\cal G^1_P\to\cal 
M_{B}$ introduced in Section \ref{sec2.9} as a functor from the same 
category $\cal G^1_P$, but targeted at the category of finite 
dimensional modules over $U^1\otimes_{\Z}B$ which assigns to an object 
$\{(t_1,T_1,\Ge_1)$, \dots, $(t_k,T_k,\Ge_k)\}$ of $\cal G^1_{P}$ the 
tensor product 
$$U(t_1,T_1)_{\Ge_1}\otimes_{B}\dots\otimes_{B}U(t_k,T_k)_{\Ge_k},
$$
See Section \ref{sL.11}
On morphisms the functor is defined by the same Boltzmann weights as
above.  

\section{Reshetikhin-Turaev Functor Based on $sl(2)$}\label{sec3}

In this section we construct a functor $\RT^2$  similar to the functor
$\RT^1$  defined in Section \ref{sec2}. It is a modified special  case of
a functor introduced by Deguchi and Akutsu in Section 5 of \cite{DA}.
The functor  constructed in \cite{DA} depends on a natural number $N$.
Our modification  is related to the case $N=2$. 

The choice of scaling of Clebsch-Gordan morphisms used by Deguchi and 
Akutsu \cite{DA} is based on orthogonality relations.  Its 
disadvantage is that at each vertex one of the adjacent edges is 
distinguished and has to be directed upwards.  

Here I use another scaling. It does not satisfy the orthogonality 
relations, but eliminate the choice of an edge at each trivalent 
vertex. Together with a deferent choice of basis vectors in the 
representations this makes formulas much simpler than in \cite{DA}. 

\subsection{Category of Colored Framed Graphs}\label{sec3.1} The source 
category for $\RT^2$ is denoted $\cal G^2$ and defined as follows. 

An object of $\cal G^2$ is a class of finite sequences of pairs. 
Each of the pairs consists of a complex number and a sign. The complex 
numbers are not allowed to be odd integers. In a pair the complex number 
may be multiplied by $-1$ with simultaneous change of the sign component 
of the pair. An object of $\cal G^2$ consists of all the sequences 
which can be obtained from each other by a sequence of operations like 
this.

A morphism of $\cal G^2$ 
$$\{(A_1,\Ge_1),\dots,(A_k,\Ge_k)\}\to\{(B_1,\Gd_1),\dots,(B_l,\Gd_l)\}$$ 
is a morphism $\GG:k\to l$ of category $\cal G$ (see Section \ref{sec1.7}) 
equipped with the additional structure described below.   

Each of 1-strata is colored with an orientation and a complex number. The 
latter is called a {\em weight}. It is prohibited to be an odd 
integer. The orientation of a 1-stratum
may be changed with a simultaneous multiplication of the weight 
by $-1$. In other words, this is an element $C$ of 
$H_1(\GG,\{vertices\};\C)$, that is a 1-chain of $\GG$ with complex 
coefficients. This chain must satisfy the following 

\begin{prop}[Admissibility Condition]\label{p3.1.A} 
The boundary $\p C$ of $C$ must involve each trivalent vertex of $\GG$ 
with coefficient $\pm 1$,  boundary vertex $(i,0)$, $i=1,\dots k$ with 
coefficient $-\Ge_iA_i$ and boundary vertex $(i,1)$, $i=1,\dots l$ with
coefficient $\Gd_iB_i$, respectively. 
\end{prop}

To the properties of the weight chain $C$, we should add that 
each 1-stratum is involved in $C$ with a coefficient, which cannot be 
an odd integer.

An internal vertex is called a {\em source} if it is involved in $\p 
C$ with coefficient $-1$, otherwise it is called a {\em sink}. In 
figures a sink is shown in with a small light disk.

\subsection{Constructing Reshetikhin-Turaev Functor}\label{sec3.2} Our 
goal is to construct a functor $\RT^2$ from $\cal G^2$ to the category 
of finite dimensional representations of Hopf algebra  $U_{\sqrt{-1}}\ 
sl(2)$.  For a brief summary on these quantum algebra and its 
representations used below, see Appendix 2.  The image of an object 
$\{(A_1,\Ge_1),\dots,(A_k,\Ge_k)\}$ of $\cal G^2$ is defined as the 
tensor product $I(\Ge_1A_1)\otimes\dots\otimes I(\Ge_kA_k)$, see 
Section \ref{A2.4}.  

Each of the factors $I(\Ge A)$ is a 2-dimensional vector space with a 
canonical basis $e_0$, $e_1$.  Thus, $I(\Ge_1A_1)\otimes\dots\otimes 
I(\Ge_kA_k)$ has a canonical basis consisting of 
$e_{i_1}\otimes\dots\otimes e_{i_k}$. The basis vectors of 
$I(\Ge_1A_1)\otimes\dots\otimes I(\Ge_kA_k)$  may be identified with 
sequences of zeros and ones. These zeros and ones are associated with 
orientation of the corresponding points.

For presentation of the images of morphisms via Boltzmann weights, one 
has to consider diagrams of graphs representing the morphisms of $\cal 
G^2$ and coloring of their 1-strata with basis vectors of the 
representations associated with the colors of the 1-strata.  It is 
convenient to identify these second colors with orientations of 
1-strata: if a 1-strata is oriented upwards, it is colored with $e_0$, 
otherwise it is colored with $e_1$.  The same orientations can be used 
for presenting of the (underlying) $\cal G^2$-coloring, i.~e., coloring 
with pairs (orientation, weight).  Recall that the orientation 
in this pair can be reversed with simultaneous multiplication of the 
weight by $-1$.  Therefore, one may use any orientation to 
specify the $\cal G^2$-color, and it is natural to use the same 
orientation for both purposes: describing the representation attached 
to the string and a basis vector chosen in the representation.  The 
weight is shown at the string.  The same picture can be used 
just for $\cal G^2$-colored graph.  To emphasize the double usage of 
orientations, we use light arrowheads.  

Hence an edge with upwards light arrowhead and number $A\in\C$ and an 
edge with downwards light arrowhead and number $-A$ denote a string 
colored with the same primary color $(\text{upwards orientation}, A)$, 
which is associated with the same representation $I(A)$.  The secondary 
colors are vectors $e_0$ and $e_1$, respectively, in the first and 
second case.  The same edge with the same orientations and 
multiplicities, but shown with dark arrowheads is colored with the 
same primary color (associated with $I(A)$), but with no secondary 
color, i.~e., with no basis vector of $I(A)$ specified.

The image of a morphism is defined via Boltzmann weights shown in 
Tables \ref{tab:BW1-sl2} and \ref{tab:BW2-sl2}.  In these tables and 
below $\im$ denotes $\sqrt{-1}$, while $\im^x$ denotes
$\exp(\frac{x\pi\im}2)$.  The Boltzmann weights are found on 
the basis of the Reshetikhin-Turaev Theorem \cite{RT1}, similarly to 
the Boltzmann weights of Tables \ref{tab:BW1} and \ref{tab:BW2}.  The 
relevant information about the Hopf algebra $U_{\im}sl(2)$ and its 
representations is placed in Appendix 2.  One can check directly that 
the Boltzmann weights of Tables \ref{tab:BW1-sl2} and 
\ref{tab:BW2-sl2} define a functor from $\cal G^2$ to the category of 
modules over $U_{\im}sl(2)$.

\begin{table}[hbtp]\centerline{\renewcommand{\arraystretch}{1.7}
\begin{tabular}{|c|c|c|c|}\hline
${\fig{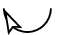}}_{A}$ & 
${\fig{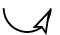}}_{A}$&
${\fig{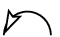}}_{A}$&
${\fig{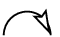}}_{A}$ \\                 
\hline   $1$ & $\im^{-A-1}$&$1$&$\im^{A+1}$ \\ \hline\hline               
${\fig{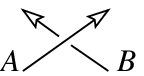}}$& ${\fig{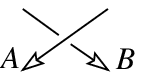}}$&
${\fig{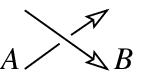}}$&
${\fig{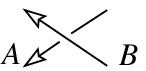}}$\\
\hline                                       
  $\im^{\frac{(1+A)(1+B)}2}$&
 $\im^{\frac{(1+A)(1+B)}2}$&$\im^{\frac{(1+A)(1+B)}2}$&
 $\im^{\frac{(1+A)(1+B)}2}$ \\
\hline\hline
${\fig{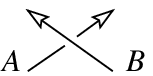}}$&
${\fig{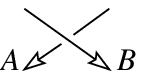}}$&
${\fig{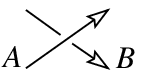}}$&
${\fig{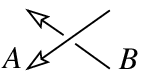}}$\\
\hline
 $\im^{\frac{-(1+A)(1+B)}2}$&
$\im^{\frac{-(1+A)(1+B)}2}$&
$\im^{\frac{-(1+A)(1+B)}2}$&
$\im^{\frac{-(1+A)(1+B)}2}$\\
\hline
\hline
\multicolumn{2}{|c}{$\fig{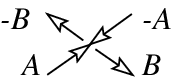}$}&
\multicolumn{2}{|c|}{$\fig{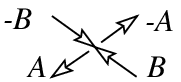}$}\\ \hline
\multicolumn{2}{|c}{$\im^{\frac{-(1-A)(1-B)}2}(\im^{1+A}+\im^{1-A})$} 
& 
\multicolumn{2}{|c|}{$\im^{\frac{(1+A)(1+B)}2}(\im^{-1+B}+\im^{-1-B})$} 
\\ 
\hline
\hline 
$\fig{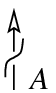}$&$\fig{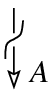}$&$\fig{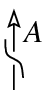}$ &
$\fig{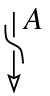}$\\ \hline
$\im^{\frac{A^2-1}4}$&$\im^{\frac{A^2-1}4}$&$\im^{\frac{1-A^2}4}$&
$\im^{\frac{1-A^2}4}$ \\ \hline
\end{tabular}}\hskip.3cm 
\caption{}\label{tab:BW1-sl2}
\end{table}

\begin{table}[hbtp]
\centerline{\renewcommand{\arraystretch}{1.2}
\begin{tabular}{|c|c|c|c|}\hline
  & 
${}^A\raisebox{-3ex}{\includegraphics{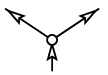}}\ {}^B$&
 \placeY{A}{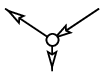}B &  \placeY{A}{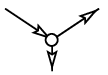}B\\ 
sink vertex
&$\scriptstyle{{\phantom{aaaaaa}}C=A+B+1}$&$\scriptstyle{{\phantom{aaaaaaaa}}
C=-A+B-1}$&
$\scriptstyle{{\phantom{aaaaaa}}C=A-B-1}$\\ 
\cline{2-4} 
&$\im^{1+C}+\im^{1-C}$&$\im^{1+B}+\im^{1-B}$&$\im^{A-B}+\im^{-A-B}$ \\ 
\hline\hline &\placeY{A}{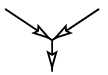}B & \placeY{A}{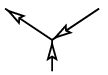}B &
\placeY{A}{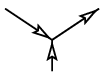}B \\ 
source vertex
&$\scriptstyle{{\phantom{aaaaaa}}C=A+B+1}$&$\scriptstyle{{\phantom{aaaaaaa}}
C=A-B-1}$&
$\scriptstyle{{\phantom{aaaaaaa}}C=-A+B-1}$\\ 
\cline{2-4} 
&$1$&$1$&$\im^{-1-A}$\\ \hline\hline
&$\scriptstyle{{\phantom{aaaaaa}}C=A+B+1}$&
$\scriptstyle{{\phantom{aaaaaaaa}}C=-A+B-1}$&
$\scriptstyle{{\phantom{aaaaaa}}C=A-B-1}$\\
sink vertex&
\placeA{A}{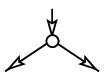}B 
&\placeA{A}{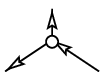}B&\placeA{A}{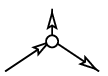}B \\
\cline{2-4}
& $\im^{1+C}+\im^{1-C}$ &$\im^{1+B}+\im^{1-B}$ 
&$-\im^{B+A}-\im^{B-A}$ \\ \hline\hline
&$\scriptstyle{{\phantom{aaaaaa}}C=A+B+1}$
&$\scriptstyle{{\phantom{aaaaaa}}
C=A-B-1}$& $\scriptstyle{{\phantom{aaaaaaa}}C=-A+B-1}$\\
source vertex&
\placeA{A}{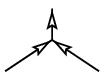}B 
&\placeA{A}{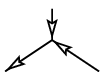}B&\placeA{A}{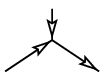}B \\
\cline{2-4}
& $1$ &$1$ 
&$\im^{1+A}$ \\ \hline
\end{tabular}}\hskip.3cm 
\caption{}\label{tab:BW2-sl2}
\end{table}

\subsection{Homological Meaning of Colorings}\label{sec3.3} 
Presentation of a secondary coloring of a generic graph by 
orientations of 1-strata suggests a homology interpretation of the 
secondary coloring as an integer 1-chain, in which any 1-stratum has 
coefficient $1$, if it is taken with the orientation specified by 
the light arrow. Denote this chain by $Z$.

An easy analysis of Table \ref{tab:BW2-sl2} shows that $Z$ has a 
curious property: it completes the weight chain $C$ of the 
corresponding primary coloring to a cycle modulo the boundary of the 
underlying graph.  In other words, $\p(C+Z)$ involves no internal 
vertices.

This suggests to associate to a secondary color (i.~e. a basis vector of
$I(A)$) the coefficient with which it appears in $C+Z$. For a basis
vector $e_i$ of the representation $I(A)$,  the complex number
$a=(-1)^iA+1$ is called the {\em charge.\/} A secondary $\cal G^2$-coloring
of a generic graph $\GG$ can be presented by a 1-chain $c=C+Z$ with complex
coefficients on $\GG$: each 1-stratum of $\GG$ is equipped with the
orientation associated with the secondary color and the charge.
Hence a secondary coloring defines an element of $H_1(\GG,\p\GG;\C)$.
However, it is not defined by this homology class: one has to add the
orientations of 1-strata. These orientations are encoded in $Z$.

Replacing the weights by charges makes Tables
\ref{tab:BW1-sl2} and \ref{tab:BW2-sl2} simpler. See Table
\ref{tab:BW-sl2c}. The difference between sink and source vertices
almost vanishes. Most formulas expressing the Boltzmann weights become
simpler. This can be explained by the algebraic nature of charge: it is
an Eigenvalue of $K$.

The homology class of a secondary coloring reduced modulo 2 depends
only on the primary coloring. It belongs to $H_1(\GG;\C/2\Z)$ and can
be described directly in terms of the primary coloring by adding 1 to
each weight and reducing modulo 2.

\begin{table}[hbtp]
\renewcommand{\arraystretch}{1.7}
\begin{tabular}{|c|c|c|c|}\hline
${\fig{lb-d.eps}}_{a}$ & 
${\fig{rb-d.eps}}_{a}$&
${\fig{ld-d.eps}}_{a}$&
${\fig{rd-d.eps}}_{a}$ \\                 
\hline   $1$ & $\im^{-a}$&$1$&$\im^{a}$ \\ \hline\hline               
${\fig{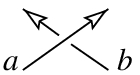}}$& ${\fig{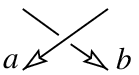}}$&
${\fig{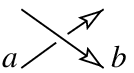}}$&
${\fig{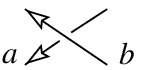}}$\\
\hline                                       
  $\im^{\frac{ab}2}$&
 $\im^{\frac{ab}2}$&$\im^{\frac{ab}2}$&
 $\im^{\frac{ab}2}$ \\
\hline\hline
${\fig{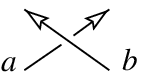}}$&
${\fig{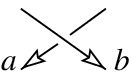}}$&
${\fig{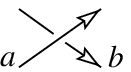}}$&
${\fig{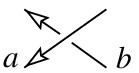}}$\\
\hline
 $\im^{\frac{-ab}2}$&
$\im^{\frac{-ab}2}$&
$\im^{\frac{-ab}2}$&
$\im^{\frac{-ab}2}$\\
\hline
\hline
\multicolumn{2}{|c}{$\fig{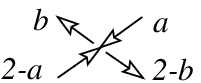}$}&
\multicolumn{2}{|c|}{$\fig{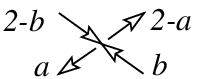}$}\\ \hline
\multicolumn{2}{|c}{$\im^{\frac{-ab}2}(\im^a-\im^{-a})$}&
\multicolumn{2}{|c|}{$\im^{\frac{ab}2}(\im^{-b}-\im^{b})$}\\ 
\hline
\hline 
$\fig{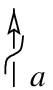}$&$\fig{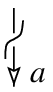}$&$\fig{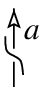}$ &
$\fig{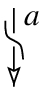}$\\ \hline
$\im^{\frac{a(a-2)}4}$&$\im^{\frac{a(a-2)}4}$&$\im^{\frac{a(2-a)}4}$&
$\im^{\frac{a(2-a)}4}$ \\ \hline
\end{tabular}
\renewcommand{\arraystretch}{1.23}
\begin{tabular}{|c|c|c|}\hline
${}^a\raisebox{-3ex}{\includegraphics{fig/y000d.eps}}\ {}^b$&
 \placeY{a}{y011d.eps}b &  \placeY{a}{y101d.eps}b\\ 
$\scriptstyle{{\phantom{aaaa}}c=a+b}$&$\scriptstyle{{\phantom{aaaaaa}}
c=-a+b}$&
$\scriptstyle{{\phantom{aaaa}}c=a-b}$\\ 
\hline 
$\im^c-\im^{-c}$&$\im^b-\im^{-b}$&$\im^{-b}(\im^a-\im^{-a})
$ \\ \hline\hline
 \placeY{a}{y111d.eps}b & \placeY{a}{y010d.eps}b &
\placeY{a}{y100d.eps}b \\ 
$\scriptstyle{{\phantom{aaaa}}c=a+b}$&$\scriptstyle{{\phantom{aaaaa}}
c=a-b}$&
$\scriptstyle{{\phantom{aaaaa}}c=-a+b}$\\ 
\hline 
$1$&$1$&$\im^{-a}$\\ \hline\hline
$\scriptstyle{{\phantom{aaaa}}c=a+b}$&
$\scriptstyle{{\phantom{aaaaaa}}c=-a+b}$&
$\scriptstyle{{\phantom{aaaa}}c=a-b}$\\
\placeA{a}{y-111d.eps}b 
&\placeA{a}{y-100d.eps}b&\placeA{a}{y-010d.eps}b \\
\hline
 $\im^c-\im^{-c}$ &$\im^b-\im^{-b}$ 
&$\im^{b}(\im^a-\im^{-a})$ \\ \hline\hline
$\scriptstyle{{\phantom{aaaa}}c=a+b}$
&$\scriptstyle{{\phantom{aaaa}}
c=a-b}$& $\scriptstyle{{\phantom{aaaaa}}c=-a+b}$\\
\placeA{a}{y-000d.eps}b 
&\placeA{a}{y-101d.eps}b&\placeA{a}{y-011d.eps}b \\
\hline
 $1$ &$1$ 
&$\im^{a}$ \\ \hline
\end{tabular}
\caption{}\label{tab:BW-sl2c}
\end{table}

\begin{figure}[htb]
\centerline{\includegraphics{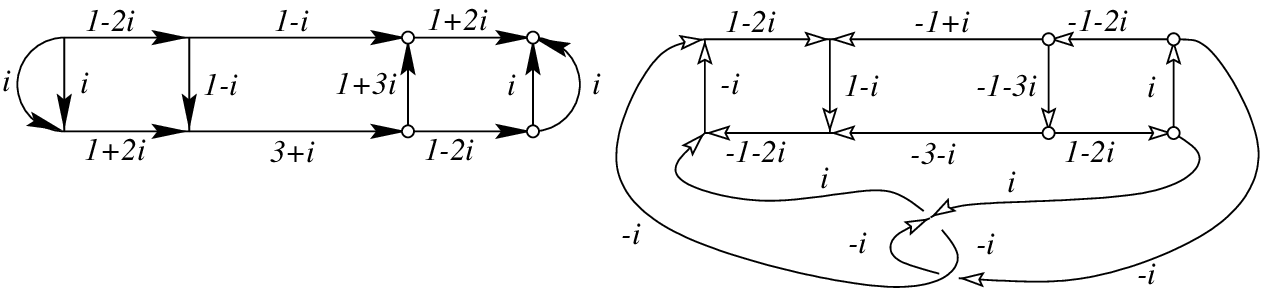}}
\caption{}
\label{fnosecond}
\end{figure}

It may happen that a diagram of a generic graph with a primary 
$\cal G^2$-coloring does not admit any secondary coloring. In Figure 
\ref{fnosecond} such a closed generic graph with a primary 
$\cal G^2$-coloring is shown on the left hand side. This property depends 
on embedding of the graph. On the right hand side of Figure 
\ref{fnosecond}	the same graph, but embedded in a different way is 
shown with a secondary coloring (to keep evident the relation to the primary
coloring, the secondary coloring is presented not by the charges, but by
the weights).

Of course, absence of secondary colorings implies that $\RT^2$ maps 
the (primarily) colored	framed graph to the zero morphism.

\subsection{Algebraization}\label{sec3.4} In this section we do with 
$\RT^2$ what is done with $\RT^1$ in Sections \ref{sec2.7} -- 
\ref{sec2.10}: we eliminate formal power series.  As in the case of 
$\RT^1$, the most obvious obstacle for this is  presence of products 
in exponents.  In the case of $\RT^1$ all the products in exponents 
are products of a multiplicity and a weight.  In formulas: Boltzmann 
weights involve terms of the form $q^{jI}$, where $j$ is the 
multiplicity component of a color and $I$ is the weight components of, 
maybe other, color.  

In Table \ref{tab:BW1-sl2} the quadratic part of exponents is made of 
values of a single parameter, weight, which are complex numbers 
prohibited to be odd integers.  Thus the solution that is used in the 
case of $\RT^1$ cannot be applied here literally.  To rewrite the 
Boltzmann weights of tables \ref{tab:BW1-sl2} and \ref{tab:BW2-sl2} in 
a formal algebraic way without exponential function, one may follow 
the approach indicated below.  

The complex parameters are substituted with parameters taken from a 
commutative ring $W$, while the Boltzmann weights are taken from a 
larger ring $B$. Besides the inclusion, $W$ and $B$ are related by a 
group homomorphism of the additive group $W^+$ of $W$ to the 
multiplicative group $B^{\times}$ of $B$. 

In the original setup of $\RT^2$, both $W$ and $B$ are $\C$ and the 
homomorphism $W^+\to B^{\times}$ is defined by $A\to\im^A$. 

In our generalization, we use the same exponential notation for the 
homomorphism $W^+\to B^{\times}$. This is justified by the condition 
that $4\in W$ must be mapped by the homomorphism to $1\in B$, i.~e. 
$\im^4=1$. This condition implies that $B$ contains a root of unity of 
degree 4: this is $\im^1$, indeed $(\im^1)^4=\im^{(4\cdot1)}=\im^4=1$. 
This element of $B$ is denoted by $\im$. 

Let us call a {\em 2-palette} a triple
$P$ consisting of $W$, $B$ and the map
$W^+\to B^{\times}:A\to\im^A$ as above: \begin{itemize}
\item $B$ is a commutative ring with unity,
\item $W$ is a subring containing the unity,
\item $W^+\to B^{\times}:A\to\im^A$ is a homomorphism with $\im^4=1$.
\end{itemize}

The category $\cal G^2_P$ of colored framed graphs replacing $\cal 
G^2$ looks as follows. Its object is a class of finite sequences of 
pairs. Each of the pairs consists of an element of $W$ and a sign. The 
former is called the {\em weight} of the color. It is not allowed to 
belong to $\{x\in W\mid \im^{2x+2}=1\}$. In a pair, the weight may 
be multiplied by $-1$ with simultaneous change of the sign. An object
is a  class of all the sequences which can be obtained from each other
by  operations like this applied to the pairs comprising the sequence.

A morphism of $\cal G^2_P$ 
$$\{(A_1,\Ge_1),\dots,(A_k,\Ge_k)\}\to\{(B_1,\Gd_1),\dots,(B_l,\Gd_l)\}$$ 
is a morphism $\GG:k\to l$ of category $\cal G$ equipped with an
additional  structure called {\em $\cal G^2_P$-coloring.} 

Each 1-stratum is colored  with an orientation and an element of $W$ 
called
the {\em  weight} of the  color. It may not belong to $\{x\in W\mid
\im^{2x+2}=1\}$. The orientation  of a 1-stratum may be changed with
simultaneous multiplication of the weight by $-1$. In other words,
this is an element $C$ of  $H_1(\GG,\{vertices\};W)$, that is a 
1-chain of
$\GG$ with coefficients  in $W$. This chain must satisfy the following

\begin{prop}[Admissibility  Condition]\label{p3.4.A} The boundary $\p 
C$ must involve each internal (i.~e., trivalent)  vertex of $\GG$ with 
coefficient $\pm 1$, boundary vertex $(i,0)$, $i=1,\dots k$ with 
coefficient $-\Ge_iA_i$ and boundary vertex $(i,1)$, $i=1,\dots l$ 
with coefficient $\Gd_iB_i$, respectively. 
\end{prop}

The generalization $\A^2$ of $\RT^2$ acts from $\cal G^2_P$ 
to a the category of finitely generated free modules over 
$B$. It assigns to the
class of sequences $\{(\pm A_1,\pm)$,\dots, $(\pm A_k,\pm)\}$ the vector
space of dimension $2^k$ generated by all the representatives of the 
class. This is the tensor product (over $B$) of $k$ copies of a free  
$B$-module of rank 2. 

It is not just a $B$-module. It admits a natural action of a Hopf 
algebra $U^2$ which can be defined exactly as $U_{\im}sl(2)$ (see
Section \ref{A2.3}), but
without relation $K=\im^H$. Namely, $U^2$ is the Hopf algebra 
generated by $H$, $K$, $X$ and $Y$ satisfying the relations:
\begin{equation}\label{relM}\begin{gathered}
\[Y,X\]=K-K^{-1},\\
[H,X]=-2X,\quad [H,Y]=2Y,\\ 
KX=-XK,\quad KY=-YK.\end{gathered}\end{equation}
with co-product defined by 
\begin{equation}\label{co-prod-M}\begin{aligned}\GD:
U^2\to&U^2\otimes U^2:
\\
&\GD(H)=H\otimes 1+1\otimes H\\ 
&\GD(K)=K\otimes K\\
&\GD(K^{-1})=K^{-1}\otimes K^{-1}\\
&\GD(X)=X\otimes K^{-1}+1\otimes X,\\
&\GD(Y)=Y\otimes 1+K\otimes Y,
\end{aligned}\end{equation}
co-unit defined by
\begin{equation}\label{co-unit-M}\begin{aligned}\Ge:
U^2\to&\C[[h]]:\\
&\Ge(K)=\Ge(K^{-1})=1,\quad\Ge(H)=\Ge(X)=\Ge(X)=0,\end{aligned}
\end{equation} and antipode
\begin{equation}\label{antipode-M}\begin{aligned}s:
U^2\to&U^2:\\
&s(H)=-H,\\
&s(K)=K^{-1},\quad s(K^{-1})=K,\\
&s(X)=-XK\quad(=KX),\\
&s(Y)=-K^{-1}Y\quad(=YK^{-1}).
\end{aligned}\end{equation}

In the case of the free $B$-module of rank 2 corresponding to 
$\{(\pm A,\pm)\}$, the generators $H$, $K$, $X$ and $Y$ act by the
following formulas:
\begin{equation}\label{XYKaction}\begin{aligned}
X:\ &(A,+)\mapsto (\im^{A+1}+\im^{-A+1})(-A,-),\qquad (-A,-)\mapsto 
0\\ Y:\ &(A,+)\mapsto0,\qquad (-A,-)\mapsto(A,+)\\ 
K:\ &(A,+)\mapsto \im^{A+1}(A,+)\qquad (-A,-)\mapsto \im^{A-1}(-A,-)\\
H:\ &(A,+)\mapsto (A+1)(A,+)\qquad (-A,-)\mapsto (A-1)(-A,-)
\end{aligned}\end{equation}
Denote this $U^2\otimes_{\Z}B$-module by $I_B(A)$.

In the case of a longer sequence $\{(A_1,\Ge_1),\dots,(A_k,\Ge_k)\}$,
the corresponding $B$-module is considered as  the tensor product of
the sequence of $U^2\otimes_{\Z}B$-modules of rank 2 corresponding  to
the elements of the sequence and the action is defined by via the
coproduct of  $U^2$. 

A morphism $\GG$ is mapped to the morphism of the  corresponding
modules defined by the matrix whose entries are constructed  as
follows. An entry of the matrix corresponds to the choice of signs at 
the boundary points of $\GG$. Consider all the extensions of these
signs  to orientations of $1$-strata of a diagram of $\GG$. As in the
case  of $\RT^2$, the extensions  are called secondary colorings. In a 
diagram the orientations comprising a secondary coloring are shown by
light  arrowheads, as in the case of $\RT^2$. For each of the
extensions, take the product of Boltzmann weights at the vertices and 
extremal points of the height function in the diagram according to the 
tables \ref{tab:BW1-sl2} and  \ref{tab:BW2-sl2}  and sum up the 
product over all the extensions. If a local picture at a vertex of the
diagram  does not appear in the tables, assign 0 to the vertex, i.~e.
disregard the  extension.

Of course, $\RT^2$ is a special case of $\A^2$ with $W=B=\C$.

\section{Relations Between the Reshetikhin-Turaev 
Functors}\label{sec4}

Although at first glance $\cal G_h^1$ and $\cal G^2$ seem to be quite 
different, more careful analysis shows deep similarities and relations 
between them.

\subsection{Things to Recolor}\label{sec4.1}There is 
a subcategory $S\cal G^1$ of $\cal G^1_{\pi\im/2}$ and a 
natural transformation of $\RT^1|_{S\cal G^1}$ to $\RT^2$.

An object $\{(j_1,J_1,\Ge_1),\dots,(j_k,J_k,\Ge_k)\}$ of $\cal 
G^1_{\pi\im/2}$ belongs to $S\cal G^1$ if $j_r+J_r=1$ for 
$r=1,\dots,k$. A morphism $\GG$ of $\cal G^1_{\pi\im/2}$ belongs to 
$S\cal G^1$ if for each 1-stratum of $\GG$ the sum of its multiplicity 
and weight is $1$.  

This choice is motivated as follows. We want to have a correspondence 
between $\cal G^1$- and $\cal G^2$-colorings which would preserve the 
Boltzmann weights. Comparison of the first lines of Tables 
\ref{tab:BW1} and \ref{tab:BW1-sl2} suggests to assume $q=\im$ and 
$-2j=-A-1$, i.~e. $j=\frac{A+1}2$. The Boltzmann weights at half-twists
suggest $-jJ=\frac{A^2-1}4$, which together with $j=\frac{A+1}2$	
imply $J=\frac{1-A}2$ and $j+J=1$. Of course, this does not pretend to
be a proof, but  rather a strong indication in favor of the choice made
above.

\subsection{Recoloring}\label{sec4.2}
Define a functor $\GF:S\cal G^1\to \cal G^2$.  An object 
$\{(j_1,J_1,\Ge_1)$,\dots, $(j_k,J_k,\Ge_k)\}$ of $S\cal G^1$	is turned
by $\GF$ to $\{(2j_1-1,\Ge_1),\dots,(2j_k-1,\Ge_k)\}$. To a morphism 
$\GG$ of
$S\cal G^1$ it assigns the isotopy class of the same oriented framed 
graphs with	weights on 1-strata obtained from the original 
multiplicities by the same formula: if the  multiplicity (in 
$\cal G^1$-coloring)	is $j$, the weight (in $\cal G^2$-coloring) 
is $2j-1$.

\begin{prop}[Consistency in Admissibility Conditions]\label{p4.2.A}	
The recoloring converts  any morphism of $S\cal G^1$	to a morphism of
$\cal  G^2$. In other words, $\GF$	turns a coloring of a generic graph
satisfying Admissibility Conditions \ref{2.2.A}	and the condition 
defining $S\cal G^1$ (the sum of multiplicity and weight for each 
1-stratum is $1$) into a coloring satisfying Admissibility Condition 
\ref{p3.1.A}. \end{prop}

\begin{prop}[Lemma]\label{p4.2.B} Any morphism of $S\cal G^1$ has no 
strong vertices.  \end{prop}

\begin{proof}Recall that a vertex is said to be strong if the 
adjacent edges are oriented either all towards it or in the 
opposite direction. In notations of Admissibility Conditions 
\ref{2.2.A}, this means that $\Ge_1=\Ge_2=\Ge_3$. By Admissibility
Conditions \eqref{j-admiss} 
and \eqref{J-admiss}, at a strong vertex $j_1+j_2+j_3=0$ and 
$J_1+J_2+J_3=1$. Hence $j_1+j_2+j_3+J_1+J_2+J_3=1$.
On the other hand, $j_r+J_r=1$ in $S\cal G^1$, hence 
$j_1+j_2+j_3+J_1+J_2+J_3=3$. \end{proof}

\begin{proof}[Proof of \ref{p4.2.A}] By Lemma at each internal vertex 
at least one of the adjacent edges is directed towards the vertex and 
at least one outwards. Therefore $j_1+j_2-j_3=0$, where $j_1$, $j_2$ 
and $j_3$ are the multiplicities of the edges adjacent to a vertex 
and numerated appropriately. Then 
$|(2j_1-1)+(2j_2-1)-(2j_3-1)|=1$, and hence the vertex is involved 
with coefficient $\pm1$ into the boundary of the weight chain 
of the coloring provided by $\GF$.
\end{proof}

\subsection{Invariance of Boltzmann Weights Under 
Recoloring}\label{sec4.3}
Categories $\cal G^1_h$, $S\cal G^1$ and $\cal G^2$ can be enriched by 
incorporating secondary colorings into objects and morphisms. The 
enriched 
categories are denoted by  $\cal G^{1\flat}_{h}$, $S\cal G^{1\flat}$ 
and $\cal G^{2\flat}$, respectively.

\begin{prop}[Consistency in Boltzmann Weights]\label{p4.3.A} Functor 
$\GF$ can be enhanced to a functor $\GF^{\flat}$ from $S\cal 
G^{1\flat}$ to $\cal G^{2\flat}$ in such a way that the Boltzmann 
weight of $\RT^1$ equals of the Boltzmann weight of $\RT^2$ at the 
same point after application $\GF^{\flat}$.  \end{prop}

\begin{proof}The construction looks as follows: the dotted arcs (i.~e., 
colored with bosons) get light arrowhead with the same orientation, 
while the solid arcs (colored with fermions) change orientation and 
$\cal G^2$-weight.  For example, $$\fig{y-111ex.eps} \text{ 
turns into }\fig{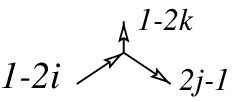}.$$ Notice that the $\RT^1$-Boltzmann 
weight at $\fig{y-111ex.eps}$ is $(q^{2i}-q^{-2i})q^{2j}$, while the 
$\RT^2$-Boltzmann weight at $\fig{y-111exd.eps}$ is 
$-\im^{2j-1}(\im^{1-2i}+\im^{2i-1})=\im^{2j}(\im^{2i}-\im^{-2i})$.  
Thus the $\RT^2$-Boltzmann weight after the transformation coincides
with the $\RT^1$-Boltzmann weight of the original graph.
One can easily check this for each of the 
Boltzmann weights of Tables \ref{tab:BW1} and \ref{tab:BW2}, except for
the weights at strong vertices, which do not appear in $S\cal G^1$ by 
Lemma \ref{p4.2.A}. Here is another example. This is one of the most 
complicated entries from Table \ref{tab:BW1}:
$\fig{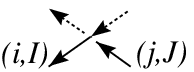} \text{ turns into }\fig{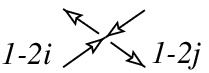}.$
(Here we present a secondary $\cal G^{2\flat}$-coloring by weights.)
The $\RT^1$-Boltzmann 
weight at $\fig{bx-rrex.eps}$ is 
$$q^{iJ+jI+i-j}(1-q^{-4i})=q^{i(1-j)+j(1-i)+i-j}(1-q^{-4i})=
q^{-2ij}(q^{2i}-q^{-2i}),$$  while the $\RT^2$-Boltzmann weight at 
$\fig{bx-rrexd.eps}$ is

\begin{multline*}
\im^{\frac12(-(1-2i)(1-2j)+(1-2i)+(1-2j)+1)}(\im^{1-2i}+\im^{2i-1})=
\im^{\frac12(-4ij+2)}\im(\im^{-2i}-\im^{2i})=\\ \im^{-2ij}(\im^{2i}-\im^{-
2i}).\end{multline*} 

Other entries of Tables \ref{tab:BW1} and \ref{tab:BW2} can be treated 
similarly, but the number of necessary checkings can be reduced, if 
one takes into account the behavior of the Boltzmann weights under the 
moves of Figure \ref{f4}.  \end{proof}

\subsection{Functorial Reformulation}\label{sec4.4}
Denote by $Vect(\C)$ the category of complex vector 
spaces. Denote by $(\RT^1)^{\sharp}$	and $(\RT^2)^{\sharp}$ the 
compositions of $\RT^1$	and $\RT^2$, respectively, with the 
appropriate forgetting functor to $Vect(\C)$.

\begin{prop}[Corollary of \ref{p4.3.A}]\label{p4.4.A} The functor 
$\GF$
which maps the subcategory $S\cal G^1$ of $\cal G^1_{\pi\im/2}$ to 
$\cal G^2$ comprises together with the identity functor of the category 
$Vect(\C)$  a natural transformation of 
$(\RT^1)^{\sharp}|_{S\cal G^1}$ to $(\RT^2)^{\sharp}$. In other words,
$$\begin{CD}S\cal G^1 @>{(\RT^1)^{\sharp}|_{S\cal G^1}}>>Vect(\C)\\ 
@VV{\GF}V	@|\\
			  \cal G^2 @>{(\RT^2)^{\sharp}}>>Vect(\C)
  \end{CD}$$ 
is a commutative diagram.
In particular, for any morphism $\GG$ of $S\cal G^1$
$$(\RT^2)^{\sharp}\GF(\GG)=(\RT^1)^{\sharp}(\GG).\qed$$
\end{prop}

Probably, this is a manifestation of a deeper relation between 
$U_{\im}gl(1|1)$ and $U_{\im}sl(2)$ and one can get rid of forgetting 
functors to $Vect(\C)$, cf. \cite{MR-P}. However, we do not elaborate 
this, since it does not seem to have immediate topological 
corollaries.

\subsection{Generalization}\label{sec4.5}
In this section we generalize the results of Sections \ref{sec4.1} -
\ref{sec4.4} replacing Reshetikhin-Turaev functors $\RT^1$ and $\RT^2$
by their generalizations $\A^1$ and $\A^2$. 

Let $P$ be a 1-palette (see Section \ref{sec2.8}, where we defined 
also $\cal G^1_P$-colored graphs). Recall that $P$ consists of a 
commutative ring
$B$, a subgroup $M$ of the multiplicative group of $B$, a
subgroup $W$ of the additive group $B^+$ of $B$ and a (bilinear)
pairing $M\times W\to M:(m,w)\mapsto m^w$. We assume that $W$ contains
the unity $1$ of $B$ and  that the pairing takes value $m$ on $(m,1)$.

Furthermore, suppose \begin{itemize}
\item $W$ is a subring of $B$, 
\item the pairing $M\times W\to M:(m,w)\mapsto m^w$ is an action of 
the multiplicative monoid of $W$ and
\item $M$ contains an element $\im $ such that $\im^4=1$, and
$\im^2=-1\ne1$. 
\end{itemize}
The second assumption means that $(m^u)^v=m^{uv}$ for any $u,v\in W$.

Denote by $Q$ the triple consisting of the rings $B$, $W$ and the map 
$W\to B:w\mapsto \im^w$. The latter is a homomorphism of the additive 
group of $W$ to the multiplicative group of $B$, since $M\times W\to 
M:(m,w)\mapsto m^w$ is bilinear. Hence $Q$ is a 2-palette (see Section 
\ref{sec3.4}) and one can consider the category $\cal G^2_Q$  of $\cal 
G^2_Q$-colored graphs. 

This construction of $Q$ can be reversed. Given {\em any\/} 2-palette 
$Q=(B,\ W,\ W^+\to B^{\times}:A\mapsto\im^A)$ 
(i.e., $Q$ consists of a commutative unitary ring $B$, its unitary 
subring $W$ and a homomorphism $W^+\to B^{\times}:A\mapsto\im^A$ such 
that $\im^4=1$), one can construct the quartiple $P$ consisting of 
\begin{itemize}
\item the same ring $B$, 
\item the image $M$ of the homomorphism $W^+\to B^{\times}$ 
defined by formula $A\mapsto\im^A$, 
\item the subring $W\subset B$ and
\item the pairing $M\times W\to M:(\im^A,B)\mapsto\im^{AB}$. 
\end{itemize}
It is clear that this $P$ is a 1-palette and satisfies 
all the conditions imposed on $P$ at the beginning of 
this section. Moreover, if one applies to this $P$ the construction 
described above, the result is the initial 2-palette $Q$.

Let $P$	be as at the beginning of Section \ref{sec4.5}.
Define a subcategory $S_{\im}\cal G^1_P$ of $\cal G^1_P$. An object 
$\{(t_1,T_1,\Ge_1)$, \dots, $(t_k,T_k,\Ge_k)\}$ belongs to 
$S_{\im}\cal G^1_P$ if $\im^{1-T_r}=t_r$ for $r=1,\dots,k$. A morphism 
$\GG$ of $\cal G^1_P$ belongs to $S_{\im}\cal G^1_P$ if for each 
1-stratum of $\GG$ its multiplicity $t$ and weight $T$ are related in 
the same way: $\im^{1-T}=t$. 

Now define a functor $\GF_P:S_{\im}\cal G^1_P\to\cal G^2_Q$.
It turns an object  $\{(t_1,T_1,\Ge_1)$,\dots, 
$(t_k,T_k,\Ge_k)\}$ of $S\cal G^1_P$ to 
$\{(1-2T_1,\Ge_1),\dots,(1-2T_k,\Ge_k)\}$. To a morphism  $\GG$ of 
$S\cal G^1_P$ it assigns the isotopy class of the same oriented framed 
graphs with charges on 1-strata obtained from the original  weights by 
the same formula: if the original weight (in  $\cal G^1_P$-coloring) 
is $T$, the new weight (in $\cal G^2_Q$-coloring) is $1-2T$.

The functor $\GF_P:S_{\im}\cal G^1_P\to\cal G^2_Q$ 
generalizes in the obvious sense the functor $\GF:S\cal G^1\to\cal 
G^2$ defined in Section \ref{sec4.2}: if $B=\C$, $W=\Z$ 
$M=\C\sminus0$, and $M\times W\to M$ is defined by $(m,w)\mapsto m^w$ 
then $S_{\im}\cal G^1_P=S\cal G^1$, $\cal G^2_Q=\cal G^2$ and 
$\GF_P=\GF$.

\begin{prop}[Consistency in Admissibility Conditions]\label{p4.5.A} 
$\GF_P$ transforms a coloring of a generic graph satisfying 
Admissibility Conditions \ref{p2.8.A} and the conditions defining 
$S_{\im}\cal G^1_P$  into a coloring satisfying Admissibility 
Condition \ref{p3.4.A}. \end{prop} 

\begin{prop}[Lemma]\label{p4.5.B} Any morphism of $S\cal G^1_P$ has no 
strong vertices.  \end{prop} 

\begin{proof}Recall that a vertex is said to be strong if the adjacent 
edges are oriented either all towards it or in the opposite direction. 
In notations of Admissibility Conditions \ref{p2.8.A}, this means that 
$\Ge_1=\Ge_2=\Ge_3$. By \eqref{t-admiss} and \eqref{T-admiss}, at a 
strong vertex $t_1t_2t_3=1$ and $T_1+T_2+T_3=1$. Belonging to 
$S_{\im}\cal G^1_P$ means that $t_i=\im^{1-T_i}$. Hence 
$t_1t_2t_3=\im^{3-(T_1+T_2+T_3)}=\im^2\ne 1$, which contradicts to 
$t_1t_2t_3=1$.  \end{proof} 

\begin{proof}[Proof of \ref{p4.5.A}] By Lemma at each internal vertex 
at least one of the adjacent edges is directed towards the vertex and 
at least one outwards. Therefore $|T_1+T_2-T_3|=1$, where $T_1$, $T_2$ 
and $T_3$ are the weights of the edges adjacent to a vertex numerated 
appropriately. Then $|(1-2T_1)+(1-2T_2)-(1-2T_3)|= |1-2|=1$ and  hence 
the vertex is involved with coefficient $\pm1$ into the boundary of 
the weight chain of the coloring provided by $\GF$. \end{proof} 

Denote by $\cal G^{1\flat}_P$, $S_{\im}\cal G^{1\flat}_P$ and $\cal 
G^{2\flat}_Q$ the categories $\cal G^{1}_P$, $S_{\im}\cal G^{1}_P$ and 
$\cal G^{2}_Q$, respectively, enhanced by incorporating secondary 
colorings to objects and morphisms. 

\begin{prop}[Consistency in Boltzmann Weights]\label{p4.5.C} Functor 
$\GF_P$ can be enhanced to a functor $\GF^{\flat}_P$ from $S_{\im}\cal 
G^{1\flat}_P$ to $\cal G^{2\flat}_Q$ in such a way that the Boltzmann 
weight of $\A^1$ equals of the Boltzmann weight of $\A^2$ at the same 
point after application $\GF^{\flat}_P$.  \end{prop} 

This is proved like Theorem \ref{p4.3.A}. Functor $\GF^{\flat}_P$ 
turns each dotted arc colored with $(\im^{1-T},T)$ into a solid one 
colored with $1-2T$ without changing orientation, 
but the arrowhead becomes light. Each solid arc colored with 
$(\im^{1-T},T)$ changes its orientation and receives weight $2T-1$. 
The proof is completed by a
straightforward comparing the Boltzmann weights.\qed

\begin{prop}[Corollary of \ref{p4.5.C}]\label{p4.5.D} The functor 
$\GF_P$ which maps the subcategory $S_{\im}\cal G^1_P$ of $\cal G^1_P$ 
to $\cal G^2_Q$ comprises together with the identity functor of the 
category $\cal M_B$  a natural transformation of $\A^1|_{S\cal G^1}$ 
to $\A^2$. In other words, $$\begin{CD}S_{\im}\cal G^1_P 
@>{\A^1|_{S\cal G^1}}>>\cal M_B\\ @VV{\GF_P}V @|\\ \cal G^2_Q 
@>{\A^2}>>\cal M_B \end{CD}$$ is a commutative diagram. In particular, 
for any morphism $\GG$ of $S_{\im}\cal G^1_P$ 
$$\A^2\GF_P(\GG)=\A^1(\GG).\qed$$ \end{prop}

\subsection{Non-Surjectivity and Surjectivity of 
Recoloring}\label{sec4.6} As it was observed in Section 
\ref{sec3.3}, there exist $\cal G^2$-colorings which do not admit any 
secondary coloring and therefore give rise to morphisms of $\cal G^2$ 
mapped by $\RT^2$ to the zero morphisms. Obviously, a morphism of 
$\cal G^2$ which does not admit	a secondary coloring, does not belong 
to the image of $\GF$. However, morphisms of this kind are not 
interesting, especially if we compare $\RT^2$ and $\RT^1$.

\begin{prop}[Partial Reversing of $\GF^{\flat}_P$]\label{p4.6.A} Any 
morphism of $\cal G^{2\flat}_P$ with a diagram having neither fragment 
$\fig{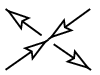}$, nor $\fig{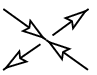}$ belongs to the image of 
$\GF^{\flat}_P$.  \end{prop} 

\begin{proof}Graphically the construction which converts a morphism of 
this kind can be described as follows: it replaces solid arcs by 
dotted ones, light arrowheads by solid ones; each weight $A$ 
by the pair made of multiplicity $\im^{\frac{A+1}2}$ and weight 
$\frac{1-A}2$.  This 
construction is inverse to the construction which proves \ref{p4.5.C}.  
One can check that fragments of Table \ref{tab:BW2-sl2}	are turned to 
fragments of the first column of Table \ref{tab:BW4}, with the
multiplicities of the adjacent edges satisfying Admissibility 
Conditions \ref{p2.8.A}. Similarly, all entries of Table 
\ref{tab:BW1-sl2}, besides $\fig{bxpb.eps}$ and $\fig{bxnb.eps}$, are 
turned into the corresponding entries 
of Table \ref{tab:BW3}.  \end{proof}

The condition of absence of $\fig{bxpb.eps}$, nor $\fig{bxnb.eps}$ 
cannot be skipped. One can prove that the right hand side 
graph of Figure \ref{fnosecond}	does not belong to the image of 
$\GF^{\flat}$. Indeed, it is impossible to find appropriate weights.
Existence of weights is equivalent to existence of secondary coloring
for the graph on the left hand side of Figure \ref{fnosecond}.

\begin{prop}[Reversing $\GF^{\flat}_P$ on Graphs Without Internal
Vertices]\label{p4.5.E} A morphism of $\cal G^{2\flat}_Q$ with a
diagram having no internal vertices belongs to the image of
$\GF^{\flat}_P$. The $\cal G^{1\flat}_P$-coloring is uniquely defined 
by the underlying orientations of components.\end{prop}

\begin{proof} Each arc of the diagram with weight $A$, where the 
orientation of the secondary color and orientation of the component 
coincide, has to be replaced by a dotted arc with the same 
orientation, multiplicity $\im^{\frac{A+1}2}$ and weight 
$\frac{1-A}2$. Each arc with weight $A$ and the orientations opposite 
to each other has to be equipped with the orientation of the 
component, multiplicity $\im^{\frac{1-A}2}$ and weight $\frac{A+1}2$. 
An easy verification shows that each entry of Table \ref{tab:BW1-sl2} 
transforms to an entry of Table \ref{tab:BW3}. \end{proof}

\subsection{Conclusion}\label{sec4.7} The results of this section 
show that $\RT^1$ and $\RT^2$ are closely related, but do not reduce 
completely to each other.  The source category $\cal G^1_h$ seems to 
be wider in several ways than the source category $\cal G^2$ of 
$\RT^2$.  First, colors have more components.  Second, admissibility 
conditions allow more kinds of local pictures.  Orientations at 
vertices in $\cal G^1$ have no counterpart in $\cal G^2$.  At last, 
the restriction of $\RT^1$ to a comparatively small part of $\cal 
G^1_{\pi\im/2}$ can be factored through a map onto an essential part 
of $\cal G^2$. 

On the other hand, $\cal G^2$ and $\RT^2$ are essentially simpler than 
$\cal G^1$ and $\RT^1$.  I have not been able to find a complete 
reduction of the formers to the latters, and doubt if such a 
reduction exists. A simpler theory has obvious advantages. Thus we 
need to consider both.

Generalizations $\A^1$ and $\A^2$ of $\RT^1$ and $\RT^2$ are also
related to each other. However the conditions on the initial
algebraic data under which the relating natural transformation
exists are quite hard. This provides an additional reason to treat
both theories.

\section{Skein Principles and Relations}\label{sec5}

\subsection{Skein Principles}\label{sec5.1}
Since $\A^1$ and $\A^2$
map a colored point to an irreducible representation, the Schur Lemma 
implies the following.

\begin{prop}[Colors' Simplicity]\label{p5.1.A} Let $c$ be 1 or 2 and 
$\GG$ be a morphism of the category $\cal G^c_P$ with 
$\p\GG=\Pi_1\times0\cup\Pi_1\times1$. Assume that $B$ is a field.
Then $\A^c(\GG)$ is a 
multiplication by a constant. The constant is zero, unless the end points 
of $\GG$ are colored with the same color.\qed\end{prop} 

\begin{prop}[Restricted Colors' Completeness]\label{p5.1.B} Let $c=1$ 
or $2$ and $\GG$ be a morphism of $\cal G^c_P$ with 
$\p\GG=\Pi_2\times0\cup\Pi_2\times1$.  Let the coloring of the end 
points $\GG$ can be extended to a coloring of the graph 
${\fig{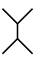}}$ satisfying the admissibility conditions.
 Assume that $B$ is a field. Then $\A^c(\GG)$ is a linear 
combination of morphisms which are the images under $\A^c$ of the 
graph ${\fig{i-gr-sml.eps}}$ equipped with colorings extending the 
same coloring of the end points and satisfying the admissibility 
conditions.  \end{prop} 

In fact, there are at most two terms in the linear combination of
\ref{p5.1.B}.

\begin{proof}[Proof of \ref{p5.1.B}] Consider the case $c=1$. The 
case $c=2$ differs only in notations. The morphism under consideration is 
a morphism $$(u_0,U_0)_{\Ge_0}\otimes(v_0,V_0)_{\Gd_0}\to
(u_1,U_1)_{\Ge_1}\otimes(v_1,V_1)_{\Gd_1}.$$
The assumptions on the colors imply that each of the tensor products is a
direct sum of two irreducible modules over the algebra $U^1$,
see Appendix 1, Section \ref{sL.11}. The
morphism under consideration must preserve these decompositions. Therefore
it is a linear combination of at most two isomorphisms between the
corresponding summands. (If there is no isomorphic summands in the
decompositions then the original morphism is trivial.) An isomorphism
between the summands of the decompositions can be realized as the image of
the graph ${\fig{i-gr-sml.eps}}$ colored in the appropriate way.
\end{proof}

The restriction on colors of $\p\GG$ in \ref{p5.1.B} makes a profound
difference from the case of quantum invariants studied in \cite{RT1} and 
\cite{RT2}.  However, it is
necessary. Indeed, by Admissibility Conditions \ref{p2.8.A}, 
${\fig{i-gr-sml.eps}}$ does not admit any coloring such that the
bottom two end points are colored with $(t,T_1,+)$ and $(t^{-1},T_2,+)$.

It is impossible to enlarge the palette to eliminate the restrictions
from \ref{p5.1.B} without sacrificing \ref{p5.1.A}. Indeed, the
tensor product of modules $(t,T_1)_+$ and $(t^{-1},T_2)_+$ is not a direct
sum of irreducible $U^1$-modules.

\begin{prop}[Killing by Splittable Closed Part]\label{p5.1.C}
Let $c=1$ or $2$ and $\GG\subset\R^2\times [0,1]$ be a colored framed 
graph representing a morphism of $\cal G^c_P$ and containing 
a closed nonempty subgraph $\GS$ which is separated from the rest of $\GG$ 
by a 2-sphere disjoint with $\GG$ and embedded in $\R^2\times [0,1]$. Then 
$\A^c(\GG)=0$.
\end{prop}

\begin{proof} First, consider the case $\GS=\GG$. Then $\GG$ is an 
endomorphism of the empty sequence. Since $\GS\ne\empt$, it can be 
decomposed as $\rd\circ\GF\circ\lb$ where $\GF$ is an endomorphism of \
 $(t,T)_+\otimes(t,T)_-$ \ or \ $I(A)\otimes I(-A)$ for $c=1$ and $c=2$ 
respectively. The first operator in 
this composition, $\lb$, maps $B$ to the invariant irreducible 
submodule which can be described as $\Ker X\cap\Ker Y$ in both cases. 
Whatever $\GF$ is, it takes this subspace to itself (since 
$\GF$ commutes with $X$ and $Y$). The operator $\rd$ annihilates this 
subspace. Therefore the whole composition vanishes.

If $\GG\ne \GS$, the splitting provides representation 
$\GG=\GS\otimes(\GG\sminus\GS)$. Thus 
$\A^c(\GG)=\A^c(\GS\otimes(\GG\sminus\GS))=
\A^c(\GS)\otimes\A^c(\GG\sminus\GS) 
=0\otimes\A^c(\GG\sminus\GS)=0$. \end{proof}

While Theorems \ref{p5.1.A} and \ref{p5.1.B} have counterparts in cases of 
other Reshetikhin-Turaev functors, Theorem \ref{p5.1.C} is more special. 
It implies other specific features of theories related to the Alexander 
polynomial.

\subsection{First Skein Relations}\label{sec5.2} Although 
\ref{p5.1.A} and \ref{p5.1.B} do not look like skein relations, they 
obviously imply existence of skein relations.  Of course, these 
relations can be proven independently by a straightforward 
calculation.  

For example, Theorem \ref{p5.1.A} implies that
$\A^1\(\fig{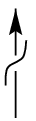}_{(t,T)}\)$ is equal to $\A^1\(\
\fig{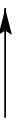}_{(t,T)}\)$ multiplied by a constant. The constant is
obvious from the construction of $\A^1$: it is $t^{-T}$. Hence, we 
have relation
\begin{equation}\label{A-r-tw-skein}
\A^1\(\ {\fig{r-htw-up.eps}_{(t,T)}}\)=t^{-T}\A^1\(\
{\fig{vert-up.eps}_{(t,T)}}\). \end{equation}
This can be used to replace any graph by the same graph with the
black-board framing.

\begin{prop}[Removing Loop Relation for $\A^1$]\label{p5.2.A}
\begin{equation}\label{Aphi-gr+-skein}
\A^1\({\fig{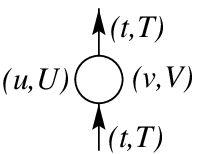}}\)=(t^{2}-t^{-2}) \A^1 \(\
{\fig{vert-up.eps}_{(t,T)}}\)
\end{equation}
for any admissible colorings of the arcs in the middle
of the left hand side such that at strong vertices all the 
orientations are of the same sign. \end{prop}

Proof of \ref{p5.2.A} is a straightforward calculation. It is
facilitated by the following observation. It suffices to calculate both
sides with any color at the both end points, since the morphism is scalar.
However a choice of the color may facilitate calculations. For example, if 
the arcs 
in the middle of the left hand side of \eqref{Aphi-gr+-skein} are oriented
upwards, then the left hand side admits only one coloring with bosonic
colors at the end points, and two colorings with fermionic colors. In the
bosonic case the Boltzmann weight at the upper vertex is $1$ and at the
bottom $t^{2}-t^{-2}$, see Table \ref{tab:BW4}.
Similar calculations prove \eqref{Aphi-gr+-skein} for other               
orientations of these arcs. If the arcs constitute an oriented cycle, the
choice of the fermionic color at the end points reduces the state sum to a
single summand, otherwise the bosonic color does this.
\qed

\begin{prop}[Junction Relations for $\A^1$]\label{p5.2.B} If $u^4v^4\ne1$, 
then \begin{multline}\label{p5.2.B1}
\A^1\(\fig{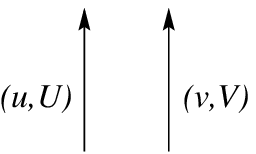}\)=\\
\frac1{u^2v^2-u^{-2}v^{-2}}\A^1\(\fig{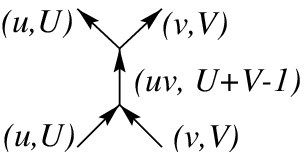} 
\)+\\
\frac1{u^{-2}v^{-2}-u^2v^2}\A^1\(\fig{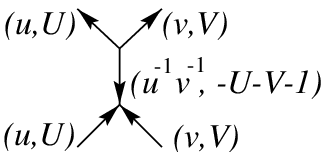}\) \end{multline}
If $u^4v^{-4}\ne1$, then
\begin{multline}\label{p5.2.B2}
\A^1\(\fig{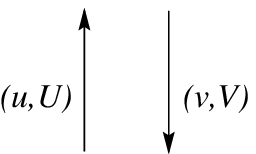}\)=\\
\frac1{u^{2}v^{-2}-u^{-2}v^{2}}\A^1\(\fig{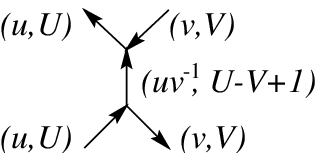} 
\)+\\ \frac1{u^{-2}v^2-u^2v^{-2}}
\A^1\(\fig{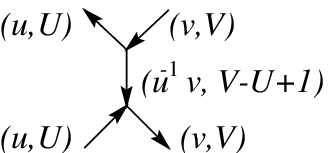}\) \end{multline}            
\end{prop}

This is proved by a direct calculation. \qed

Relations \eqref{p5.2.B1} and \eqref{p5.2.B2} are so similar that it
would be natural to unite them. The only tool we need for this is the
following

{\bf\it Agreements on Hiding Inessentials.} {\it
\begin{itemize}
\item If the orientations of some edges are not shown in a diagram,
the edges are oriented, but in a non-specified way.
\item If several diagrams are involved in the same
equality and have end points, the edges adjacent to the equally
positioned end points are colored (and, in particular, oriented) in the 
same way in all the diagrams.
\item If only the multiplicity of a color is shown, the other ingredients 
of the color are ignored because they are not important. One may recover 
them in any way respecting the admissibility conditions. (However, missing 
orientations at strong vertices are counter-clockwise) 
\item If a formula involves a sum over colors of an edge without specified 
orientation, the sum runs over all the colors (and, in particular, 
orientations) making admissible triples at the end points of the edge with 
the colors of the adjacent edges. \end{itemize} } 

Under this agreement \eqref{p5.2.B1} and \eqref{p5.2.B2} unite 
as follows:
\begin{equation}\label{p5.2.B4}
\A^1\(\fig{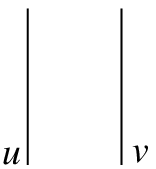}\)=\sum_{(w,W)}
\frac1{w^{2}-w^{-2}}\A^1\(\fig{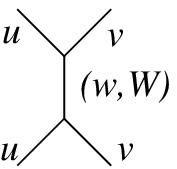}\) \end{equation}
Relation \eqref{Aphi-gr+-skein},
which is already written partially under the agreement, 
can be rewritten in a more concise form:
\begin{equation}\label{Aphi-gr-skein}
\A^1\({\fig{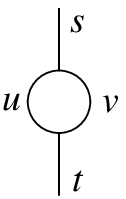}}\)=\Gd_{s,t}(t^{2}-t^{-2}) \A^1 \(\
{\fig{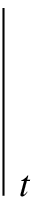}}\)
\end{equation}

\subsection{First Skein Relations for the Second 
Functor}\label{sec5.3}
In a similar way, but simpler, one can prove the following counterparts of 
\eqref{A-r-tw-skein}, \ref{p5.2.A} and \ref{p5.2.B} for $\A^2$:
\begin{equation}\label{2A-r-tw-skein}
\A^2\(\ {\fig{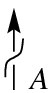}}\)=\im^{\frac{A^2-1}4}\A^2\(\
{\fig{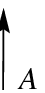}}\) \end{equation}

\begin{prop}[Removing Loop Relation for $\A^2$]\label{p5.3.A}
\begin{equation}\label{Aphi-gru-skein}
\A^2\({\fig{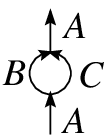}}\)=(\im^{1+A}+\im^{1-A}) \A^2 \(\
{\fig{vert-A.eps}}\)
\end{equation}
for any admissible colorings of the arcs in the middle
of the left hand side. \qed
\end{prop}

\begin{prop}[Junction Relations for $\A^2$]\label{p5.3.B} If 
$A+B\not\in 2\Z$, then 
\begin{multline}\label{p5.3.B2}
\A^2\(\fig{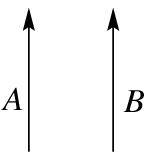}\)=\\
\frac1{\im^{A+B}-\im^{-A-B}}\A^2\(\fig{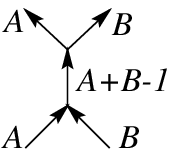} 
\)+\\
\frac1{-\im^{A+B}+\im^{-A-B}}\A^2\(\fig{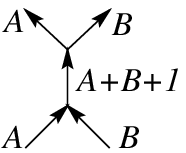}\) 
\end{multline} \qed \end{prop}

\section{Alexander Invariant of Closed Colored Framed 
Graphs}\label{sec6}

\subsection{Definition of Conway Function via Reshetikhin-Turaev 
Functor}\label{sec6.1} An isotopy class of a $\cal G^1$-colored framed link 
can be considered as a morphism of $\cal G^1_h$ from the empty set 
to itself. The tensor product of the empty family of modules is naturally 
identified with the ground field $\C$, since this is the unity for the 
tensor multiplication of modules. Thus the isotopy class of a colored 
framed link defines, via the Reshetikhin-Turaev functor, a complex number 
depending on $q$. In case of other polynomial link invariant, such as the 
Jones polynomial, this is a way to relate the polynomial and the functor: 
the analogous function of $q$ defined by the functor is the corresponding 
polynomial. 

However, as follows from Theorem \ref{p5.1.C}, in the case of $U_qgl(1|1)$ 
this does not work: the number constructed in this way is identically 
zero. This is well-known, see \cite{KS}, \cite{RS1}. The standard way 
to obtain the Conway
function from the Reshetikhin-Turaev functor involves an auxiliary 
geometric construction. One of the strings in the link should be pulled 
and cut, to convert the link into a tangle connecting $\Pi_1\times0$ 
with
$\Pi_1\times1$ in $\R^2\times[0,1]$. By \ref{p5.1.A}, the corresponding 
morphism is multiplication by a number depending on $q$ and the colors of 
the components. This number divided by $q^{2j}-q^{-2j}$, where $j$ is the 
multiplicity of the cut string's color, is the value of the Conway 
function of the original link evaluated at $q^{j_1}$, \dots $q^{j_n}$ 
where $j_1$, \dots, $j_n$ are the multiplicities of the link components' 
colors, provided the weights of all the colors are zero. See Rozansky and 
Saleur \cite{RS1}, \cite{RS2}, Murakami \cite{Mu1}, Deguchi 
and Akutsu \cite{DA}.

\subsection{Why It Does Not Matter Where To Cut}\label{sec6.2}
Here I show that the Conway functor defines in a similar
way an invariant for $\cal G^1_P$-colored framed graphs in $R^3$. 
The proof is based on the following Lemma.

\begin{figure}[htb]
\centerline{\includegraphics{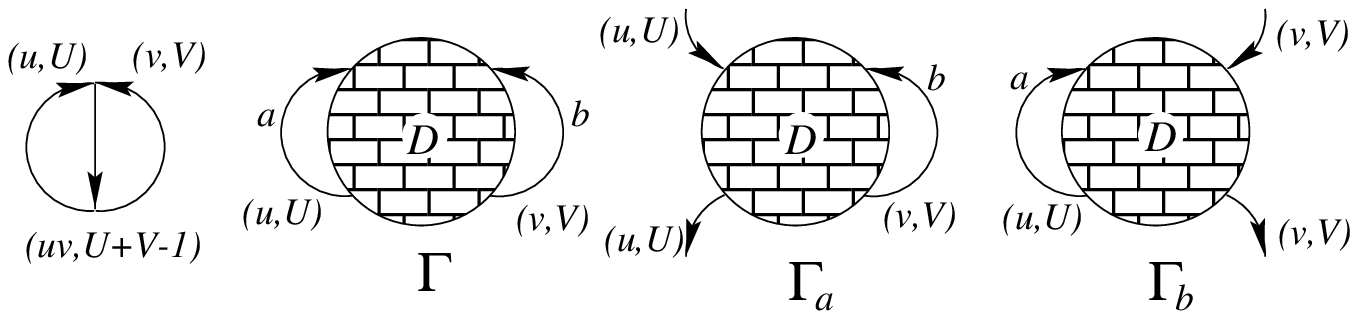}}
\caption{}
\label{f6.2.A}
\end{figure}

\begin{prop}[Lemma]\label{p6.2.A} Let  $P=(B,M,W,M\times W\to M)$ be a 
1-palette (see Section \ref{sec2.8}). Assume that $B$ is a field. Let 
$\GG$ be a $\cal G^1_P$-colored framed
closed generic graph embedded in $\R^3$ such that its projection to $\R^2$
is generic. Let $D\subset \R^2$ be a disk containing all the
projection of $\GG$ except two arcs, $a$ and $b$, colored with $(u,U)$ and
$(v,V)$, respectively. Let $\GG_a$ and $\GG_b$ be colored framed
graphs
which are obtained from $\GG$ by cutting at $a$ and $b$, respectively, and
moving the end points upwards and downwards to a position such that one of 
them becomes the uppermost and the other one lowermost. While moving, the 
framing at the end points is kept black-board. See Figure
\ref{f6.2.A}. Then  $\A^1(\GG_a)$
is multiplication by $(u^2-u^{-2})\GD$ and $\A^1(\GG_b)$
is multiplication by $(v^2-v^{-2})\GD$ with the same $\GD\in B$. 
\end{prop}

\begin{proof} First, assume that these arcs are colored
in such a way that the coloring can be extended
to an  admissible coloring of a planar theta graph containing $a$ and
$b$.
 By \ref{p5.1.B}, the common part of all three graphs $\GG$,
$\GG_a$ and $\GG_b$, that is the part whose projection is covered by $D$,
is mapped by $\A^1$ to a linear combination of the images of two
differently colored copies of ${\fig{i-gr-sml.eps}}$. This linear
combination looks as follows:
$\Ga\A^1\({\fig{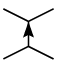}}\)+
\Gb\A^1\({\fig{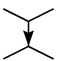}}\)$.

Then
$$\A^1(\GG_a)=\Ga\A^1\({\fig{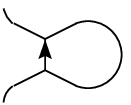}}\) +
\Gb\A^1\({\fig{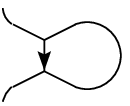}}\)$$
By \eqref{Aphi-gr+-skein}, $\A^1\({\fig{lf-g-smu.eps}}\)=
\A^1\({\fig{lf-g-smd.eps}}\)=
(u^2-u^{-2})\A^1\({\fig{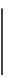}}\)$ and hence     
$$\A^1(\GG_a)=(\Ga+\Gb)(u^2-u^{-2}) \A^1 \({\fig{vert.eps}}\).$$
Thus $\A^1(\GG_a)$ is multiplication by $(\Ga+\Gb)(u^2-u^{-2})$.

Similarly, $\A^1(\GG_b)$ is multiplication by 
$(\Ga+\Gb)(v^2-v^{-2})$.
This proves the statement: take  $\GD=\Ga+\Gb$.

The assumption about the extensibility of the colors over a theta 
graph was used at the very beginning of the proof. There are 
situations when it does not hold true. For example, with the 
orientations as in Figure \ref{f6.2.A}, it can happen that $u=v^{-1}$. 
Then \ref{p5.1.B} cannot be applied and the fragment of $\GG$ covered 
by $D$ cannot be replaced by the linear combination. However in this 
case one can first make a small kink on one of the arcs, say $b$, and 
then expand $D$ to hide under $D$ the new crossing point. See Figure 
\ref{f6.2.A2}. Outside $D$ we see the arc $b$ on the previous place, 
but oppositely oriented. Now the assumption holds true, so we can 
apply the arguments above. 

\begin{figure}[htb]
\centerline{\includegraphics{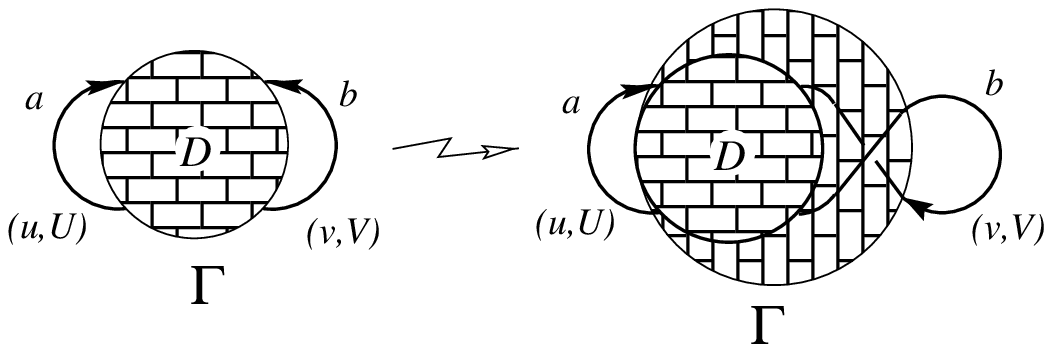}}
\caption{}
\label{f6.2.A2}
\end{figure}
 The result is not exactly what we desire: we proved a
statement in which in place of $\GG_b$ we have a graph isotopic to
the graph symmetric to $\GG_b$ with respect to a horizontal line.
However, it is easy to show that the symmetry does not change the
morphism. This follows from isotopy invariance of the image with
respect to $\A^1$ and application of the following isotopy:
$\fig{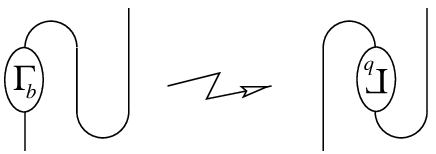}$.
\end{proof}

Here is a counter-part of Lemma \ref{p6.2.A} concerning $\A^2$. It has
a similar proof.

\begin{prop}[Lemma]\label{p6.2.B} Let  $P=(B,W,W^+\to B^{\times})$ be 
a 2-palette (see Section \ref{sec3.4}). Assume that $B$ is a field. 
Let $\GG$ be a $\cal G^2_P$-colored framed closed generic graph 
embedded in $\R^3$ such that its projection to $\R^2$ is generic. Let 
$D\subset \R^2$ be a disk containing all the projection of $\GG$ 
except two arcs, $a$ and $b$, colored with $A$ and $B$, respectively. 
Let $\GG_a$ and $\GG_b$ be colored framed graphs which are obtained 
from $\GG$ by cutting $a$ and $b$, respectively, and moving the end 
points upwards and downwards to a position such that one of them 
becomes the uppermost and the other one lowermost. While moving, the 
framing at the end points is kept black-board. Then  $\A^2(\GG_a)$ is 
multiplication by $\im(\im^A+\im^{-A})\GD$ and $\A^2(\GG_b)$ is 
multiplication by $\im(\im^B+\im^{-B})\GD$ with the same $\GD\in B$. 
\qed\end{prop}

\subsection{Definition of Alexander Invariants}\label{sec6.3}
Let $P=(B,M,W,M\times W\to M)$  be a 1-palette and $B$ be a
field.  Let  $\GG$ be a $\cal G^1_P$-colored framed closed generic
graph in $\R^3$. Cut $\GG$ at  a point on any outermost arc $a$ of its
diagram, where the framing is  black-board and move the new end points
up and down keeping the framing  black-board during the movement. The
result represents a morphism $\GG_a$  of $\cal G^1_{P}$. Let the color
of $a$ is $(t,T)$. Then  $\dfrac1{t^2-t^{-2}}\A^1(\GG_a)$ is
multiplication by some element of  $B$. Denote this element by
$\UD^1(\GG)$ and call it the {\em  $gl(1|1)$-Alexander invariant} of
$\GG$. By  Lemma \ref{p6.2.A}, $\UD^1(\GG)$  does not depend of the
choice of $a$. Of course, it depends on the  coloring. 

Let  $P=(B,W,W^+\to B^{\times})$,  be a 2-palette
and let $B$ be a
field. Let $\GG$ be a $\cal G^2_P$-colored framed closed generic graph
in $\R^3$. Cut $\GG$ at  a point on any outermost arc $a$ of its
diagram,  where the framing is  black-board and move the new end points
up and down keeping the framing  black-board during the movement. The
result represents a morphism $\GG_a$  of $\cal G^2_{P}$. Let the color
of $a$ is $A$. Then  $\dfrac1{\im^{1+A}+\im^{1-A}}\A^2(\GG_a)$ is
multiplication by some element of  $B$. Denote this element by
$\UD^2(\GG)$ and call it the {\em  $sl(2)$-Alexander invariant} of
$\GG$. By  Lemma \ref{p6.2.B}, $\UD^2(\GG)$  does not depend of the
choice of $a$. Of course, it depends on the  coloring. 

When the choice is obvious from the context, we will call
$gl(1|1)$- and $sl(2)$-Alexander invariants just the {\em Alexander
invariant.\/} The Alexander invariant $\UD^c(\GG)$ is invariant under
(ambient) isotopy of $\GG$. This follows  from isotopy invariance of
$\A^c(\GG_a)$ and an obvious construction  turning an isotopy of $\GG$
into an isotopy of $\GG_a$. 

\subsection{Vertex State Sum Representation}\label{sec6.4} Geometric 
operation of cutting and moving the cut points is not necessary if one 
needs just to calculate $\UD^c(\GG)$. Instead one can use the following 
calculation.

First, consider the $gl(1|1)$-case. Let $\GG\subset\R^3$ be a closed
$\cal G^1$-colored graph. Choose its generic projection  to $\R^2$. On
the boundary  of the exterior domain of the diagram choose an arc.
Choose for this arc  either bosonic or fermionic color and consider all
the extensions of this  choice to the whole set of 1-strata of $\GG$.
Choose a generic projection  of the diagram to a line. For each
distribution of bosons and fermions on  the 1-strata of $\GG$ put
according to Tables \ref{tab:BW3} and  \ref{tab:BW4} the Boltzmann
weights at the critical points of the  projection and vertices of the
diagram and make the product of them. Sum  up all the products for all
the distributions. If on the selected arc the  fermion have been
chosen, multiply the sum by $-1$. If the selected arc is  colored with
multiplicity $t$, divide the result by $t^2-t^{-2}$. If the 
orientation of this arc defines the clockwise direction of moving
around  the diagram of $\GG$, multiply by $t^2$. Otherwise, multiply 
by $t^{-2}$. The final result is $\UD^1(\GG)$. 

To prove that this is really $\UD^1(\GG)$, one can cut $\GG$ at a
generic  point on the chosen arc, move the cut points, and calculate
$\UD^1(\GG)$  according to the definition applied to this cut. Moving
the cut points  created two new critical points. To calculate the image
of the graph under  $\A^1$, it is sufficient to calculate the image of
one of the basis  vectors, either boson or fermion, and divide the
result by $t^2-t^{-2}$. The  calculation described above does this in
terms of Boltzmann weights.

Now consider the $sl(2)$-case. Let $\GG\subset\R^3$ be a closed
$\cal G^2$-colored graph. Choose its generic projection  to $\R^2$. On
the boundary  of the exterior domain of the diagram choose an arc. Put
on this arc a light arrowhead defining the clockwise direction
of moving around the diagram.   Consider all the extensions of this 
choice of light arrowhead to secondary colorings of the diagram. Choose
a generic projection  of the diagram to a line. For each distribution
of orientations on  the 1-strata of $\GG$ put according to
Tables \ref{tab:BW1-sl2} and  \ref{tab:BW2-sl2} the Boltzmann weights at
the critical points of the  projection and vertices of the diagram and
make the product of them. Sum  up all the products for all the
distributions.  If the selected arc is  colored with weight $A$, divide
the sum by $1+\im^{-2A}$.  The result is $\UD^1(\GG)$. 

\subsection{Vanishing on Splittable}\label{sec6.5} As follows from Theorem 
\ref{p5.1.A}, the Alexander invariant vanishes on any closed
$\cal G^c_P$-colored 
framed graph $\GG$, which can be presented as a union of two closed 
subgraphs separated from each other by a sphere embedded in $\R^3$.

\section{Special Properties of $gl(1|1)$-Alexander Invariant}\label{sec7}

\subsection{First Examples: Unknot and Theta Graph}\label{sec7.1}
The $gl(1|1)$-Alexander polynomial of the unknot $U$ with the 0-framing
colored with color $(t,T)$ is 
$\dfrac1{t^2-t^{-2}}.$  More generally, for the unknot $U_k$ with 
framing $k/2$, $k\in \Z$ colored with $(t,T)$,
$$\UD^1(U_k)=\frac{t^{-kT}}{t^2-t^{-2}}.$$

For any planar theta graph $\Theta$ with planar framing and arbitrary 
coloring with counter-clockwise orientations at strong vertices  
$\UD^1(\Theta)=1$. This follows from \ref{p5.2.A}.

\subsection{$gl(1|1)$-Alexander Invariant of Tetrahedron's
1-Skeleton}\label{sec7.2} 
Consider now a planar 1-skeleton of a tetrahedron with the planar framing. 
Denote it by $\GT$. This graph plays an important role in the face model, 
see Section \ref{sec10}. Contrary to the case of the theta graph, the 
$gl(1|1)$-Alexander invariant of $\GT$ depends of orientation, and
we have to begin with classification of the orientations. 

\begin{figure}[ht]
\centerline{\includegraphics{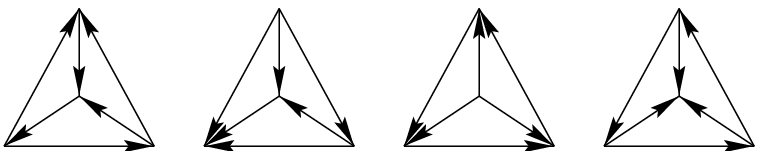}}
\caption{}
\label{ortetr}
\end{figure}

\begin{prop}[Orientations of $\GT$]\label{p7.2.A} Up to
homeomorphism, there are 4 orientations on $\GT$.
They are classified by the numbers of repulsing and attracting vertices.
All the orientations without strong vertices are
homeomorphic. All the orientations having both a repulsing and an
attracting vertices also comprise one homeomorphism type. There is a
single type with one repulsing and no attracting ones, and a single type
with one attracting and no repulsing vertices. See Figure
\ref{ortetr}.\qed \end{prop}

\begin{prop}[Orientations Unextendable to Colorings]\label{p7.2.B} If
an orientation of $\GT$ has only one strong vertex then it cannot be
extended to an admissible $\cal G^1_P$-coloring of $\GT$. \end{prop}

\begin{proof}Indeed, consider the orientation with  repulsing and no 
attracting vertices. Assume that an admissible coloring exists. Denote by 
$T_1$, $T_2$, $T_3$ the weights of edges adjacent to the repulsing vertex. 
Let $U_1$, $U_2$, $U_3$ be the weights of the edges opposite to the edges 
with weights $I$, $J$, $K$, respectively. Assume that the edge with weight 
$L$ is directed towards the end point of the edge with weight $K$. See 
Figure \ref{reputetr}. 

\begin{figure}[htb]
\centerline{\includegraphics{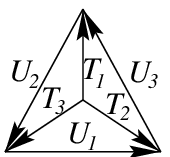}}
\caption{}
\label{reputetr}
\end{figure}
Then, by Admissibility Condition \ref{T-admiss}, $T_1+T_2+T_3=-1$,
$U_1=T_2+U_3-1$, $U_2=T_3+U_1-1$, $U_3=T_1+U_2-1$. The sum of all these
four equations is $0=-4$.

Similarly one can prove that the orientation with an attracting and no
repulsing vertices cannot be extended to an admissible coloring.
\end{proof}

The other two orientations can be extended to colorings. 

\begin{prop}[Lemma on Characteristic Edge]\label{p7.2.B1} For an 
orientation of $\GT$ extendable to an admissible $\cal G^1_P$-coloring, 
there exists a unique edge such that reversing its 
orientation turns the orientation of $\GT$ to an orientation with a 
repulsing vertex and without an attracting one. \qed
\end{prop}

The edge specified by Lemma \ref{p7.2.B1} is said to be {\em 
characteristic}. In Figure \ref{ortetr} the bottom edges in the two 
leftmost graphs are characteristic. It can be defined also as the only
edge connecting two odd vertices, see Section \ref{sec2.2}.
                                                           
\begin{prop}[$gl(1|1)$-Alexander Invariant of $\GT$]\label{p7.2.C} Let
$\GT$ be an  colored planar 1-skeleton of a tetrahedron with planar
framing. Then  $\UD^1(\GT)$ is equal to $t^2-t^{-2}$, where $t$ is the
multiplicity of the  color of the characteristic edge of $\GT$.
\end{prop} 

I am not aware about any conceptual proof of \ref{p7.2.C}. It is proved by 
a straightforward calculation, and I left it to the reader. To save 
efforts, I would recommend to choose the edge for the cut in such a way 
that there is a color extending from this edge uniquely to a coloring of 
$\GT$. In the case with no strong vertices, one can take the edge 
preceding to the characteristic one in the 4-cycle and color it with 
boson. In the case of two strong vertices, color with boson the edge 
connecting the strong vertices. In both cases all the four edges meeting 
the chosen one must be colored with fermion and the disjoint edge with 
boson. \qed 

The calculations of this section demonstrate that the Alexander invariant
is much simpler than the Jones polynomial. Indeed, the Jones polynomial
for the same elementary graphs are much more complicated. See e.~g.
\cite{BHMV}.

\subsection{Functorial Change of Colors}\label{sec7.3}  Let
$P_1=(B_1,M_1,W_1,M_1\times W_1\to M_1)$ and
$P_2=(B_2,M_2,W_2,M_2\times W_2\to M_2)$ 
be 1-palettes,
let $S$ be a subring of $B_1$ containing $M_1$ and $W_1$  and
$\Gb:S\to B_2$  be a ring homomorphism such that
$\Gb(M_1)\subset M_2$, $\Gb(W_1)\subset W_2$ and
$\Gb(m)^{\Gb(w)}=\Gb(m^w)$ for any $m\in M_1$ and $w\in W_1$.

Let $\GG_1$ be a morphism of $\cal G^1_{P_1}$ such that the  multiplicity of 
each 1-stratum of $\GG_1$ does not belong to $\Gb^{-1}\{t\in M_2\mid
t^4=1\}$. Then replacing on each 1-stratum the 
multiplicity and weight with their images under $\Gb$
(and preserving the other components of the coloring) turns
$\GG_1$ to a morphism $\GG_2$ 
of $\cal G^1_{P_2}$. Each entry of the matrix representing $\A^1(\GG_2)$ in 
the standard bases is equal to the image of the corresponding entry of the 
matrix representing $\A^1(\GG_1)$ under $\Gb$. This follows immediately
from the definition of $\A^1$.

Assume now that $B_1$ in $B_2$ are fields. Then there is a similar
relation  for the Alexander invariant: if $\GG_1$ be a $\cal
G^1_{P_1}$-colored closed framed  generic graph in $\R^3$, then
replacing the multiplicities and weights of 1-strata of  $\GG_1$ with
their images under $\Gb$  turns $\GG_1$ to a $\cal G^1_{P_2}$-colored
graph $\GG_2$  and $\UD^1(\GG_2)$ is equal to the image of 
$\UD^1(\GG_1)$ under  $\Gb$.

\subsection{Universal Colors}\label{sec7.4} As we have seen in Section
\ref{sec7.2}, restrictions on the orientations and weights imposed by 
Admissibility Conditions \ref{p2.8.A} are quite subbtle. The
multiplicities are more straightforward: given orientations and
weights, the multiplicities are essentially just a 1-cycle with
coefficients in the multiplicative group $M$. The only non-trivial
condition on this cycle is that it takes on each 1-stratum a value
which is not a fourth root of unity, and, in particular, is not 1.
This follows from Admissibility Condition \eqref{t-admiss}.
For  example, a graph consisting of two disjoint circles and a segment,
which  joins them, does not admit a coloring because 
any 1-cycle vanishes on the segment. 

Among all the 1-cycles on a graph $\GG$ with all the coefficient groups
there is the {\em universal 1-cycle} such that any other 1-cycle can be
obtained from this one by a change of coefficients (i.~e., as in the
previous  section).  The coefficient group for the universal cycle is
$H^1(\GG;\Z)$. The  universal cycle is constructed as follows. An
(oriented) 1-stratum $e$ of  $\GG$ can be considered as a 1-cocycle  of
$\GG$. (On  a 1-chain $c$ this 1-cocycle takes value equal to the
coefficient with which $e$ enters in  $c$.) It is easy to check that
assigning to a 1-stratum the cohomology  class defined in this way is a
cycle. Denote it by $u$. 

Let $M$ be an abelian group. Fix some $x\in H_1(\GG;M)$. Evaluation 
cohomology classes on $x$ defines a homomorphism $y\mapsto y\cap x
:H^1(\GG;\Z)\to M$, and $x$ is the image of $u$ under the homomorphism 
$H_1(\GG;H^1(X;Z))\to H^1(\GG;M)$ induced by $y\mapsto y\cap x
:H^1(\GG;\Z)\to M$.

Let $P=(B,M,W,M\times W\to M)$ be a 1-palette with $W=\Z$.
Let $\GG$ be a framed closed
generic graph in $\R^3$ equipped with orientations. Let its diagram is
equipped with weights satisfying Admissibility Condition
\eqref{T-admiss}. 
If these orientations and weights extends to a $\cal G^1_P$-coloring
on the diagram of $\GG$, it extends  also to  the quartiple
which is made of
$$B=Q(H^1(\GG;\Z)),\quad M=H^1(\GG;\Z),\quad W=\Z. $$
Moreover, $\UD^1(\GG)$ with any 
$\cal G^1_P$-coloring such that $W=\Z$ can be expressed, in the way
shown in Section \ref{sec6.5}, in 
terms of the Alexander invariant of $\GG$ with a coloring whose 
multiplicity component is the universal 1-cycle. 

Notice that the Alexander invariant with the universal multiplicities is a 
ratio of formal linear combinations of cohomology 1-classes of the
graph.  To avoid confusion of addition in $H^1(\GG)$ and the formal
addition, it  makes sense to write the addition in $H^1(\GG)$ as
multiplication.

\subsection{Dependence of the Alexander Invariant on Framings and 
Weights: Case of Links}\label{sec7.5} Relation \eqref{A-r-tw-skein} 
describes dependence of $\A^1(\GG)$ and thereupon $\UD^1(\GG)$ of 
framing. As follows from this relation, if $T=0$, the           
$\UD^1(\GG)$ does not change when the framing changes. In the case of 
link one can choose a coloring such that all the weight components of 
the colors are zero, eliminating dependence of $\UD^1(\GG)$ on 
framing. 

Introduction non-zero weights and framings does not enrich the 
Alexander invariant: still it can be calculated if the linking matrix and  
the Alexander invariant of the same link with the zero weights are 
known. By the linking matrix I mean the matrix comprised of the 
self-linking numbers (i.~e., the framing numbers) and pairwise linking 
numbers of the components. 

These numbers are easy to calculate out of a link diagram. The linking
number of two components is one half of the sum of the signs (local writhe
numbers) of all the crossing points, where the projections of the
components intersect. The self-linking number is the sum of the signs of
all the crossing points, where the projection of the component intersects
itself, plus one half of the sum of the signs of the half-twists on the
projection of the component.

Since several components can be colored with the same color, we need a 
reduced linking matrix. For a framed colored link $L$, denote by 
$\lk_{L}((u,U)(v,V))$ the sum of  linking numbers of all the components of 
$L$ colored with $(u,U)$ and all the components colored with $(v,V)$. 
Denote by $\lk_L((u,U)(u,U))$ the sum of all the self-linking numbers of 
the components of $L$ colored with $(u,U)$ plus the doubled sum all their 
pairwise linking numbers. To calculate these numbers for a link presented 
by its diagram, one can use the same rule as above: if $(u,U)\ne(v,V)$ 
then $\lk_L((u,U)(v,V))$ is the half of the sum of the signs of all the 
crossing points where the colors $(u,U)$, $(v,V)$ meet; 
$\lk_L((u,U)(u,U))$ is the sum of signs of the crossing points where 
both strings colored with $(u,U)$ plus one half of the sum of the signs of 
the half-twist on strings colored with $(u,U)$. 

\begin{prop}\label{p7.5.A}
Let $L$ be a $\cal G^1_P$-colored framed link in $\R^3$ and $L_0$ the
$\cal G^1_P$-colored 
framed link obtained from $L$ by replacing the weight components of the 
color with zero. Then
\begin{equation}\label{reducing-weights}
\UD^1(L)=\prod_{(u,U)}u^{-2\sum_{(v,V)}V\lk_L((u,U),(v,V))}\UD^1(L_0)
\end{equation}
\end{prop}
\begin{proof}
This follows from a straightforward calculation of the contribution made
by the weights of colors in the Boltzmann weights. The weights of colors
contribute only to the Boltzmann weights at half-twist signs and crossing
points. At a positive half-twist of a component colored with $(u,U)$ it is
$u^{-U}$, at a positive crossing point with colors $(u,U)$ and $(v,V)$ it
is $u^{-V}v^{-U}$. The product of all these factors 
is $\prod_{(u,U)}u^{-2\sum_{(v,V)}V\lk_L((u,U),(v,V))}$.\end{proof}

\subsection{Dependence of the Alexander Invariant on
Weights: General Case}\label{sec7.6}
In presence of trivalent vertices, it is impossible to make all the weights
of colors equal to zero, since if at a trivalent vertex two of the three
adjacent edges are colored with the zero weight, then, due to the
admissibility condition, the third one is colored with weight $\pm1$.
However, it is possible by a formula similar to \eqref{reducing-weights}
to relate the Alexander invariants of a framed graph colored in two
ways with colors different only in their weight components.

Let $\GG$ be a framed graph equipped with two $\cal G^1_P$-colorings 
with colors different only in their weight components.Denote  by $T^k_e$ 
the value of the weight of the $k$th (with $k=1,2$) coloring taken on a 
oriented (with the  orientation shared by the colorings) 1-stratum $e$. 
Then $\sum(T^1_e-T^2_e)e$, where $e$ runs over the set of 1-strata of
$\GG$, is a 1-cycle with integer coefficients. This follows from the 
admissibility conditions \eqref{T-admiss}. Let us call this cycle by 
the {\it weight cycle\/} and denote it by $\GG_W$. Another cycle is 
made in the same way, but with $T^1_e-T^2_e$ replaced by the 
multiplicities (shared by the colors) of 1-strata. Let us call this 
cycle by the {\it multiplicity cycle\/} and denote by $\GG_M$. 

Consider a field of lines on $\GG$ normal to the framing. Choose on each
of these lines a short segment centered at the corresponding point of
$\GG$ such that the end points of the segments comprise a generic graph
$\tilde\GG$ which is a two-fold covering space of $\GG$. Equip
$\tilde\GG$ with the orientation such that the covering projection maps
it to the orientation of $\GG$ incorporated into the colorings. Equip each
1-stratum of $\tilde\GG$ with the multiplicity of its image
under the projection. The oriented 1-strata of $\tilde\GG$ equipped with
these halves of multiplicities comprise a cycle with coefficients in
$M$. Denote this cycle by $\tilde\GG_M$ 
It can be thought of as
the multiplicity cycle $\GG_M$ pushed away from $\GG$ in both directions
normal to
the framing. The cycle $\tilde\GG_M$ is disjoint from $\GG$, and hence 
from $\GG_W$.
Therefore the linking number $\lk(\tilde\GG_M,\GG_W)$ of this cycle and
the weight cycle associated with the pairing $M\times
W\to M:(m,w)\mapsto m^w$ is well-defined.
In the case when $\GG$ is a link and $T^2_e=0$, 
$$\lk(\tilde\GG_M,\GG_W)=\prod_{(u,U)}u^{2\sum_{(v,V)}V\lk_L((u,U),(v,V
))},$$
cf. \eqref{reducing-weights}.

Denote by $\GG_1$ and $\GG_2$ the graph $\GG$ equipped with the colorings
under consideration. Then

\begin{equation}\label{comparing-weights}
\UD^1(\GG_2)=\lk(\tilde\GG_M,\GG_W)^{-1}\UD^1(\GG_1)
\end{equation}

Formula \eqref{comparing-weights} generalizes \eqref{reducing-weights} and
is proved similarly.

\subsection{Relation of the $gl(1|1)$-Alexander Invariant to the Conway 
Function}\label{sec7.7} Recall that for an oriented link $L$ with (linearly) 
ordered connected components $L_1$, \dots, $L_k$ the Conway function 
$\nabla(L)(t_1,\dots,t_k)$ is a rational function in variables $t_1$, 
\dots, $t_k$ defined by the following axioms (see Turaev \cite{Tu3}): 

\begin{prop}\label{p7.7.A}
$\nabla(L)$ does not change under (ambient) isotopy of $L$.
\end{prop}

\begin{prop}\label{p7.7.B}
If $L$ is the unknot, then $\nabla(L)(t)=(t-t^{-1})^{-1}$.
\end{prop}

\begin{prop}\label{p7.7.C}
If  $L$ consists of more than one connected components, then $\nabla(L)$
is a Laurent polynomial, i.~e.
$\nabla(L)(t_1,\dots,t_k)\in\Z[t_1,t_1^{-1},\dots,t_k,t_k^{-1}]$.
\end{prop}

Denote by $\tilde\nabla(L)$ the function $\nabla(L)(t,t,\dots,t)$ of one
variable.

\begin{prop}\label{p7.7.D}
$\tilde\nabla(L)(t)$ does not depend on the ordering of the connected
components of $L$. \end{prop}

\begin{prop}[Conway Skein Relation]\label{p7.7.E} If
$L_{\fig{poscros.eps}}$, $L_{\fig{negcros.eps}}$ and
$L_{\fig{smtcros.eps}}$ are links coinciding outside a ball and in the
ball looking as the subindices in their notations, then
\begin{equation}\label{conway-skein}
\tilde\nabla(L_{\fig{poscros.eps}})=\tilde\nabla(L_{\fig{negcros.eps}}) +
(t-t^{-1})\tilde\nabla(L_{\fig{smtcros.eps}}).
\end{equation}
\end{prop}

\begin{prop}[Doubling Axiom]\label{p7.7.F} If $L'$ is a link obtained
from a  link $L=L_1\cup\dots\cup L_k$ by replacing the  $L_i$
by its $(2,1)$-cable, then
\begin{equation}\label{doubling}
\nabla(L')(t_1,\dots,t_r)=(T+T^{-1})\nabla(L)(t_1,\dots,t_{i-1},t_i^2,
t_{i+1},\dots,t_k) \end{equation}
where $T=t_i\prod_{j\ne i}t_j^{\lk(L_i,L_j)}$.
\end{prop}

\begin{prop}[Alexander Invariant Versus Conway Function]\label{p7.7.G} 
Let $M_k$ be a multi\-pli\-ca\-tive free abelian group generated by 
$t_1$, \dots, $t_k$, let $B_k$ be the field $Q(M_k)$ of rational 
functions of $t_1$,\dots, $t_k$ with integer coefficients. Denote by 
$P_k$ the quartiple $(B_k,M_k,\Z,M_k\times\Z\to\Z:(m,w)\mapsto m^w)$. 
Let $L$ be a $\cal G^1_{P_k}$-colored framed link with components 
$L_1$, \dots, $L_k$ colored with multiplicities $t_1$, \dots, $t_r$ 
and zero weights (i.~e. $T_1=T_2=\dots=T_k=0$). 
Then  $\UD^1(L)$ and the Conway function of the same $L$ (oriented 
with the orientation taken from the coloring) are related as follows: 
\begin{equation}\label{inv-polynom}
\UD^1(L)=\nabla(L)(t_1^2,\dots,t_k^2).
\end{equation}
\end{prop}

\begin{proof}It suffices to prove that the Alexander invariant of 
$\cal G^1_P$-colored framed link with zero weights satisfies the 
following counterparts for the axioms \ref{p7.7.A} - \ref{p7.7.F}. 

\numthm{7.7.A.bis}
\begin{mynumthm}\label{p7.7.A'}
$\UD^1(L)$ does not change under (ambient) isotopy of $L$.
\end{mynumthm}

\numthm{7.7.B.bis}
\begin{mynumthm}\label{p7.7.B'}
If $L$ is the unknot $\cal G^1_P$-colored with multiplicity $t$, then 
$\UD^1(L)(t)=(t^2-t^{-2})^{-1}$. \end{mynumthm}

\numthm{7.7.C.bis}
\begin{mynumthm}\label{p7.7.C'}
If  $L$ consists of more than one connected components 
$\cal G^1_{P_k}$-colored as 
in \ref{p7.7.G}, then $\UD^1(L)$ is a Laurent polynomial in 
multiplicities, i.~e.
$$\UD^1(L)\in\Z[t_1,t_1^{-1},\dots,t_k,t_k^{-1}].$$
\end{mynumthm}

(This is a weakened counterpart of \ref{p7.7.C}. The true counterpart would 
be $\UD^1(L)\in\Z[t_1^2,t_1^{-2},\dots,t_k^2,t_k^{-2}]$. However 
\ref{p7.7.C'} is also sufficient for proving uniqueness of the invariant 
satisfying the axioms, because in the proof imitating Turaev's one 
\cite{Tu3} it is used anyway in combination with the Doubling Axiom, which 
allows one to get into the domain of Laurent polynomial by additional 
doubling.) 

Let $L$ be a $\cal G^1_{P_k}$-colored framed link as in \ref{p7.7.G}.
Denote by $L^t$  the same framed link, but with all the components
colored with $(t,0)$ and  orientations taken from the coloring of $L$. 

\numthm{7.7.D.bis}
\begin{mynumthm}\label{p7.7.D'}$\UD^1(L^t)(t)=\UD(L)^1(t,\dots,t)$.
\end{mynumthm}

\numthm{7.7.E.bis.}
\begin{mynumthm}[Conway Skein Relation]\label{p7.7.E'} If
$L_{\fig{poscros.eps}}$, $L_{\fig{negcros.eps}}$ and
$L_{\fig{smtcros.eps}}$ are links coinciding outside a ball and in the
ball looking as the subindices in their notations, then
\begin{equation}\label{conway-skein'}
\UD^1(L_{\fig{poscros.eps}}^t)=\UD^1(L^t_{\fig{negcros.eps}}) +
(t^2-t^{-2})\UD^1(L^t_{\fig{smtcros.eps}}).
\end{equation}
\end{mynumthm}

\numthm{7.7.F.bis.}
\begin{mynumthm}[Doubling Axiom]\label{p7.7.F'} If $L'$ is a link obtained
from a  link $L=L_1\cup\dots\cup L_k$ by replacing the  $L_i$
by its $(2,1)$-cable, then
\begin{equation}\label{doubling'}
\UD^1(L')(t_1,\dots,t_r)=(T+T^{-1})\UD^1(L)(t_1,\dots,t_{i-1},t_i^2,
t_{i+1},\dots,t_k) \end{equation}
where $T=t_i^2\prod_{j\ne i}t_j^{2\lk(L_i,L_j)}$.
\end{mynumthm}

As is mentioned above \ref{p7.7.A'}, \ref{p7.7.B'}, \ref{p7.7.D'} holds true.

To prove \ref{p7.7.C'}, observe, first, that for any
endomorphism of a one-point object of 
$\cal G^1_{P}$,
the image under $\A^1$ is multiplication by a Laurent polynomial in
$t_1$, \dots, $t_k$. Indeed, the Boltzmann weights for  all the
generators of $\cal G^1_{P_k}$ are Laurent polynomials.
To calculate $\UD^1(L)$, one can take $\A^1(L_a)$ with $a$ from 
different
components. Let $a$ be an arc on $L_1$ and $b$ be an arc on $L_2$.
Then both $(t_1^2-t_1^{-2})^{-1}\A^1(L_a)$ and 
$(t_2^2-t_2^{-2})^{-1}\A^1(L_b)$ are multiplication
by $\UD^1(L)$. The first of these morphisms is a multiplication by an 
element of $\Z[M_k]$ divided by $t_1^2-t_1^{-2}$, the second a 
multiplication by an 
element of $\Z[M_k]$ divided by $t_2^2-t_2^{-2}$. Since $t_1^2-t_1^{-2}$ 
and $t_2^2-t_2^{-2}$ are relatively prime in $\Z[M_k]$, $\UD^1(L)$ is an 
element of $\Z[M_k]$.

The Conway Skein Relation \ref{p7.7.E'} follows easily
by comparing the Boltzmann weights at positive and negative crossing
points, see the first column of Table \ref{tab:BW3}.

Consider now the Doubling Axiom \ref{p7.7.F'}. Since the framing does not 
matter when the weights of colors vanish, we may assume that $L_i$ has 
selflinking number $\frac12$ and the $(2,1)$-cable has selflinking number 
2. Under this choice, the $(2,1)$-cable $L_i'$ of $L$ can be obtained as 
the boundary of the framing surface $F$ (which is a M\"obius band) of 
$L_i$ and its framing as overlapping with $F$ along the whole $L_i'$. 
Squeezing two parallel arcs of $L_i'$ along $F$ apply Junction Relation 
\ref{p5.2.B}
$$\UD^1(L')=(t_i^2-t_i^{-2})^{-1}\UD^1(L_+)
-(t_i^2-t_i^{-2})^{-1}\UD^1(L_-),$$
where $L_+$ is a graph obtained from $L'$ by squeezing two parallel arcs on 
$L_i'$ with the orientation induced by the orientation of $L_i'$, and 
$L_-$ the same graph, but with opposite orientation on the new arc. In 
$L_+$ and $L_-$ the new arc has framing along $F$ and weight $-1$. Its 
multiplicity is $t_i^2$ in $L_+$, and $t_i^{-2}$ in $L_-$. 

Expand now the new arcs along $L_i$ and apply Removing Loop Relation
\ref{p5.2.A}
$$\UD^1(L_+)=(t_i^2-t_i^{-2})\UD^1(L^1),\qquad 
\UD^1(L_-)=-(t_i^2-t_i^{-2})\UD^1(L^2)$$
where $L^1$ is $L$ with $L_i$ equipped with multiplicity $t_i^2$ and 
weight $-1$, and $L^2$ is $L$ with $L_i$ equipped with the opposite
orientation, multiplicity $t_i^{-2}$ and weight $-1$.
Hence $$\UD^1(L')=\UD^1(L^1)+\UD^1(L^2).$$

Further, by \eqref{reducing-weights},
$$\UD^1(L^1)=(t_i^2)^{2\lk_L(L_i,L_i)}\prod_{j\ne 
i}t_j^{2\lk(L_j,L_i)} \UD^1(L^1_0)=t_i^2\prod_{j\ne 
i}t_j^{2\lk(L_j,L_i)} \UD^1(L^1_0),$$
where $L^1_0$ differs from  $L$ only by the multiplicity of $L_i$, which
is $t_i^2$. Similarly, by \eqref{reducing-weights},
$$\UD^1(L^2)=(t_i^{-2})^{2\lk_L(L_i,L_i)}\prod_{j\ne 
i}^kt_j^{-2\lk(L_j,L_i)} 
\UD^1(L^2_0)=t_i^{-2}\prod_{j\ne i}^kt_j^{-2\lk(L_j,L_i)} 
\UD^1(L^2_0),$$
where $L^2_0$ differs from $L$ only by the multiplicity of $L_i$, which
$t_i^{-2}$ and orientation of $L_i$.

It is easy to see that squaring the multiplicity $t_i$ of a link component
corresponds to replacing $t_i$ with $t_i^{2}$ in all the Boltzmann
weights, while simultaneous reversing orientation and reversing the 
multiplicity of a link component does not change the Boltzmann
weights. Therefore $\UD^1(L^1_0)$ and $\UD^1(L_0^2)$ are obtained from
$\UD^1(L)$ by replacing $t_i$ with $t_i^2$. This corresponds to
replacement $t_i$ with $t_i^2$ in $\UD^1(L)$. The coefficient
$t_i^2\prod_{j\ne i}t_j^{2\lk(L_j,L_i)}$ equals $T$.
\end{proof}

\section{Special Properties of $sl(2)$-Alexander Invariant}\label{sec8}

\subsection{First Calculations: Unknot and Theta Graph}\label{sec8.1}
The $sl(2)$-Alexander polynomial  of the unknot $U$ with the 0-framing
colored with color $A$ is 
$\dfrac1{\im^{1+A}+\im^{1-A}}.$ For the unknot $U_k$ with 
framing $k/2$, $k\in \Z$ colored with $A$,
$$\UD^2(U_k)=\frac{\im^{k\frac{A^2-1}4}}{\im^{1+A}+\im^{1-A}}.$$

For any planar theta graph $\Theta$ with planar framing and arbitrary  
coloring, $\UD^2(\Theta)=1$. This follows from \ref{p5.3.A}.

\subsection{$sl(2)$-Alexander Invariant of Tetrahedron's 
1-Skeleton}\label{sec8.2} Recall that a planar 1-skeleton of a 
tetrahedron with the planar framing is denoted by $\GT$. In a $\cal 
G^2_P$ coloring of $\GT$, two vertices of $\GT$ are sources and the 
other two vertices are sinks. Indeed, the augmentation of the 0-cycle 
$\text{stocks}-\text{sources}$ is zero since this 0-cycle is the 
boundary of the weight chain.

\begin{prop}[$sl(2)$-Alexander Invariant of $\GT$]\label{p8.2.A} Let
$\GT$ be an  colored planar 1-skeleton of a tetrahedron with planar
framing. Then  $\UD^2(\GT)$ is equal to $\im^{1+A}+\im^{1-A}$, where 
$A$
is the weight of the edge connecting the two sink vertices of $\GT$.
\end{prop} 

It is proved by a straightforward calculation, and I left it to the
reader. To save efforts, I would recommend to choose the edge for the
cut in such a way that there exists a secondary coloring, i.e. an
orientation, of this edge  extending uniquely to a secondary coloring
of $\GT$. One can choose an edge connecting a source with a sink and
orient it from the source to the sink. 
\qed 

\subsection{Functorial Change of Colors}\label{sec8.3}  Let
$P_1=(B_1,W_1,W_1^+\to B_1^{\times})$ and
$P_2=(B_2,W_2,W_2^+\to B_2^{\times})$ 
be 2-palettes,
let $S$ be a subring of $B_1$ containing  $W_1$  and
$\Gb:S\to B_2$  be a ring homomorphism such that
$\Gb(W_1)\subset W_2$ and
$\Gb(\im^w)=\im^{\Gb(w)}$ for any $w\in W_1$.

Let $\GG_1$ be a morphism of $\cal G^2_{P_1}$ such that the weight of 
each 1-stratum of $\GG_1$ does not belong to $\Gb^{-1}\{x\in W_2\mid
\im^{2x+2}=1\}$. Then replacing on each 1-stratum the 
weight with its image under $\Gb$
(and preserving the orientation component of the coloring) turns
$\GG_1$ to a morphism $\GG_2$ 
of $\cal G^2_{P_2}$. Each entry of the matrix representing $\A^2(\GG_2)$ in 
the standard bases is equal to the image of the corresponding entry of the 
matrix representing $\A^2(\GG_1)$ under $\Gb$. This follows immediately
from the definition of $\A^2$.

Assume now that $B_1$ in $B_2$ are fields. Then there is a similar
relation  for the $sl(2)$-Alexander invariant: if $\GG_1$ be a $\cal
G^2_{P_1}$-colored closed framed  generic graph in $\R^3$, then
replacing the weights of 1-strata of  $\GG_1$ with
their images under $\Gb$  turns $\GG_1$ to a $\cal G^2_{P_2}$-colored
graph $\GG_2$  and $\UD^2(\GG_2)$ is equal to the image of 
$\UD^2(\GG_1)$ under  $\Gb$.

\subsection{Effect of Adding Numbers Divisible by 4 to Weights:
the Case of a Link}\label{sec8.4}

\begin{prop}\label{p8.4.A}Let $L$ be a $\cal G^2_P$-colored framed
link with components $L_1$, \dots, $L_n$ colored with weights
$A_1$,\dots, $A_n$, respectively. Let $p_1$, \dots, $p_n$ be  
integers
divisible by 4 and $L'$ be a $\cal G^2_P$-colored framed link made of
the same components $L_1$, \dots, $L_n$ as $L$ but colored with weights
$p_1+A_1$,\dots, $p_n+A_n$. Then
\begin{equation}\label{sl-red-weights}
\UD^2(L')=\UD^2(L)\im^{\sum_{i=1}^n\sum_{j=1}^n(1+A_i)p_j\lk(L_i,L_j)}
\end{equation}
\end{prop}

\begin{proof}
This follows from a straightforward calculation of the contribution made
by change of the weights in the Boltzmann weights. 
Only the Boltzmann weights at half-twist signs and crossing
points are effected. At a positive half-twist of the $i$th component 
the Boltzmann weight is multiplied by $\im^{\frac{A_ip_i}2}$, at a
positive crossing point of $i$th and $j$th components the Boltzmann
weight is multiplied by $\im^{\frac{(1+A_i)p_j+(1+A_j)p_i}2}$. 
The product of all these factors is
$\prod_{i=1}^n\im^{\sum_{j=1}^n(1+A_i)p_j\lk(L_i,L_j)}$.\end{proof}

\subsection{Effect of Adding Numbers Divisible by 4 to Weights:
General Case}\label{sec8.5}
Let $\GG$ be a framed graph equipped with two $\cal G^2_P$-colorings
such that  on each 1-stratum of $\GG$ their weights differ by 
an integer divisible by 4. Denote  by $A^k_e$ the value of the 
weight of the $k$th (with $k=1,2$) coloring taken on a oriented (with the  
orientation shared by the colorings) 1-stratum $e$. 
Then $\sum(A^1_e-A^2_e)e$, where $e$ runs over the set of 1-strata of
$\GG$, is a 1-cycle with coefficients in $4\Z$. 
This follows from the admissibility conditions \eqref{p3.1.A}. 
Let us denote it by $\GG_W$.
Another cycle is with coefficients in $M/2\Z$. This is the weight cycle
$\sum(1+A^1_e)e=\sum(1+A^2_e)e\mod2$. Denote it by $\GG_M$.

Consider a field of lines on $\GG$ normal to the framing. Choose on each
of these lines a short segment centered at the corresponding point of
$\GG$ such that the end points of the segments comprise a generic graph
$\tilde\GG$ which is a two-fold covering space of $\GG$. Equip
$\tilde\GG$ with the orientation such that the covering projection maps
it to the orientation of $\GG$ incorporated into the colorings. Equip each
1-stratum of $\tilde\GG$ with the multiplicity of its image
under the projection. The oriented 1-strata of $\tilde\GG$ equipped with
these halves of multiplicities comprise a cycle with coefficients in
$M$. Denote this cycle by $\tilde\GG_M$ 
It can be thought of as
the multiplicity cycle $\GG_M$ pushed away from $\GG$ in both directions
normal to
the framing. The cycle $\tilde\GG_M$ is disjoint from $\GG$, and hence 
from $\GG_W$.
Therefore the linking number $\lk(\tilde\GG_M,\GG_W)$  associated
with the multiplication pairing $M/2\Z\times 4\Z\to M/4\Z$ is well-defined.
In the case when $\GG$ is a link $L=L_1\cup\dots\cup L_n$ considered in
the preceding section, 
$$\lk(\tilde\GG_M,\GG_W)=\sum_{i=1}^n\sum_{j+1}^n(1+A_i)p_j\lk(L_i,L_j)
\mod4,$$
cf. \eqref{sl-red-weights}.

Denote by $\GG_1$ and $\GG_2$ the graph $\GG$ equipped with the colorings
under consideration. Then

\begin{equation}\label{sl-comp-weights}
\UD^2(\GG_2)=\im^{\lk(\tilde\GG_M,\GG_W)}\UD^2(\GG_1)
\end{equation}

Formula \eqref{sl-comp-weights} generalizes \eqref{sl-red-weights} and
is proved similarly.

\subsection{Relation of the $sl(2)$-Alexander Invariant to the Conway 
Function}\label{sec8.6}
\begin{prop}\label{p8.6.A}
Let $L$ be an oriented framed link with components
$L_1$, \dots, $L_n$ equipped with $\cal G^2_P$-coloring. Let $L_i$ be
colored with weight $A_i$. Then
\begin{equation}\label{e8.6.1}
\UD^2(L)=\nabla(L)(\im^{1+A_1},\dots,\im^{1+A_n})\im^{\sum_{i,j=1}^n
\frac{A_iA_j-1}2\lk(L_i,L_j)}
\end{equation}
\end{prop}

\begin{proof}
Let us construct (cf. Section \ref{sec4.5}) a 1-palette $R$
It consists of the ring $B$ involved in $P$, the image
$M$ of the homomorphism $W^+\to B^{\times}:A\to\im^A$, the subring
$W\subset B$ and the pairing $M\times W\to M:
(\im^A,B)\mapsto\im^{AB}$. By \ref{p4.5.E} any secondary $\cal
G^2_P$-coloring of the diagram of $L$ is the image of a secondary $\cal
G^1_R$-coloring under $\GF^{\flat}_Q$. Hence by \ref{p4.5.D} 
$\UD^2(L)=\UD^1(L')$, where $L'$ is the same oriented framed link $L$, 
but equipped with the corresponding $\cal G^1_R$-coloring. The $\cal 
G^1_R$-color of $L_i$ includes multiplicity 
$\im^{\frac{1+A_i}2}$ and weight $\frac{1-A_i}2$. Denote by $L^0$ the 
same link $L$ equipped with $\cal G^1_R$-coloring with the same 
orientations and multiplicities but zero weights. By 
\ref{p7.5.A},
$$\UD^1(L')= \UD^1(L^0)\prod_{i=1}^n
\im^{\frac{1+A_i}2\(-2\sum_{j=1}^n\frac{1-A_j}2\lk(L_i,L_j)\)} 
$$
By \ref{p7.7.G}
$$\UD^1(L^0)=\nabla(L)(\im^{1+A_1},\dots,\im^{1+A_n})$$  
Hence,
\begin{multline}
\UD^2(L)= \UD^1(L')=  \\
         \UD^1(L^0)\prod_{i=1}^n
\im^{-(1+A_i)\sum_{j=1}^n\frac{1-A_j}2\lk(L_i,L_j)}= \\
        \nabla(L)(\im^{1+A_1},\dots,\im^{1+A_n})\im^{\sum_{i,j=1}^n
\frac{A_iA_j-1}2\lk(L_i,L_j)}.
\end{multline}
\end{proof}

 \section{Graphical Skein Relations}\label{sec9}

\subsection{Another Look of Skein Relations}\label{sec9.1} The Alexander 
invariant for colored framed generic graphs provides terms in which one 
can rewrite skein relations of Section \ref{sec5}, and write similar, but 
more complicated ones. 

Notice, first, that since
$$\UD^1\(\fig{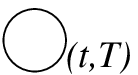}\)=\frac1{t^2-t^{-2}},$$
the skein relations \eqref{Aphi-gr-skein} and  \eqref{p5.2.B4}
can be interpreted in the following way:

\begin{equation}\label{skein2-g}
\A^1\(\ {\fig{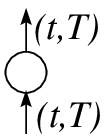}}\)= \frac1{\UD^1\(\fig{circle-t.eps}\)}
\A^1 \(\fig{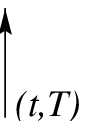}\)
\end{equation}

\begin{equation}\label{skein3-g}
\A^1\(\fig{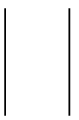}\)=\sum_{(t,T)}\UD^1\(\fig{circle-t.eps}\)
\A^1\(\fig{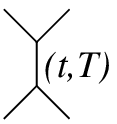}\) \end{equation}

Similarly, since
$$\UD^2\(\fig{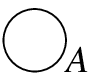}\)=\frac1{\im^{1+A}+\im^{1-A}},$$
the skein relations \eqref{Aphi-gr-skein} and  \eqref{p5.3.B2}
can be interpreted in the following way:

\begin{equation}\label{skein2-g2}
\A^2\(\ {\fig{phi-gruS.eps}}\)= \frac1{\UD^2\(\fig{circle-A.eps}\)} 
\A^2 \(\fig{vert-A.eps}\)
\end{equation}

\begin{equation}\label{skein3-g2}
\A^2\(\fig{ii.eps}\)=
\sum_{A}\UD^2\(\fig{circle-A.eps}\)
\A^2\(\fig{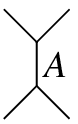}\)
\end{equation}

This more graphical form is used in transition to the face models.

\subsection{New Relations}\label{sec9.2} The relations of this section
are formulated in the graphical form and used in the transition 
to the face models, too. If 
we did not interpret the coefficients as the Alexander invariants of the 
graphs related to the graphs involved, it would be difficult even to 
formulate those relations. Another advantage of the graphical formulation 
is that we can unite the relations for $\A^1$ and $A^2$.

\begin{prop}[Flip Relation]\label{p9.2.A} For $c=1,2$
\begin{equation}\label{skein4}
\A^c\(\fig{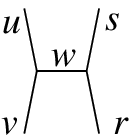}\)\
=\sum_{t}\ \UD^c\(\fig{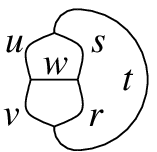}\ \)\ \UD^c\(\fig{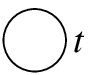}\)\
\A^c\(\fig{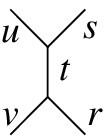}\)
\end{equation}
provided the colorings satisfy the admissibility condition and the sum
runs over nonempty set of colors $t$ (i.~e., there exists a 
color $t$ 
involved together with the given colors in an admissible coloring of the 
graphs in the right hand side of \eqref{skein4}). \end{prop}

\begin{proof} By Theorem \ref{p5.1.B}, the 
left hand side of
\eqref{skein4} can be presented as the following linear combination:
\begin{equation}\label{p9.2.A1}
\A^c\(\fig{h-gr.eps}\)\ =\sum_{t}\ x_t \A^c\(\fig{i-gr.eps}\)
\end{equation}

To find the coefficient at $\A^c\(\fig{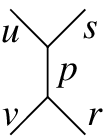}\)$,
adjoin $\fig{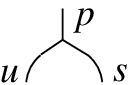}$ to the top and $\fig{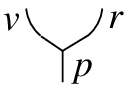}$ to the bottom
of each of the graphs involved in the equality. 
$\A^c\(\fig{h-gr.eps}\)$ turns into 
$$\A^c\(\fig{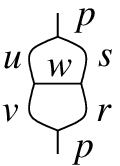}\)=\UD^c\(\fig{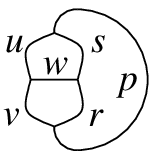}\ \)\ 
\(\UD^c\(\fig{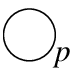}\)\)^{-1}\id.$$ 
By \eqref{skein2-g}, the summand corresponding to 
color $p$ on the right hand side  turns into morphism 
$x_p\(\A^c\(\fig{circle-p.eps}\)\)^{-2}\id$, while 
the other summands annihilate.
Hence we obtain 
$$\UD^c\(\fig{tetr13.eps}\ \)\(\UD^c\(\fig{circle-p.eps}\)\)^{-1}
=x_p\(\UD^c\(\fig{circle-p.eps}\)\)^{-2}$$  and
$$x_p=\UD^c\(\fig{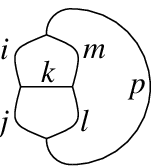}\)\ \UD^c\(\fig{circle-p.eps}\)$$
\end{proof}
 
The next skein relation is similar to \ref{p9.2.A} and admits a similar 
proof. Moreover, it can be deduced from \ref{p9.2.A}, \eqref{A-r-tw-skein} 
and \eqref{2A-r-tw-skein}. 

\begin{prop}[Crossing Relation]\label{p9.2.B} For $c=1,2$
\begin{equation}\label{skein5}
\A^c\(\fig{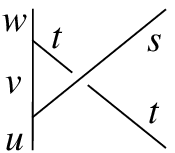}\)\
=\sum_r \UD^c\(\fig{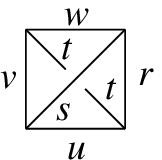}\) \UD^c\(\fig{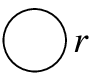}\)
\A^c\(\fig{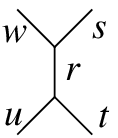}\)
\end{equation}
provided the colorings satisfy the admissibility condition and the sum
runs over nonempty set of colors $r$ (i.~e., there exists a 
color $r$ comprising together with the other colors an 
admissible coloring 
of the graphs in the right hand side of \eqref{skein5}).\qed \end{prop}

\begin{prop}[Triangle Relation]\label{p9.2.C} For $c=1,2$
\begin{equation}\label{skein6}
\A^c\(\fig{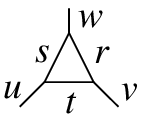}\)=\UD^c\(\fig{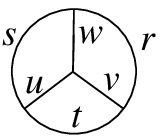}\)\A^c\(\fig{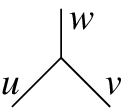}\) 
\end{equation}
provided the colorings satisfy the admissibility condition.
\end{prop}

This also follows from \ref{p9.2.A}, \eqref{A-r-tw-skein} and 
\eqref{2A-r-tw-skein}.\qed

\subsection{Skein Relations for the Alexander
Invariants}\label{sec9.3} All the skein relations for $\A^c$ discussed in
Sections \ref{sec9.1} and \ref{sec9.2}
can be thought of as skein relations for $\UD^c$. One should just
replace $\A^c$  by $\UD^c$ 
and understand $\UD^c(\ \text{nonclosed diagram of graph}\ )$ in the usual 
way, i.~e. as the Alexander invariant of some closed graph whose
diagram contains the shown diagram as a fragment, which is the only
variable part of the diagram: in the other terms of the relation which 
involve nonclosed diagram, the unshown parts of the diagrams are assumed to 
be the same.

\section{Face State Sums}\label{sec10}

\subsection{Colorings of a Diagram}\label{sec10.1} Let 
$\GG\subset\R^3$
be a framed generic closed graph. Fix a diagram for $\GG$.

Let $P=(B,M,W,M\times W\to M)$ be a 1-palette. By a $\cal 
G^1_P$-color of a 2-strata of the diagram\footnote{Recall that a 
2-stratum of a diagram is a connected component of the complement in 
$\R^2$ of the projection of $\GG$.}   we will mean an orientation of 
this 2-stratum, an element $t$ of $M$ with $t^4\ne1$ (called the 
{\em multiplicity} of the 2-stratum), and an element $T$ of $W$ (called 
its {\em weight}). 

A {\it $\cal G^1_P$-coloring of the diagram\/} of $\GG$ is an 
assignment to each of the 1- and 2-strata of the diagram a color such 
that \begin{enumerate} 
\item at each triple vertex the colors of the adjacent edges 
satisfy Admissibility Condition \ref{p2.8.A}, 
\item at each crossing point the colors of the opposite 
edges (which belong to the image of the same 1-stratum of $\GG$) 
coincide, 
\item for each 1-stratum its color and the colors of the two adjacent 
2-strata satisfy the following condition \ref{p10.1.A} (cf. \ref{p2.8.A}): 
\end{enumerate} 

\begin{prop}[Admissibility of 1-Palette Colors at a
1-Stratum]\label{p10.1.A} 
Let $t$ be the multiplicity and $T$ weight components of a 1-stratum's
color, $(t_1,T_1)$ and  $(t_2,T_2)$ the corresponding components of the 
colors of the adjacent 2-strata, let $\Ge_i=1$, if the orientation of the
$i^{\text{th}}$ adjacent 2-stratum induces the orientation of the
1-stratum and $\Ge_i=-1$ otherwise. Then
\begin{equation}\label{p10.1.A1}
 t t_1^{\Ge_1} t_2^{\Ge_2}=1,
\end{equation}
\begin{equation}\label{p10.1.A2}
 T+\Ge_1 T_1+\Ge_2 T_2=-\Ge_1\Ge_2.
\end{equation}
\end{prop}

A 1-stratum of a diagram is said to be {\em strong\/} with respect 
to a $\cal G^1_P$-coloring of the diagram if the orientations of the 
adjacent 2-strata both induce the orientation of the 1-stratum (i.e.,  
$\Ge_1=\Ge_2=1$ in notations of \ref{p10.1.A}).	One can see that 
the union of strong 1-strata is a 1-submanifold of the plane. The 
boundary of this submanifold is the set of images of all strong 
vertices of the graph. The orientations of strong 1-strata comprise a  
natural orientation of their union.

Now let $P=(B,W,W^+\to B^{\times})$ be a 2-palette. By a $\cal 
G^2_P$-color of a 2-strata of the diagram we will mean a 
pair made of an orientation 
of this 2-stratum, and an element $A$ of $B\sminus\{x\in W\mid 
\im^{2x+2}=1\}$ defined up to simultaneous reversing of the 
orientation and multiplication of $A$ by $-1$. The component of the 
color belonging to $B$ is called the {\em weight} of the 2-stratum. 

A {\it $\cal G^2_P$-coloring of the diagram\/} of $\GG$ is an 
assignment to each of the 1- and 2-strata of the diagram a color such 
that \begin{enumerate} 
\item at each triple vertex the colors of the adjacent edges 
satisfy Admissibility Condition \ref{p3.4.A}, 
\item at each crossing point the colors of the opposite 
edges (which belong to the image of the same 1-stratum of $\GG$) 
coincide, 
\item for each 1-stratum its color and the colors of the two adjacent 
2-strata satisfy the following condition \ref{p10.1.B} (cf. 
\ref{p3.4.A}): \end{enumerate} 

\begin{prop}[Admissibility of 2-Palette Colors at a
1-Stratum]\label{p10.1.B} Let a 1-stratum be colored with weight $A$ 
and let $A_1$ and  $A_2$ be the weight components of the 
colors of the adjacent 2-strata. Let $\Ge_i=1$, if the orientation of 
the $i^{\text{th}}$ adjacent 2-stratum induces the orientation of the
1-stratum and $\Ge_i=-1$ otherwise. Then
\begin{equation}\label{p10.1.B1}
 A+\Ge_1 A_1+\Ge_2 A_2=\pm1.
\end{equation}
\end{prop}

The sign in the right hand side of \eqref{p10.1.B1} depends on the 
orientation of the 1-stratum. The orientation for which the sign is 
plus is said to {\em defined\/} by the $\cal G^2_P$-coloring. 

A coloring of $\GG$ defines a coloring of the 1-strata of the diagram. 
A coloring of the diagram extending this one is called a {\it face 
extension\/} of the coloring of $\GG$. Making a face extension, one 
should care only about the last of the three conditions from the 
definitions \ref{p10.1.A} or \ref{p10.1.B} of $\cal G^1_P$- or $\cal 
G^2_P$ colorings of a diagram. 

\subsection{Alexander Invariant of Colored Strata}\label{sec10.2}
Given a diagram $\cal G^c_P$-coloring, we associate a 
$\cal G^c_P$-colored framed generic graph to
each strata of the diagram.

For a 2-stratum colored with $u$ this is \ ${\fig{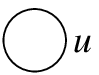}}$,\ 
i.~e. an unknot with zero framing colored with $u$:
$${\includegraphics{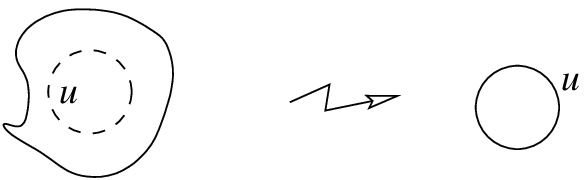}}$$

Recall that the Alexander invariant of \ ${\fig{circle-u.eps}}$\  
is $(t^{2}-t^{-2})^{-1}$ in the case $c=1$, $u=(t,T)$, and 
$(\im^{1+A}+\im^{1-A})^{-1}$ in the case $c=2$, $u=A$.

\begin{figure}[htp]
\centerline{\includegraphics{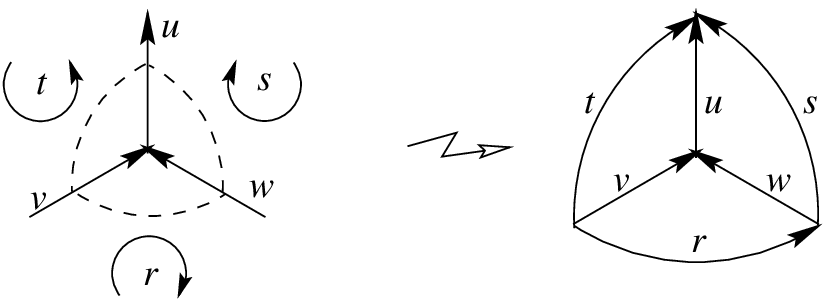}}
\caption{}
\label{ogettetr}
\end{figure}

For a 3-valent vertex whose adjacent edges are colored with $u$, $v$ 
and $w$ and adjacent 2-strata colored with  $r$, $s$ and $t$ as shown 
in Figure \ref{ogettetr}, this is the 1-skeleton of a tetrahedron 
embedded in a plane with the planar framing. It is obtained by 
selecting small arcs of the 1-strata adjacent to the vertex and 
joining their free end points with three arcs in the plane. The 
orientations of the 1-strata included in the colors induce 
orientations on the former three arcs. The orientations of the 
adjacent 2-strata included into the colors of these 2-strata induce 
orientations on the latter three arcs as on the boundary of the parts 
of the 2-strata separated from the vertex by the arcs. See Figure 
\ref{ogettetr}. If $c=1$ and the initial 3-valent vertex is strong, 
the orientation at it is inherited by the corresponding vertex of the 
tetrahedron. At other strong vertices, the tetrahedron is equipped 
with the counter-clockwise orientation. 

\begin{figure}[htp]
\centerline{\includegraphics{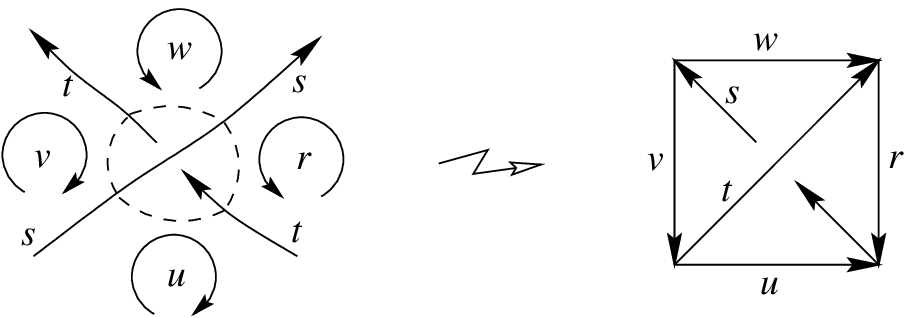}}
\caption{}
\label{ogetetr2}
\end{figure}

For a crossing point of strings colored with  $s$ and $t$ with the 
adjacent 2-strata colored with $u$, $v$, $w$ and $r$ as in Figure 
\ref{ogetetr2}, this is again the 1-skeleton of a tetrahedron. However 
this time it is not embedded in the plane of the picture, but 
presented by a diagram which looks like a square with diagonals and 
the framing is black-board. The diagram is obtained  from a regular 
neighborhood of the crossing point in the same way as in the case of a 
3-valent vertex: select small arcs of the 1-strata adjacent to the 
crossing point and join their free end points with four arcs 
surrounding the crossing point in the plane. As above, the 
orientations of the 1-strata induce the orientations on the former 
arcs, while the orientations of the adjacent 2-strata induce 
orientations on the former four arcs in the same way as above: as on 
the boundary of the parts of the 2-strata separated from the crossing 
point by the arcs. See Figure \ref{ogetetr2}. At strong vertices the 
graph is equipped with the counter-clockwise orientation. 

The Alexander invariant of this 1-skeleton of tetrahedron differs from the 
Alexander 
invariant of the 1-skeleton of tetrahedron corresponding to a 3-valent 
vertex, because 
of difference between their framings. The latter 
is isotopic to the former, but with non-black-board 
framing: 
$${\fig{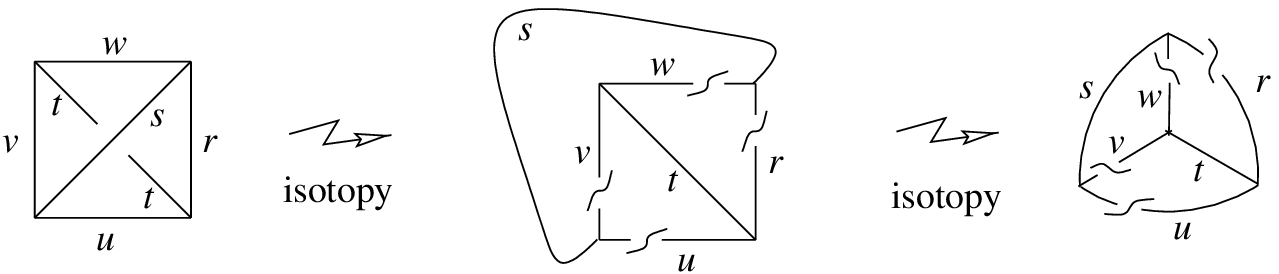}}$$ 
Thus here the Alexander invariant involves an 
additional factor. If $c=1$ the graph has strong vertices and they 
do not belong 
to the same diagonal of the original square diagram, the factor $-1$ also 
appears, since the isotopy turns one strong vertex upside down.

We have not associated any graph with 1-strata. This could be done in
a natural way, but would  not contribute to the state sum, which is
our main goal here. The reason is that the graph associated with a 
1-stratum should be a flat theta graph and the Alexander invariant of a 
flat theta graph is equal to 1.

The last ingredient that we have to introduce here is the {\it twist
factor.\/} It is associated to each  half-twist sign on a 1-stratum.
For $c=1$ and a stratum colored with multiplicity $t$ and weight $T$, 
this is $t^{-T}$, if the half-twist is positive, and $t^{T}$, if it is
negative. For $c=2$ and a stratum colored with weight $A$, 
this is $\im^{\frac{A^2-1}4}$, if the half-twist is positive, and 
$\im^{\frac{1-A^2}4}$, if it is negative.

\subsection{Spot Faces, Shadow Domain and Contour}\label{sec10.3} Let 
$\GG\subset\R^3$ be a framed $\cal G^c_P$-colored generic graph. 
Consider its diagram. Fix several faces of the diagram. They will be 
referred to as {\it spot faces.\/} Fix colors of the spot faces in 
such a way that Admissibility Conditions \ref{p10.1.A} are satisfied 
for any 1-stratum adjacent to two spot faces. 

Let $C$ be a smooth 1-dimensional closed submanifold of $\R^2$ 
which
is in a general position with respect to the diagram of $\GG$. This means
that $C$ contains no vertex of the diagram, is transversal to its
1-strata and does not pass over the signs
${\fig{plhalftw.eps}}$ and
${\fig{nghalftw.eps}}$ of half-twists. Let $\cal S$ be an
open subset of $\R^2$ with boundary $C$. Assume that each of the spot
faces intersects $\cal S$. The set $\cal S$ will be called {\it shadow
domain\/} and $C$ the {\it contour.\/}

The shadow domain is decomposed to its intersections with strata of 
the diagram of $\GG$. Some of the strata have fixed colors. These are 
all the 1-strata and the intersections of $\cal S$ with the spot 
faces. Denote by $Col_P(\cal S)$ the set of all the $\cal 
G^c_P$-colorings of 1- and 2-strata of $\cal S$ satisfying all the 
conditions of Section \ref{sec10.1} and extending this coloring of 
1-strata and spot faces. 

Contour $C$ and the part of the diagram of $\GG$ contained in the 
complement of $\cal S$  constitute a diagram of a generic graph. Denote 
this graph by $\GG_{C}$. Although its embedding in $\R^3$ is not 
specified, it  is well-defined up to ambient isotopy, since its diagram is 
given. The black-board framing of $C$ together with the framing signs on 
the part of the diagram of $\GG$ define a framing of $\GG_{C}$. 

The intersection points with the diagram divide $C$
into arcs, which are 1-strata of $\GG_{C}$. Each of these
1-strata lies on the boundary of one of the 2-strata of $\cal S$.  For
$c\in Col_P(\cal S)$ assign  to each arc of $C$ the color of
the adjacent 2-stratum of $\cal S$. The coloring of $\GG_{C}$
obtained in this way satisfied the admissibility condition. Denote the
skein class of this colored framed generic graph by $\GG_{C,c}$.

\subsection{State Sums}\label{sec10.4} Denote by $\Gt(\cal S)$ the
product of the twist factors of all the half-twist signs in $\cal S$.
For $x\in Col_P(\cal S)$ and a stratum $\GS$ of $\cal S$ denote by
$\UD^c(\GS,x)$ the Alexander invariant of the colored framed
generic graph associated  to $\GS$ for coloring~$x$. Denote
$$\Gt(S)\sum_{x\in Col_P(\cal S)}\UD^c(\GG_{C,x})
\prod_\GS\UD^c(\GS,x)^{\chi(\GS)}$$
by $Z(\cal S)$. Here the products are taken over all the strata $s$ of
$\cal S$ and $\chi$ is the Euler characteristic.

\begin{prop}[Moving Contour Theorem]\label{p10.4.A}
$Z(\cal S)$ does not change under an isotopy of $C$ in $\R^2$ during
which the spot faces remain intersecting the shadow domain and
$Col_P(\cal S)\ne\empt$. \end{prop}

In the simplest and most interesting cases by an isotopy of $C$
one can shrink the shadow domain to a subset of the spot faces. Then
the state sum is reduced to a single term and it is easy to relate it
to the Alexander invariant of the original graph $\GG$.

On the other hand, by isotopy of $C$ one can make the shadow domain
engulfing most of the diagram. Then graph $\GG_{C}$ may become pretty
standard and it will be easy to calculate the factors 
$\GD^c(\GG_{C,x})$.

Thus, Theorem \ref{p10.4.A} provides a way to express the Alexander 
invariant of an arbitrary colored framed generic graph through state sums 
involving the well-known Alexander invariants of standard unknotted 
graphs. The state sum is shaped by a diagram of the original graph. 

Theorem \ref{p10.4.A} is analogous to Kirillov-Reshetikhin theorem 
\cite{KR} 
for invariants based on the quantized $sl(2)$ (i.~e., generalizations of 
the Jones polynomial). There are though two important points where they 
differ. First, the terms of the $sl(2)$ state sum are more complicated 
functions of colors. For example, in our state sums the factors  
corresponding to edges are invisible (since the Alexander invariant of a 
theta graph is 1). The factors corresponding to other strata are also much 
simpler in our case. Second, the condition about non-vanishing of the set 
of colorings does not appear in the $sl(2)$ case, while here it is 
crucial.

\begin{proof}[Proof of \ref{p10.4.A}] Under a generic isotopy of $C$
the picture changes
topologically only when $C$ passes through vertices of the diagram of
$\GG$, through the signs of half-twists or becomes tangent to a branch of
the diagram. Therefore we have to check invariance of $Z(\cal S)$ only
under the moves shown in Figure \ref{f5}.
\begin{figure}[tbhp]
\centerline{\includegraphics{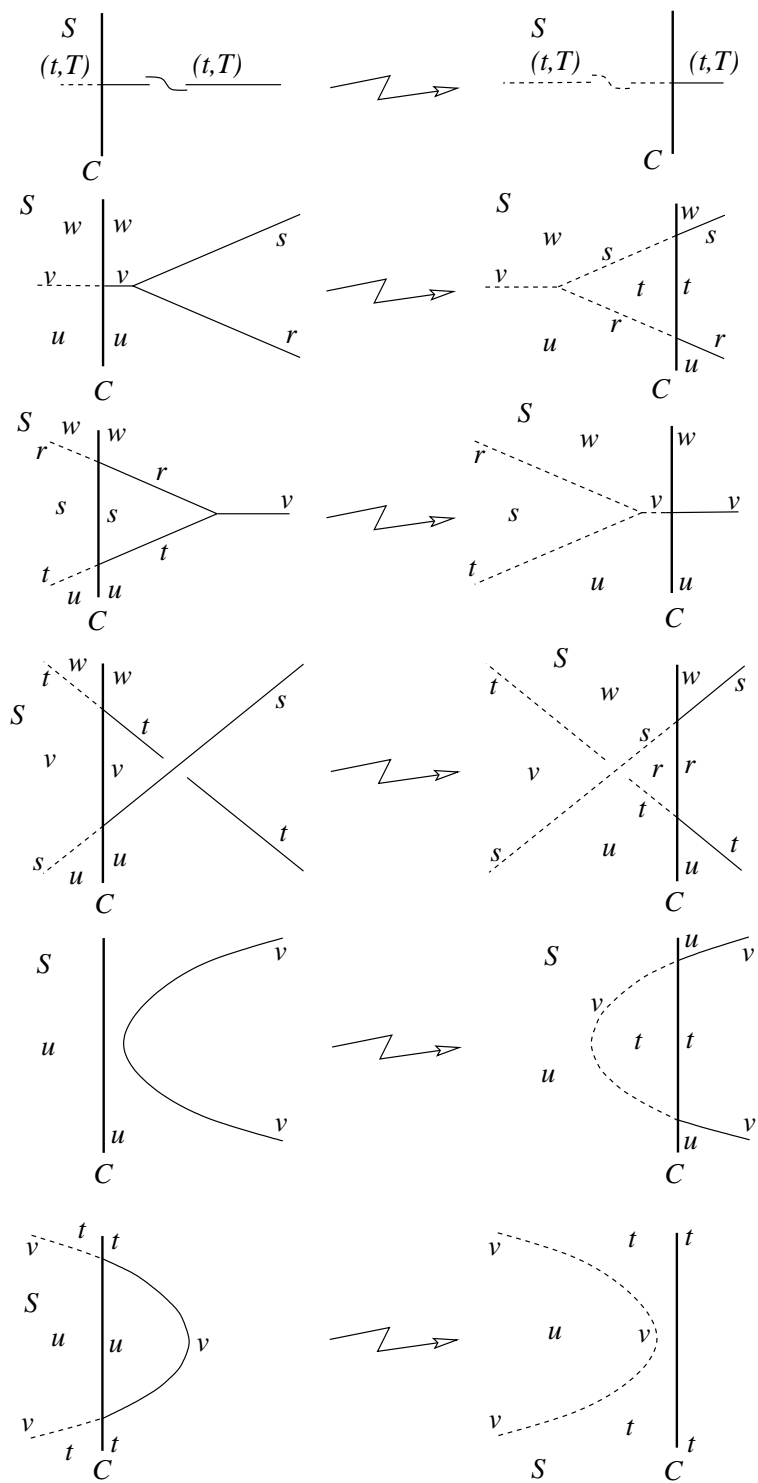}}
\caption{Moves of $\GG_{C,x}$ under generic isotopies of
$C$. The fragments of $\GG_{C,x}$ are shown in solid, while
pieces of the diagram of $\GG$ which are in the shadow domain are
shown by dashed arcs.}\label{f5} \end{figure}

Under the first move the set $Col_P(\cal S)$ of colorings does not
change (in the sense that there is a natural one-to-one correspondence
between the colorings before and after the move). Consider the case 
$c=1$. 
By \eqref{A-r-tw-skein}, $\UD^1(\GG_{C,x})$ for each $x\in Col_P(\cal
S)$ is divided by $t^{-T}$, where $t$ is the multiplicity and $T$ the
weight of the color of the arc involved in the move. On the other
hand, the total twist factor $\Gt(\cal S)$ of the shadow domain is
multiplied by $t^{-T}$, since an additional half-twist appeared on
it. Therefore $Z(\cal S)$ does not change under the first move.	
Similar arguments prove this in the case $c=2$.

To prove invariance with respect to the second move, decompose the
sets of colorings before and after the move to the classes of
colorings coinciding outside the fragment that changes. Such a class
before the move consists of a single coloring, say, $x$, while a
class after the move contains as many colorings as colors can be
put in the new triangle. Denote this color by $t$, and 
the corresponding coloring by $x_t$. Let the adjacent strata are 
colored as shown in Figure \ref{f5}. The Alexander invariants
$\UD^c(\GG_{C,x})$ for colorings in these classes before and after the 
move are related by \eqref{skein4} (with $\A^c$ replaced by $\UD^c$, 
see Section
\ref{sec9.3}). Applicability of \eqref{skein4} follows from the assumption
that $Col_P(\cal S)\ne\empt$ during the movement.  Rewrite this relation 
as follows:

$$\UD^c(\GG_{C,x})=\sum_t
\UD^c\( \fig{tetr12.eps}\ \)
\UD^c\(\fig{circlet.eps}\) \UD^c(\GG_{C,x_t}).$$

On the other hand,
$$\prod_\GS\UD^c(\GS,x_t)^{\chi(\GS)}=\prod_\GS\UD^c(\GS,x)^{\chi(\GS)} 
\UD^c\( \fig{tetr12.eps}\ \)\UD^c\(\fig{circlet.eps}\)$$
since in the shadow new strata appear, and these strata are:
\begin{enumerate}
\item A new trivalent vertex, which contributes
$\UD^c\( \fig{tetr12.eps}\ \)$.
\item New edges without half-twists. They do not contribute.
\item A new triangle colored with $t$. It contributes
$\UD^c\(\fig{circlet.eps}\)$.
\end{enumerate}
Therefore
\begin{multline}\UD^c(\GG_{C,x})\prod_\GS\UD^c(\GS,x)^{\chi(\GS)}\\=
\sum_t
 \UD^c\( \fig{tetr12.eps}\ \)\UD^c\(\fig{circlet.eps}\)
\UD^c(\GG_{C,x_t})\prod_\GS\UD^c(\GS,x)^{\chi(\GS)}\\=
\sum_t\UD^c(\GG_{C,x_t}) \prod_\GS\UD^c(\GS,x_t)^{\chi(\GS)}
\end{multline}

Invariance with respect to the third, fourth, fifth and sixth moves is 
proved
in the same way as invariance with respect to the second move, but using
instead of \eqref{skein4}
\begin{itemize}
\item \eqref{skein6} in the case of the third move,
\item \eqref{skein5} in the case of the
fourth move,
\item \eqref{skein3-g} and \eqref{skein3-g2} in the case of the
fifth move,  and
\item \eqref{skein2-g} and \eqref{skein2-g2} in the case of the
sixth move.\end{itemize} \end{proof}

\subsection{Face State Sums for Alexander Invariant}\label{sec10.5}
Moving Contour Theorem \ref{p10.4.A} does not give an explicit recipe
for calculating the Alexander invariant. Here we deduce 
explicit formulas for this. The formulas depend of choice of spot faces 
and movement of contour to which Moving Contour Theorem is applied.

Of course, the spot faces are needed for making the state sum finite. 
In our case they play also a role which does not emerge in the $sl(2)$ 
case: to make the state sum non-zero. For example, if one chooses a 
single spot face (as it is usual in $sl(2)$ case), $Z(\cal S)$ 
vanishes. Indeed, by Moving Contour Theorem, the contour can be made 
disjoint from the projection of the graph under consideration. This 
makes $\GG_C$ splittable and annihilates $\UD^c(\GG_C)$ by 
\ref{sec6.5}. 

The next possibility for the choice of spot faces is to choose two faces 
next to each other. This works as follows.

\begin{prop}[Face State Sum at Arc]\label{p10.5.A}Let
$\GG\subset\R^3$ be a $\cal G^c_P$-colored framed closed gene\-ric
graph, let $e$ be a 1-stratum of its diagram and $s$ be the
color of $e$. Let $f_1$, $f_2$ be the 2-strata of the
diagram adjacent to $e$. Color them with colors $u$ and $v$ such that 
the admissibility condition on $e$ is satisfied. Denote by 
$Col_P(u,v)$ the set of all the
admissible $\cal G^c_P$-colorings of the diagram  of $\GG$ which are
face extensions of the coloring of $\GG$ and have the colors described
above on $f_1$ and $f_2$. Let $\Gt$ be the product of the twist factors
of all the half-twist signs on the diagram. 
If $Col_P(u,v)\ne\empt$ then
\begin{equation}\label{FSSatArc}
\UD^c(\GG)=\frac{\Gt\UD^c\(\ \fig{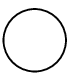}s\ \)}
{\UD^c\(\fig{circle.eps}u\ \)\UD^c\(\fig{circle.eps}v\ \)}
  \sum_{x \in Col_P(u,v)}\prod_{\GS}\UD^c(\GS,x)^{\chi(\overline\GS)}
\end{equation}
where $\overline \GS=\GS$ if $\GS$ is bounded in $\R^2$, otherwise
$\overline\GS$ is obtained by adding to $\GS$ the point of infinity in
the one-point compactificaion $S^2$ of $\R^2$.
\end{prop}

\begin{proof} First, take for a contour $C$ a small circle
bounding in $\R^2$ a disk $\cal S$, which intersects the diagram only
in an arc contained in $e$. The disk is the shadow domain. Then
$Col_P(u,v)$ consists of a single coloring, say $x$,
$$\UD^c(\GG_{C,x})=\frac{\UD^c(\GG)}{\UD^c\(\fig{circle.eps}s\ \)}
$$
and    $\prod_\GS\UD^c(\GS,x)^{\chi(\GS)} 
\UD^c\(\fig{circle.eps}u\ \)\UD^c\(\fig{circle.eps}v\ \)$. Therefore 
$$Z(\cal S)= 
\UD^c(\GG)\frac{\UD^c\(\fig{circle.eps}u\ \)\UD^c\(\fig{circle.eps}v\ \)}
{\UD^c\(\ \fig{circle.eps}s\ \)}.$$

Now expand $C$ till the moment when $\cal S$ engulfs the whole
diagram. 
If $Col_P(u,v)\ne\empt$ then we can apply Theorem \ref{p10.4.A},
and conclude that $Z(\cal S)$ does not change. At the final moment,
$\UD^c(\GG_{C,x})=\UD^c\(\fig{circle.eps}t\ \)$ where $t$ is the 
color of $C$. Therefore,
$$Z(\cal S)=\Gt\sum_{x\in
Col_P(u,v)}\UD^c\(\fig{circle.eps}t\ \)
\prod_\GS\UD^c(\GS,x)^{\chi(\GS)}=\Gt\sum_{x\in
Col_P(u,v)}\prod_\GS\UD^c(\GS,x)^{\chi(\overline \GS)}$$
Comparing with the initial value of $Z(\cal S)$, we get the desired
result. \end{proof}

\section{Appendix 1. Quantum $gl(1|1)$ and Its
Irreducible Representations}\label{sLast}

\subsection{Superalgebra $gl(1|1)$}\label{sL.1} Recall that in algebra
prefix {\it super-\/} means that the object under consideration is
$\Z_2$-graded or, at least, is related to a $\Z_2$-graded object.

For example, a vector superspace of (super-)dimension $(p|q)$ is a
$\Z_2$-graded vector space $V=V_0\oplus V_1$ whose $0$-summand $V_0$ is of
dimension $p$ and $1$-summand $V_1$ is $q$-dimensional. Denote by
$\C^{p|q}$ the superspace $\C^p\oplus\C^q$. The $0$-summand of a
superspace is called its {\it bosonic part\/} and vectors belonging to it
are called {\it bosons\/}, the $1$-summand is called {\it fermionic
part\/}, its elements are called {\it fermions.\/}

Endomorphisms of a superspace $\C^{p|q}$ comprise a Lie superalgebra,
which is denoted by $gl(p|q)$. It consists of
$(p|q)\times(p|q)$-matrices  $M=\left(\begin{matrix}A
&B\\C&D\end{matrix}\right)$. This is a {\it super}-algebra, its bosonic
part consists of  $\left(\begin{matrix}A &0\\0&D\end{matrix}\right)$,
the fermionic part consists of
$\left(\begin{matrix}0 &B\\C&0\end{matrix}\right)$. On bosons and
fermions the Lie (super)-brackets
are defined as supercommutator $[X,Y]=XY-(-1)^{\deg X\deg Y}YX$, where
$\deg X$  is 0 if $X$ is a boson and 1 if $X$ is a fermion. This
extends linearly to the whole $gl(p|q)$

The {\it supertrace\/} $\str M$ of $M=\left(\begin{matrix}A
&B\\C&D\end{matrix}\right)$ is  defined to be $\tr A-\tr D$.
This is a natural bilinear non-degenerate
form $$gl(p|q)\times gl(p|q)\to\C:(x,y)\mapsto\str(xy).$$

The Lie super-algebra $gl(1|1)$ is generated by the following 4 elements:
\begin{equation}\label{klpsi}
E=\left(\begin{matrix}1&0\\0&1\end{matrix}\right),\quad
G=\(\begin{matrix}0&0\\0&1\end{matrix}\),\quad
Y=\(\begin{matrix}0&1\\0&0 \end{matrix}\),\quad
X=\(\begin{matrix}0&0\\1&0 \end{matrix}\).\end{equation}

subject to the following relations
\begin{equation}
\begin{gathered}\label{relgl11}
\{X,Y\}(=XY+YX)=E,\\
X^2=Y^2=0,\\
[G,X]=X,\quad[G,Y]=-Y,\\
[E,G]=[E,X]=[E,Y]=0.\end{gathered}\end{equation}
In fact, these relations define $gl(1|1)$.

\subsection{Quantization}\label{sL.2} The universal enveloping algebra
$Ugl(1|1)$ admits a deformation (see P.~P.~Kulish \cite{Ku}) 
with parameter $h$ resulting the quantum superalgebra
$U_qgl(1|1)$. Only the first of relations \eqref{relgl11} changes: its
right hand side $E$ is replaced by ``quantum $E$'', that is
$\dfrac{q^E-q^{-E}}{q-q^{-1}}$ with $q=e^h$. To make sense of this 
non-algebraic expression,
one adjoins to the algebra formal power series in $h$.

Denote
$$q^{E}=e^{hE}=\sum_{n=0}^\infty\dfrac{h^nE^n}{n!}$$
by $K$. Since $E$ is central, $K$ is central, too.

Algebra $A$ of formal power series in $h$ generated by $X$, $Y$, 
$E$ and $G$ and satisfying relations
\begin{equation}\label{relglq11}\begin{gathered}
\{X,Y\}=\frac{K-K^{-1}}{q-q^{-1}},\\
X^2=Y^2=0,\\
[E,X] =0,\quad [E,Y]=0,\quad [E,G]=0,\\
[G,X]=X,\quad [G,Y]=-Y.\end{gathered}\end{equation}
can be equipped with a co-product
\begin{equation}\label{co-prod}\begin{aligned}\GD:A\to& A\otimes A:\\
&\GD(E)=1\otimes E+E\otimes1\\
&\GD(G)=1\otimes G+G\otimes1\\
&\GD(X)=X\otimes K^{-1}+1\otimes X,\\
&\GD(Y)=Y\otimes 1+K\otimes Y
\end{aligned}\end{equation}
co-unit
\begin{equation}\label{co-unit}\Ge:A\to\C:\qquad
\Ge(X)=\Ge(Y)=\Ge(E)=\Ge(G)=0,\end{equation}
and antipode
\begin{equation}\label{antipode}\begin{aligned}s:A\to&
A:\\
&s(E)=-E,\quad s(G)=-G,\quad s(X)=-XK,\quad
s(Y)=-YK^{-1}\end{aligned}\end{equation}
so that it turns to a Hopf (super-)algebra. It is denoted by
$U_qgl(1|1)$ and called {\it quantum supergroup $U_qgl(1|1)$.\/}

One can easily check that
\begin{equation}\begin{aligned}
\GD(K)&=K\otimes K,\quad \GD(K^{-1})=K^{-1}\otimes K^{-1},\\
\quad s(K)&=K^{-1},\quad \Ge(H)=1.
\end{aligned}\end{equation}

\subsection{$R$-Matrix}\label{sL.3}
Recall that a Hopf algebra $A$ is said to be {\it quasitriangular\/} if it
is equipped with a {\it universal $R$-matrix \/} $R$ which is an
invertible element of $A^{\otimes2}=A\otimes A$ satisfying the
following three conditions:
\begin{equation}\label{R1}
P\circ\GD(a)=R\GD(a)R^{-1}\text{ for any }a\in A
\end{equation}
\begin{equation}\label{R2}
(\GD\otimes\id_A)(R)=R_{13}R_{23}
\end{equation}
\begin{equation}\label{R3}
(\id_A\otimes\GD)(R)=R_{13}R_{12}
\end{equation}
Here $P$ is the permutation homomorphism $a\otimes b\mapsto (-1)^{\deg
a\deg b}b\otimes
a:A\otimes A\to A\otimes A$, $R_{12}=R\otimes1\in A^{\otimes3}$,
$R_{13}=(\id\otimes P)(R_{12})$, $R_{23}=1\otimes R$.

Quantum supergroup $U_qgl(1|1)$ is known to be quasitriangular with the
universal $R$-matrix
\begin{equation}\label{R-matrix}
R=\(1+(q-q^{-1})(X\otimes Y)(K\otimes K^{-1}) \)q^{-E\otimes 
G-G\otimes E} \end{equation}

Each quasitriangular  Hopf superalgebra has a unique element
$$u=\sum_i (-1)^{\deg \Ga_i\deg \Gb_i}s(\Gb_i)\Ga_i,$$
where $R=\sum_i \Ga_i\otimes \Gb_i$. In
$U_qgl(1|1))$
$$u=q^{2EG}\(1+(q-q^{-1})KYX\).$$
Then
$$s(u)=q^{2EG}\(1-(q-q^{-1})K^{-1}XY\),$$
cf. \cite{RT1} and \cite{RS2}.

Recall that a universal twist  of a quasitriangular Hopf
superalgebra $A$ is an element $v\in A$ such that $v^2=us(u)$.
In $U_qgl(1|1)$ this means
$$v^2=us(u)=q^{4EG}\(1+(q-q^{-1})(KYX-K^{-1}XY)\),$$
cf. \cite{RT1} and \cite{RS2}. One can easily check that this equation
in $U_qgl(1|1)$ has a solution
$$v=q^{2EG}\(K-(q-q^{-1})XY\).$$

\subsection{$(1|1)$-Dimensional Irreps}\label{sL.4} By an
$U_qgl(1|1)$-module  (of dimension $(p|q)$) we will mean a vector
superspace $V$ of dimension $(p|q)$ equipped with a homomorphism of
$U_qgl(1|1)$ to the superalgebra of endomorphisms of $V$ (for instance
$gl(p|q)$).

This homomorphism is supposed to respect the $\Z_2$-grading of the
superstructure. In particular, the images of bosonic $E$ and $G$ map
bosons to bosons and fermions to fermions, while the images of fermions
$X$ and $Y$ map bosons to fermions and fermions to bosons.

Recall that a module is said to be {\it cyclic\/} if it is generated as a
module by one vector, and {\it irreducible\/} if it does not contain a
proper submodule.

There are two families of irreducible $(1|1)$-dimensional
$U_ggl(1|1)$-modules. In each of them the modules are parametrized by
two parameters, which are denoted by $j$ and $J$ and run over
$\C\sminus\{\pi n\sqrt{-1}/2h : n\in\Z\}$ and $\C$ respectively.

The module of the first family corresponding to $(j,J)$ is described by
the following formulas:
\begin{equation}\label{(t,w)_{b}}\begin{aligned}(j,J)_+:&U_qgl(1|1)\to
gl(1|1):\\ &E\mapsto
\(\begin{matrix}2j&0\\0&2j\end{matrix}\),\qquad
G\mapsto
\(\begin{matrix}\frac{J-1}2&0\\0&\frac{J+1}2  \end{matrix}\),\\
&X\mapsto\(\begin{matrix}0&0\\\frac{q^{2j}-q^{-2j}}{q-q^{-1}}&0 
\end{matrix}\),\qquad Y\mapsto\(\begin{matrix}0&1\\0&0 \end{matrix}\)
\end{aligned}
\end{equation}

Recall that $\dfrac{q^n-q^{-n}}{q-q^{-1}}$ is denoted by $[n]_q$ and
called {\it quantum $n$\/}. Hence
$(j,J)_+:X\mapsto\(\begin{matrix}0&0\\ {[}2j{]}_q&0\end{matrix}\)$.

The standard basis vectors in $\C^{1|1}$ are denoted by $e_0$ (boson) and 
$e_1$ (fermion). On these vectors the generators of
$U_qgl(1|1)$ act via $(t,w)_+$ as follows:
\begin{equation}\label{E_in_(t,w)_+}
Ee_0=2j e_0,\quad Ee_1=2j e_1
\end{equation}
\begin{equation}\label{G_in_(t,w)_+}
Ge_0=\tfrac{J-1}2 e_0,\quad Ge_1=\tfrac{J+1}2e_1
\end{equation}
\begin{equation}\label{X_in_(t,w)_+}
 Xe_0=[2j]_qe_1,\quad Xe_1=0
\end{equation}
\begin{equation}\label{Y_in_(t,w)_+}
 Ye_0=0,\quad Ye_1=e_0
\end{equation}

One can easily find that $u\in U_qgl(1|1)$ (see Section \ref{sL.2}) acts
in $(j,J)_+$ as multiplication by $q^{2j(J+1)}$  and $v$, as
multiplication by $q^{2jJ}$. Hence, $v^{-1}u$ acts as multiplication by
$q^{2j}$.

Another family of irreducible $U_qgl(1|1)$-modules, denoted by $(j,J)_-$,
is obtained from $(j,J)_+$ by switching the bosonic and fermionic 
subspaces. We change also the signs of the parameters in 
anticipation of \ref{L.5.A}.

\begin{equation}\label{(t,w)_{c}}\begin{aligned}(j,J)_-:&U_qgl(1|1)\to
gl(1|1):\\ &E\mapsto\(\begin{matrix}-2j&0\\0&-2j\end{matrix}\),\qquad
G\mapsto
\(\begin{matrix}\frac{-J+1}2&0\\0&\frac{-J-1}2\end{matrix}\),\\
& X\mapsto\(\begin{matrix}0&-[2j]_q\\0&0 \end{matrix}\) ,\qquad
Y\mapsto\(\begin{matrix}0&0\\1&0 \end{matrix}\)
\end{aligned}
\end{equation}

Operator $u\in U_qgl(1|1)$ acts in $(j,J)_-$,  as multiplication by
$q^{2j(J-1)}$ and $v$, as multiplication by $q^{2jJ}$, so $v^{-1}u$
acts as multiplication by $q^{-2j}$.

\begin{prop}[Classification of $(1|1)$-Dimensional Irreps] An
irredu\-c\-ible $U_qgl(1|1)$-module of dimension $(1|1)$ 
 is isomorphic to
a module of one of the families \eqref{(t,w)_{b}}, \eqref{(t,w)_{c}}. No
module of one of these families is equivalent to a module of another one.
\qed \end{prop}

\subsection{Duality}\label{sL.5} Recall that for a
finite-dimensional $A$-module $V$, where $A$ is a Hopf superalgebra
over a field $k$, the dual space $V\spcheck=\Hom_{k}(V,k)$ is equipped
with a structure of $A$-module by the formula
$(Lc)(v)=c((-1)^{\deg L\deg c}s(L)(v))$, where $L\in A$, $c\in
V\spcheck$, $v\in V$ and $s$ is the antipode of $A$. This representation
is said to be {\it dual\/} to $V$.

\begin{prop}[Dual $(1|1)$-Dimensional Modules]\label{L.5.A} Modules
$(j,J)_+$ and $(j,J)_-$ are dual to each other. There is a canonical
isomorphism $(j,J)_+\spcheck\to(j,J)_-$, which is described by the
formulas: $e^0\mapsto e_0$, $e^1\mapsto q^{2j}e_1$, where $e^i$ is 
the generator
of $(j,J)_+\spcheck$ from the basis super-dual to the standard basis
$e_0$, $e_1$ of $(j,J)_+$, it is defined by $e^i(e_j)=(-1)^{ij}\Gd^i_j$.
\end{prop}

The duality allows one to define various pairings and co-pairings.
Cf. \cite{RT1}.
First, there is an obvious pairing
\begin{equation}\label{d-right-2dim}
\ld:(j,J)_-\times(j,J)_+\to\C:\quad
(x,y)\mapsto x(y).
\end{equation}
Here we identify $(j,J)_-$ with $((j,J)_+)\spcheck$ by the isomorphism of
\ref{L.5.A}. This pairing acts as follows:
$$\ld:(e_a,e_b)\mapsto(-1)^{a}q^{-2ja}\Gd_{ab}.$$

Another pairing is obtained from this one by a  ``quantum
super-transposition'':
\begin{equation}\label{d-left-2dim}
\rd:(j,J)_+\times(j,J)_-\to\C:\quad
(x,y)\mapsto (-1)^{\deg x\deg y}y(v^{-1}ux).
\end{equation}
The factor $(-1)^{\deg x\deg y}$ appears here because of transposition of
$x$ and $y$ (this is the contribution of super-mathematics). On the basis
vectors the pairing acts as follows:
$$\rd:(e_a,e_b)\mapsto q^{2j(1-a)}\Gd_{ab}.$$

There is an obvious (super-) co-pairing
\begin{equation}\label{b-right-2dim}
\lb:\C\to (j,J)_+\otimes(j,J)_-:\quad 1\mapsto e_0\otimes e_0-e_1\otimes
q^{2j}e_1, \end{equation}
and ``quantum super-transposed'' co-pairing
\begin{multline}\label{b-left-2dim}
\rb:\C\to (j,J)_-\otimes(j,J)_+:\quad 1\mapsto
\sum_{a=0}^1(-1)^ae_a\otimes (u^{-1}ve_a)\\= q^{-2j}e_0\otimes
e_0+e_1\otimes e_1.
\end{multline}

\subsection{Tensor Product $(i,I)_+\otimes(j,J)_+$}\label{sL.6}

\begin{prop}[Lemma. Action of Generators]\label{L.6.A}Generators $E$,
$G$, $X$, $Y$ of $U_qgl(1|1)$ act in $(i,I)_+\otimes(j,J)_+$ as 
follows: \begin{equation}\label{act in (i,I)b otimes (j,J)b}
\begin{matrix}
{\begin{aligned} E:\quad
&e_0\otimes e_0\mapsto 2(i+j)(e_0\otimes e_0)\\
&e_1\otimes e_1\mapsto 2(i+j)(e_1\otimes e_1)\\
&e_0\otimes e_1\mapsto 2(i+j)(e_0\otimes e_1)\\
&e_1\otimes e_0\mapsto 2(i+j)(e_1\otimes e_0)
\end{aligned}}&
{\begin{aligned}G:\quad
&e_0\otimes e_0\mapsto \tfrac{I+J-2}2(e_0\otimes e_0)\\
&e_1\otimes e_1\mapsto \tfrac{I+J+2}2(e_1\otimes e_1)\\
&e_0\otimes e_1\mapsto \tfrac{I+J}2(e_0\otimes e_1)\\
&e_1\otimes e_0\mapsto \tfrac{I+J}2(e_1\otimes e_0)
\end{aligned}}\end{matrix}\end{equation}
$$
{\begin{aligned}X:\quad
&e_0\otimes e_0\mapsto [2j]_q(e_0\otimes
e_1)+q^{-2j}[2i]_q(e_1\otimes e_0)\\
&e_1\otimes e_1\mapsto 0\\
&e_0\otimes e_1\mapsto q^{-2j}[2i]_q(e_1\otimes e_1)\\
&e_1\otimes e_0\mapsto -[2j]_q(e_1\otimes e_1)\\
Y:\quad&e_0\otimes e_0\mapsto 0\\
&e_1\otimes e_1\mapsto (e_0\otimes e_1)-q^{2i}e_1\otimes e_0\\
&e_0\otimes e_1\mapsto q^{2i}(e_0\otimes e_0)\\
&e_1\otimes e_0\mapsto e_0\otimes e_0
\end{aligned}}
$$
\end{prop}
\begin{proof}  Recall that the action of $U_qgl(1|1)$ in the tensor
product is defined via the coproduct $\GD:U_qgl(1|1)\to U_qgl(1|1)\otimes
U_qgl(1|1)$. One can easily check the formulas above using
\eqref{co-prod}.
\end{proof}

\begin{prop}[Decomposition in Generic Case]\label{L.6.B} Let $i$, $j$
and $k=i+j$ be complex
numbers which are not of the form $\pi n\sqrt{-1}/2h$ with $n\in\Z$. Let
$I$, $J$  be arbitrary complex numbers.
Then $(i,I)_+\otimes(j,J)_+$ is isomorphic to
$$(i+j,I+J-1)_+\oplus(-i-j,-I-J-1)_-.$$
There is an isomorphism
$$(i+j,I+J-1)_+\oplus(-i-j,-I-J-1)_-\to(i,I)_+\otimes(j,J)_+
$$
with matrix
$$\(\begin{matrix}1&0&0&0\\
 0&0&1&0\\
0&\dfrac{[2j]_q}{[2i+2j]_q}&0&1\\
0&q^{-2j}\dfrac{[2i]_q}{[2i+2j]_q}&0&-q^{2i}
                 \end{matrix}\)$$
with respect to the bases $\{(e_0,0),(e_1,0),(0,e_0),(0,e_1)\}$ of
$(i+j,I+J-1)_+\oplus(-i-j,-I-J-1)_-$ and $\{e_0\otimes
e_0,e_1\otimes e_1,e_0\otimes e_1, e_1\otimes e_0\}$ of
$(i,I)_+\otimes(j,J)_+$.
\end{prop}

\begin{proof}
Consider the submodule generated by $e_0\otimes e_0$. Since
$E(e_0\otimes e_0)=2(i+j)(e_0\otimes e_0)$ and $i+j\not\equiv0\mod(\pi
\sqrt{-1}/2h)$, this is a cyclic $(1|1)$-dimensional representation. Since
$Y(e_0\otimes e_0)=0$ and $e_0\otimes e_0$  is a boson, this is a
vector representation. Since $G(e_0\otimes e_0)=(\frac{I+J}2-1)(e_0\otimes 
e_0)$, this is $(i+j,I+J-1)_+$. Choose for its basis
$e_0\otimes e_0$ and
$X(1/[2i+2j]_q)(e_0\otimes e_0)=
[2j]_q/[2i+2j]_q(e_0\otimes e_1)+
q^{-2j}[2i]_q/[2i+2j]_q(e_1\otimes e_0)$.

Consider now the submodule generated by $e_1\otimes e_1$. Since
$E(e_1\otimes e_1)=2(i+j)(e_1\otimes e_1)$ and $i+j\not\equiv0\mod(\pi
\sqrt{-1}/2h)$, this is a cyclic $(1|1)$-dimensional representation. Since
$X(e_1\otimes e_1)=0$ and $e_1\otimes e_1$  is a boson, this is a
covector representation. Since $G(e_1\otimes 
e_1)=(\frac{I+J}2+1)(e_1\otimes 
e_1)$ and $(j,J)_-(G)(e_0)=(\frac{-J+1}2)e_0$ (see \eqref{(t,w)_{c}}),
this is $(-i-j,-I-J-1)_-$. Take $e_1\otimes e_1$ and $Y(e_1\otimes
e_1)= e_0\otimes e_1-q^{2i}(e_1\otimes e_0)$ for its basis.
\end{proof}

\begin{prop}\label{L.6.C}
In representation $(i,I)_+\otimes (j,J)_+$ with the basis $e_0\otimes
e_0$, $e_0\otimes e_1$, $e_1\otimes e_0$,  $e_1\otimes e_1$, the
$R$-matrix is
$$ q^{-iJ-j I}\(\begin{matrix}q^{i+j}&0&0&0 \\
                                  0&q^{j-i}&0&0 \\
                       0&q^{i-j}(q^{2i}-q^{-2i})&q^{i-j}&0 \\
                                  0&0&0&q^{-i-j} \end{matrix}\)$$
\end{prop}

\begin{proof}
Recall (see \eqref{R-matrix}) that the universal $R$-matrix is
$$R=\(1+(q-q^{-1})(X\otimes Y)(K\otimes K^{-1}) \)
q^{-E\otimes G-G\otimes E}.$$

It is easy to see that $q^{-E\otimes G-G\otimes E}$ is
$$q^{-iJ-j I} \(\begin{matrix} q^{i+j}&0&0&0 \\
                                  0&q^{j-i}&0&0 \\
                                  0&0&q^{i-j}&0 \\
                                  0&0&0&q^{-i-j}
   \end{matrix}\)
$$
Since $Xe_1=0$ and $Ye_0=0$, out of all the basis vectors only
$e_0\otimes e_1$ is not annihilated by $X\otimes Y$.
Further, $q^{-E\otimes G-G\otimes E}(e_0\otimes e_1)=q^{-iJ-j
I}q^{j-i}(e_0\otimes e_1)$ and
$$(X\otimes Y)(e_0\otimes e_1)=[2i]_q(e_1\otimes e_0),$$
so
$$(q-q^{-1})(X\otimes Y)(K\otimes K^{-1})q^{-E\otimes G-G\otimes 
E}(e_0\otimes
e_1)= q^{-iJ-j I}q^{i-j}(q^{2i}-q^{-2i})(e_1\otimes e_0).$$
All together this gives the desired formula.
\end{proof}

$R$-matrix will be used always composed with the operator of
super-transposition $$P:(i,I)_+\otimes (j,J)_+\to(j,J)_+\otimes
(i,I)_+ $$ with matrix $$ \(\begin{matrix}1&0&0&0\\ 0&0&1&0\\ 0&1&0&0\\
0&0&0&-1 \end{matrix}\) $$ (Here $-1$ at the right bottom corner is due to
super-rule: the transposition exchanges two fermions $e_1$ and $e_1$ in
$e_1\otimes e_1$ and fermions skew commute.)

\begin{multline} 
P\circ R=q^{-iJ-j I}\(\begin{matrix}q^{i+j}&0&0&0 \\
                          0&q^{i-j}(q^{2i}-q^{-2i})&q^{i-j}&0 \\
                                  0&q^{j-i}&0&0 \\
                                  0&0&0&-q^{-i-j}\end{matrix}\)
				\\=q^{-iJ-jI-i-j}
				\(\begin{matrix}q^{2i+2j}&0&0&0 
                            \\ 0&q^{4i}-1&q^{2i}&0 \\
                                  0&q^{2j}&0&0 \\
                                  0&0&0&-1\end{matrix}\)
\end{multline}

Similarly one can calculate the corresponding objects related to tensor 
products $(i,I)_-\otimes(j,J)_-$, $(i,I)_+\otimes(j,J)_-$ and 
$(i,I)_-\otimes(j,J)_+$.

\subsection{A $q$-Less Subalgebra}\label{sL.10} Consider
super-subring $U^1$ of $U_qgl(1|1)$ generated by bosons 
$K=q^{E}$, $K^{-1}$, $F=2G$ and fermions $Z=(q-q^{-1})X$, 
$Y$, see \ref{sL.2}. This is a superalgebra over $\Z$ defined by relations: 
\begin{equation}\label{relA}\begin{gathered}
\{Z,Y\}=K-K^{-1},\\
Z^2=Y^2=0,\\
[K,Z] =0,\quad [K,Y]=0\\
[F,Z]=2Z\quad [F,Y]=-2Y.\end{gathered}
\end{equation} 
The coproduct of $U_qgl(1|1)$ induces in $U^1$ a coproduct
\begin{equation}\label{co-prodA}\begin{aligned}\GD:U^1\to& U^1\otimes U^1:\\
&\GD(K)=K\otimes K\\
&\GD(K^{-1})=K^{-1}\otimes K^{-1}\\
&\GD(F)=1\otimes F+F\otimes1,\\
&\GD(Z)=Z\otimes K^{-1}+1\otimes Z,\\
&\GD(Y)=Y\otimes 1+K\otimes Y.
\end{aligned}\end{equation}
$U^1$ has a co-unit
\begin{equation}\label{co-unitA}\Ge:U^1\to\Z:\qquad
\Ge(F)=\Ge(Y)=\Ge(Z)=0, \quad\Ge(K)=\Ge(K^{-1})=1\end{equation}
and antipode
\begin{equation}\label{antipodeA}\begin{aligned}s:U^1\to&
U^1:\\
&s(K)=K^{-1},\quad s(F)=-F,\quad s(Z)=-ZK,\quad
s(Y)=-YK^{-1}\end{aligned}\end{equation}
so that $U^1$ is a Hopf (super-)algebra over $\Z$.

\subsection{Irreducible Representations of $U^1$}\label{sL.11} Let 
$P=(B,M,W,M\times W\to M)$ be a 1-palette (see Section \ref{sec2.8}). 
That is let $B$ be a commutative ring , $M$ a subgroup of the 
multiplicative group of $B$, and $W$ be a subgroup of the additive 
group of $B$ equipped with a (bilinear) pairing $M\times W\to 
M:(m,w)\mapsto m^w$ such that $1\in W$ and $m^1=m$ for each $m\in M$. 

For $t\in M$, $t^4\ne1$ and $T\in W$ consider the action of 
$U^1$ on $B^{(1|1)}$ defined by 
\begin{equation}\label{(t,T)+}\begin{aligned}
 K\mapsto 
\left (\begin{matrix}t^2&0 \\ 0&t^2\end{matrix}\right ), \quad
&F\mapsto\(\begin{matrix}T-1&0\\0&T+1\end{matrix}\)  \\
 Y\mapsto\(\begin{matrix}0&1\\0&0 \end{matrix}\), \quad
Z&\mapsto\(\begin{matrix}0&0\\ t^2-t^{-2}&0 \end{matrix}\).
\end{aligned}
\end{equation}
This is an irreducible representation. Denote it by $U(t,T)_+$. The dual 
representation is an action defined as follows in the same space:
\begin{equation}\label{(t,T)-}\begin{aligned}
K\mapsto\(\begin{matrix}t^{-2}&0\\0&t^{-2}\end{matrix}\),&\quad
F\mapsto\(\begin{matrix}-T+1&0\\0&-T-1\end{matrix}\)\\           
Y\mapsto\(\begin{matrix}0&0\\1&0 \end{matrix}\),\quad
Z&\mapsto\(\begin{matrix}0&t^{-2}-t^2\\0&0 \end{matrix}\).
\end{aligned}
\end{equation}
It is denoted by $U(t,T)_{-}$.

\begin{prop}[Duality]\label{sL.11.A} Modules
$(t,T)_+$ and $(t,T)_-$ are dual to each other. There is a canonical
isomorphism $(t,T)_+\spcheck\to(t,T)_-$, which is described by the
formulas: $e^0\mapsto e_0$, $e^1\mapsto t^{2}e_1$, where $e^i$ is 
the generator
of $(t,T)_+\spcheck$ from the basis super-dual to the standard basis
$e_0$, $e_1$ of $(t,T)_+$, it is defined by $e^i(e_j)=(-1)^{ij}\Gd^i_j$.
\qed\end{prop}

\begin{prop}[Decomposition of $(u,U)_+\otimes_B(v,V)_+$]\label{L.11.B} Let
$u$, $v$ and $t=uv$ belong to $M\sminus\{t\in M\mid t^4=1\}$ and
$t^2-t^{-2}$ be invertible in $B$. Let
$U,V\in W$.
Then $(u,U)_+\otimes(v,V)_+$ is isomorphic to
$(uv,U+V-1)_+\oplus(u^{-1}v^{-1},-U-V-1)_-$.
There is an isomorphism
$$(t,U+V-1)_+\oplus(t^{-1},-U-V-1)_-\to(u,U)_+\otimes(v,V)_+
$$
with matrix
$$\(\begin{matrix}t^2-t^{-2}&0&0&0\\
 0&0&1&0\\
0&v^2-v^{-2}&0&1\\
0&v^{-2}(u^2-u^{-2})&0&-u^{2}
                 \end{matrix}\)$$
with respect to the bases $\{(e_0,0),(e_1,0),(0,e_0),(0,e_1)\}$ of the
$B$-module
$$(t,U+V-1)_+\oplus(t^{-1},-U-V-1)_-$$ and $\{e_0\otimes
e_0,e_1\otimes e_1,e_0\otimes e_1, e_1\otimes e_0\}$ of
$(u,U)_+\otimes_B(v,V)_+$.\qed
\end{prop}

\begin{prop}[Decomposition of $(u,U)_-\otimes_B(v,V)_-$]\label{L.11.C}
Let $u$, $v$ and
$t=uv$ belong to $M\sminus\{t\in M\mid t^4=1\}$  and
$t^2-t^{-2}$ be invertible in $B$. Let $U,V\in W$.
Then $(u,U)_-\otimes_B(v,V)_-$ is isomorphic to
$(t^{-1},-U-V-1)_+\oplus(t,U+V-1)_-$. There is an isomorphism
$$(t^{-1},-U-V-1)_+\oplus(t,U+V-1)_-\to(u,U)_-\otimes_B(v,V)_-
$$
with matrix
$$\(\begin{matrix} 0&0&1&0\\
                  t^2-t^{-2}&0&0&0\\
                  0&v^{2}(u^2-u^{-2})&0&u^{-2}\\
                  0&-v^2+v^{-2}&0&1
                  \end{matrix}\)$$
with respect to the bases $\{(e_0,0),(e_1,0),(0,e_0),(0,e_1)\}$ of the
$B$-module
$$(t^{-1},-U-V-1)_+\oplus(t,U+V-1)_-$$ and $\{e_0\otimes
e_0,e_1\otimes e_1,e_0\otimes e_1, e_1\otimes e_0\}$ of
$(u,U)_-\otimes_B(v,V)_-$.
\end{prop}

\begin{prop}[Decomposition of $(u,U)_+\otimes_B(v,V)_-$]\label{L.11.D}
Let $u$, $v$ and
$t=uv^{-1}$ belong to $M\sminus\{t\in M\mid t^4=1\}$  and
$t^2-t^{-2}$ be invertible in $B$. Let $U,V\in W$.
Then $(u,U)_+\otimes_B(v,V)_-$ is isomorphic to
$(t,U-V+1)_+\oplus(t^{-1},V-U+1)_-$. There is an isomorphism
$$(t,U-V+1)_+\oplus(t^{-1},V-U+1)_-\to(u,U)_+\otimes_B(v,V)_-$$
with matrix
$$\(\begin{matrix} 1&0&v^2-v^{-2}&0\\
-u^{2}&0&-v^{2}(u^2-u^{-2})&0\\
                  0&0&0&t^2-t^{-2}\\
                  0&1&0&0
\end{matrix}\)$$
with respect to the bases $\{(e_0,0),(e_1,0),(0,e_0),(0,e_1)\}$ of the
$B$-module
$$(t^{-1},V-U-1)_+\oplus(t,U-V-1)_-$$ and $\{e_0\otimes
e_0,e_1\otimes e_1,e_0\otimes e_1, e_1\otimes e_0\}$ of
$(u,U)_+\otimes(v,V)_-$. \qed
\end{prop}

\section{Appendix 2. Representations of Quantum $sl(2)$ 
at $\sqrt{-1}$}\label{A2}

\subsection{Algebra $sl(2)$}\label{A2.1}
Recall that  Lie algebra $sl(2)$ is generated by the following 3 
elements: \begin{equation}\label{A2.1.1}
H=\left(\begin{matrix}1&0\\0&-1\end{matrix}\right),\quad
X^+=\(\begin{matrix}0&1\\0&0 \end{matrix}\),\quad
X^-=\(\begin{matrix}0&0\\1&0 \end{matrix}\).\end{equation}

subject to the following relations
\begin{equation}\label{relsl2}
[X^+,X^-]=H,\quad
[H,X^+]=2X^+,\quad[H,X^-]=-2X^-.
\end{equation}
In fact, these relations define $sl(2)$.

\subsection{Quantization}\label{A2.2} The universal enveloping algebra 
$Usl(2)$ admits a deformation resulting the quantum superalgebra 
$U_qsl(2)$. Only the first of relations \eqref{relsl2} changes: its right 
hand side $H$ is replaced by ``quantum $H$'', that is 
$\dfrac{q^H-q^{-H}}{q-q^{-1}}$  where $q=e^h$. To make sense of this 
non-algebraic expression, one adjoins formal power series in $h$ to the 
algebra. Denote
$$q^{H}=e^{hH}=\sum_{n=0}^\infty\dfrac{h^nH^n}{n!}$$
by $K$. Thus the deformed system of relations is:

\begin{equation}\label{relslq2}\begin{gathered}
K=q^H,\quad
\[X^+,X^-\]=\frac{K-K^{-1}}{q-q^{-1}},\\
[H,X^+]=2X^+,\quad[H,X^-]=-2X^-.\end{gathered}\end{equation}
The algebra $A$ over the field $\C[[h]]$ of formal power series in $h$ 
which consists of formal power series in $h$ with coefficients generated 
by $H$, $X^+$ and $X^-$ subject to \eqref{relslq2} can be equipped 
with a co-product
\begin{equation}\label{sl2co-prod}\begin{aligned}\GD:A\to& A\otimes A:\\
&\GD(H)=H\otimes 1+1\otimes H\\
&\GD(X^+)=X^+\otimes 1+K\otimes X^+,\\
&\GD(X^-)=X^-\otimes K^{-1}+1\otimes X^-
\end{aligned}\end{equation}
co-unit
\begin{equation}\label{sl2co-unit}\Ge:A\to\C[[h]]:\quad
\Ge(K)=\Ge(K^{-1})=1,\quad\Ge(H)=\Ge(X^+)=\Ge(X^-)=0,\end{equation}
and antipode
\begin{equation}\label{sl2antipode}\begin{aligned}s:A\to&
A:\\
&s(H)=-H,\quad  
s(X^+)=-K^{-1}X^+,\quad
s(X^-)=-X^-K\end{aligned}\end{equation}
so that it turns to a Hopf algebra. It is denoted by
$U_qsl(2)$.

\subsection{Simplification at square root of $-1$}\label{A2.3} Consider 
$U_qsl(2)$ at a special value $q=\sqrt{-1}$ corresponding to 
$h=\frac{\pi\sqrt{-1}}2$. To abbreviate formulas, but keep letter $i$ 
available for other purposes, let us denote $\sqrt{-1}$ by $\im$ and 
$\exp(a\frac{\pi\sqrt{-1}}2)$ by $\im^a$.

In case $q=\im$ the generators may be changed to simplify relations. 
Indeed, let 
$$X=X^-,\quad Y=2X^+\sqrt{-1}. $$
Then relations 
\begin{equation}\label{relsli2}\begin{gathered}
K=\im^H,\quad\[Y,X\]=K-K^{-1},\\
[H,X]=-2X,\quad [H,Y]=2Y, \text{ and hence }
KX=-XK,\quad KY=-YK.\end{gathered}\end{equation}
define the algebra $U_{\im}{sl}(2)$. This is a Hopf 
algebra with co-product defined by 
\begin{equation}\label{co-prod-i}\begin{aligned}\GD:
U_{\im}sl(2)\to&U_{\im}sl(2)\otimes
U_{\im}sl(2):
\\
&\GD(H)=H\otimes 1+1\otimes H\\ 
&\GD(K)=K\otimes K\\
&\GD(K^{-1})=K^{-1}\otimes K^{-1}\\
&\GD(X)=X\otimes K^{-1}+1\otimes X,\\
&\GD(Y)=Y\otimes 1+K\otimes Y,
\end{aligned}\end{equation}
co-unit defined by
\begin{equation}\label{co-unit-i}\begin{aligned}\Ge:
U_{\im}sl(2)\to&
\C[[h]]:\\
&\Ge(K)=\Ge(K^{-1})=1,\quad\Ge(H)=\Ge(X)=\Ge(X)=0,\end{aligned}
\end{equation} and antipode
\begin{equation}\label{antipode-i}\begin{aligned}s:
U_{\im}sl(2)\to&
U_{\im}sl(2):\\
&s(H)=-H,\\
&s(K)=K^{-1},\quad s(K^{-1})=K,\\
&s(X)=-XK\quad(=KX),\\
&s(Y)=-K^{-1}Y\quad(=YK^{-1}).
\end{aligned}\end{equation}

\subsection{Family of Modules of Dimension 2}\label{A2.4}
The relations \eqref{relsli2} are satisfied by matrices
\begin{equation}\begin{aligned} H=
\(\begin{matrix}a+1 &0\\0&a-1 \end{matrix}\),&\quad K=
\(\begin{matrix}\im^{a+1} &0\\0&-\im^{a+1} \end{matrix}\),\\
X=\(\begin{matrix}0&0\\ \im^{a+1}-\im^{-a-1}&0 
\end{matrix}\),&\quad
Y=\(\begin{matrix}0&1\\
0&0 \end{matrix}\),\end{aligned}\end{equation}
where $a\in\C$.

The space of this representation is a two-dimensional vector space over 
the quotient field of the ring of formal power series $\C[[h]]$. Recall 
that this field consists of formal Laurent series, i.~e. $\C[[h]][h^{-1}]$. 

Denote this representation of $U_{\im}sl(2)$ by $I(a)$.
The standard basis vectors in $\(\C[[h]][h^{-1}]\)^{2}$ are denoted by 
$e_0$, $e_1$. On these vectors the generators of
$U_{\im}sl(2)$ act via $I(a)$ as follows:
\begin{equation}\label{H_in_a}
He_0=(a+1)e_0,\quad He_1=(a-1)e_1
\end{equation}
\begin{equation}\label{K_in_a}
Ke_0=\im^{a+1} e_0,\quad Ke_1=-\im^{a+1} e_1
\end{equation}
\begin{equation}\label{X_in_a}
 Xe_0=(\im^{a+1}-\im^{-a-1})e_1,\quad Xe_1=0
\end{equation}
\begin{equation}\label{Y_in_a}
  Ye_0=0,\quad Ye_1=e_0
\end{equation}

This representation is irreducible, unless $\im^{a+1}-\im^{-a-1}= 0$, i.~e. 
$a$ is an odd integer.

\subsection{Duality}\label{A2.6} Recall that for any finite-dimensional 
representation $V$ of a Hopf algebra $A$ over a field $k$ the dual 
space $V\spcheck=\Hom_{k}(V,k)$ is equipped with a structure of $A$-module 
by the formula $(Tc)(v)=c(s(T)(v))$, where $T\in A$, 
$c\in V\spcheck$, $v\in V$ and $s$ is the antipode of $A$. This 
representation is said to be {\it dual\/} to $V$. 

\begin{prop}[Duality Between $I(a)$ and 
$I(-a)$]\label{A2.6.A} 
Modules $I(a)$ and $I(-a)$ are dual to 
each other. There is an isomorphism $$f:I(a)\spcheck\to I(-a),$$ 
which is described by the formulas:
$$f:e^0\mapsto e_1,\quad f:e^1\mapsto\im^{a-1}e_0,$$ where 
$e^i$ is 
the generator of $I(a)\spcheck$ from the basis dual to the standard basis 
$e_0$, $e_1$ of $I(a)$, it is defined by 
$e^i(e_j)=\Gd^i_j$. \qed\end{prop}

\begin{rem}\label{A2.6.B} The dual module $I(a)\spcheck$ can be 
continuously deformed to $I(a)$ in the family of modules $I(x)$: the 
complex numbers $-a$ and $a$ can be connected by a path in 
$\C\sminus(2\Z+1)$. 
\end{rem}

\subsection{$R$-Matrix}\label{A2.7} Quantum group $U_qsl(2)$ is known to 
be quasitriangular with the universal $R$-matrix
\begin{equation}\label{sl2R-matrix}
R=q^{H\otimes 
H/2}\sum_{n=0}^\infty\frac{q^{n(n-1)/2}}{[n]_q!}
\((q-q^{-1})X^-\otimes X^+\)^n
\end{equation}

In all the representations of $U_qsl(2)$ considered in this paper, 
$(X^+)^2=(X^-)^2=0$ and $q=\im$. Factorization by these relations reduces 
$R$-matrix to

\begin{equation}\label{R-matrix at i}
R=\im^{H\otimes H/2}\(1+X\otimes Y\).
\end{equation}

\begin{prop}\label{A2.7.A}
In representation $I(a)\otimes I(b)$ with the basis $e_0\otimes
e_0$, $e_0\otimes e_1$, $e_1\otimes e_0$,  $e_1\otimes e_1$, the
$R$-matrix is
$$ \im^{\frac{ab+1}2}\(\begin{matrix}\im^{\frac{a+b}2}&0&0&0 \\
                   0&\im^{\frac{-a+b}2+1}&0&0 \\                 
0&\im^{\frac{a-b}2}(\im^{a}+\im^{-a})&\im^{\frac{a-b}2-1}&0 
\\ 0&0&0&\im^{\frac{-a-b}2}\end{matrix}\)$$
\qed \end{prop}

$R$-matrix will be used always composed with the operator of
transposition $$P:I(a)\otimes I(b)\to I(b)\otimes I(a)$$ with matrix 
$$ \(\begin{matrix}1&0&0&0\\ 0&0&1&0\\ 0&1&0&0\\
0&0&0&1 \end{matrix}\) $$ 
The composition is
$$ P\circ R=\im^{\frac{ab+1}2}\(\begin{matrix}\im^{\frac{a+b}2}&0&0&0 \\
 0&\im^{\frac{a-b}2}(\im^{a}+\im^{-a})&\im^{\frac{a-b}2-1}&0 
\\                    0&\im^{\frac{-a+b}2+1}&0&0 \\ 
 0&0&0&\im^{\frac{-a-b}2}\end{matrix}\)
$$

\subsection{Tensor Product $I(a)\otimes I(b)$}\label{A2.8}

\begin{prop}[Action of Generators]\label{A2.8.A}Let $a$, $b$ be complex 
numbers. The generators 
$H$, $K$, $X$ and $Y$ act in $I(a)\otimes I(b)$ as follows:
\begin{equation}\label{act in I(a) otimes I(b)} 
\begin{aligned} H:\quad
&e_0\otimes e_0\mapsto (a+b+2)(e_0\otimes e_0)\\
&e_1\otimes e_1\mapsto (a+b-2)(e_1\otimes e_1)\\
&e_0\otimes e_1\mapsto (a+b)(e_0\otimes e_1)\\
&e_1\otimes e_0\mapsto (a+b)(e_1\otimes e_0)\\
 K:\quad
&e_0\otimes e_0\mapsto -\im^{a+b}(e_0\otimes e_0)\\
&e_1\otimes e_1\mapsto -\im^{a+b}(e_1\otimes e_1)\\
&e_0\otimes e_1\mapsto \im^{a+b} (e_0\otimes e_1)\\
&e_1\otimes e_0\mapsto \im^{a+b}(e_1\otimes e_0)\\
X:\quad&e_0\otimes e_0\mapsto \im(\im^{b}+\im^{-b})(e_0\otimes 
e_1)+(\im^{a}+\im^{-a})\im^{-b}(e_1\otimes e_0)\\
&e_1\otimes e_1\mapsto 0\\
&e_0\otimes e_1\mapsto (-\im^{a}-\im^{-a})\im^{-b}(e_1\otimes e_1)\\
&e_1\otimes e_0\mapsto \im(\im^{b}+\im^{-b})(e_1\otimes e_1)\\
Y:\quad
&e_0\otimes e_0\mapsto 0\\
&e_1\otimes e_1\mapsto e_0\otimes e_1-\im^{a+1}e_1\otimes e_0\\
&e_0\otimes e_1\mapsto \im^{a+1}(e_0\otimes e_0)\\
&e_1\otimes e_0\mapsto e_0\otimes e_0\\
\end{aligned}
\end{equation}
\end{prop}
\begin{proof}  Recall that the action of $U_{\im}sl(2)$ in the tensor
product is defined via the coproduct $\GD:U_{\im}sl(2)\to 
U_{\im}sl(2)\otimes  $. One can easily check the 
formulas above using \eqref{co-prod-i}.
\end{proof}
\begin{prop}[Decomposition in Generic Case]\label{A2.8.B} Let $a$ and $b$ 
be complex numbers such that neither $a$, nor $b$, nor $a+b+1$ is an 
odd integer. 
Then $I(a)\otimes I(b)$ is isomorphic to $I(a+b+1)\oplus I(a+b-1)$.
There is an isomorphism
$$I(a+b+1)\oplus I(a+b-1)\to I(a)\otimes I(b)
$$
with matrix
$$\(\begin{matrix}{\im^{-a-b}-\im^{a+b}}&0&0&0\\
 0&0&0&1\\
0&\im(\im^{b}+\im^{-b})&1&0\\
0&\im^{-b}(\im^{a}+\im^{-a})&\im^{a-1}&0
                 \end{matrix}\)$$
with respect to the bases 
$$\{(e_0,0),(e_1,0),(0,e_0),(0,e_1)\} \text{ of } I(a+b+1)\oplus 
I(a+b-1)$$ and 
$$\{e_0\otimes e_0,e_1\otimes e_1,e_0\otimes e_1, e_1\otimes e_0\} \text{ 
of } I(a)\otimes I(b).$$ \qed \end{prop}

\end{document}